\documentclass{report}

\usepackage{amsfonts,amsmath,eufrak}
\usepackage{pstricks,pst-node,cite}
\usepackage{amscd}
\usepackage{latexsym, amsmath, amssymb}

\newcommand{\Inv}{\mbox{\rm Inv}}

\newtheorem{thm}{Theorem}[chapter]
\newtheorem{prop}[thm]{Proposition}
\newtheorem{cor}[thm]{Corollary}

\newtheorem{lem}[thm]{Lemma}

\newenvironment{proof}{\noindent {\bf Proof: }}{\QED\medskip}
\def\QED{{\hspace*{\fill}{$\Box$}\quad}
    \vskip 0pt plus20pt}

\ifx\blackandwhite\undefined
  \psset{linecolor=blue}\newrgbcolor{darkred}{.8 0 0}\newrgbcolor{emgreen}{0 .8
  0}\newrgbcolor{pup}{.8 0 .8}\newrgbcolor{xxx}{.8 .8 .8}
\else
  \psset{linecolor=black}\newgray{darkred}{.7}\newgray{emgreen}{0.3}\newgray{pup}{.5}\newgray{xxx}{.5}
\fi

\psset{linewidth=.4pt,dash=3pt 3pt,doublesep=.05,dotsize=1pt 5}
\SpecialCoor

\newcommand{\vv}{\pspicture[.4](-.6,-.5)(.6,.5)
\psbezier(.5;45)(.25;45)(.25;315)(.5;315)
\psbezier(.5;225)(.25;225)(.25;135)(.5;135)
\endpspicture}
\newcommand{\hh}{\pspicture[.4](-.6,-.5)(.6,.5)
\psbezier(.5;45)(.25;45)(.25;135)(.5;135)
\psbezier(.5;225)(.25;225)(.25;315)(.5;315)
\endpspicture}
\newcommand{\btwohvert}{\pspicture[.4](-.6,-.7)(.6,.7)
\qline(0, .25)(.35, .6)\qline(0, .25)(-.35, .6)
\qline(0,-.25)(.35,-.6)\qline(0,-.25)(-.35,-.6)
\psline[doubleline=true](0,-.25)(0,.25)
\endpspicture}
\newcommand{\btwohhoriz}{\pspicture[.4](-.8,-.5)(.8,.5)
\qline( .25,0)( .6,.35)\qline( .25,0)( .6,-.35)
\qline(-.25,0)(-.6,.35)\qline(-.25,0)(-.6,-.35)
\psline[doubleline=true](-.25,0)(.25,0)
\endpspicture}

\newcommand{\doubleloop}{\pspicture[.4](-.6,-.5)(.6,.5)
\pscircle[doubleline=true](0,0){.4}
\endpspicture}

\newcommand{\middlearrow}{\lput{:U}{\pspicture(0,0)(0,0)
\psline[arrows=->,arrowscale=1.5](2.2pt,0)(2.3pt,0)\endpspicture}}

\begin{document}
\pagenumbering{roman}

\begin{center}
{\Large \bf Graphical Calculus on Representations of Quantum Lie Algebras}

\vskip 0.2in

\centerline{By}

\vskip 0.1in

\centerline{\large Dongseok KIM}

\centerline{\large B. S. (Kyungpook National University, Korea) 1990}

\centerline{\large M. S. (Kyungpook National University, Korea) 1992}

\centerline{\large M. A. (University of Texas at Austin) 1998}

\vskip 0.3in

\centerline{\large DISSERTATION}

\vskip 0.1in

\centerline{\large Submitted in partial satisfaction of the
requirements for the degree of}

\vskip 0.1in

\centerline{\large DOCTOR OF PHILOSOPHY}

\vskip 0.1in

\centerline{\large in}

\vskip 0.1in

\centerline{\large MATHEMATICS}

\vskip 3.2in

\centerline{\large UNIVERSITY OF CALIFORNIA}

\vskip 0.1in

\centerline{\large  DAVIS}

\vskip 0.1in

\centerline{\large March 2003}
\end{center}

\newpage
\large
\tableofcontents


\newpage
\large
{\Large \bf ACKNOWLEDGEMENTS} \\

I would like to thank my adviser, Greg Kuperberg. He has been
exceptionally patient to show me how to look Mathematics in many
different prospects. He has shared many ideas with me, as well as
guided me to the subject of representation theory. None of these
could be done without his supervision. I would like to thank
Professors Dimtry Fuchs, Joel Hass, Albert Schwartz and William
Thurston for teaching me. I would like to thank Professor Khovanov
for taking time to teach me many things that I should know and
helping me to untie some of my research. I would like to thank my
friends in the math department for helping me get through the
period of exhaustive study.

I would like to thank my family. They have always supported my
education. I especially thank my wife Youngmi, for standing
through our stay in Unites States. I can not forget my parent who
has been incredibly supportive for all. My children, Richard and
Hannah have been the biggest motivation for me to stay
in focus.


\newpage
\begin{center}
\underline{\bf \Large Abstract}
\end{center}

\medskip
\large

The main theme of this thesis is the representation theory of
quantum Lie algebras. We develop graphical calculation methods.
Jones-Wenzl projectors for $\mathcal{U}_q(\mathfrak{sl}(2,\mathbb{C}))$
are very powerful tools to find not only
invariants of links but also invariants of 3-manifolds.
We find single clasp expansions of generalized Jones-Wenzl projectors for
simple Lie algebras of rank $2$. Trihedron coefficients
of the representation theory for $\mathcal{U}_q(\mathfrak{sl}(2,\mathbb{C}))$
has significant meaning and it is called $3j$ symbols.
Using single clasp expansions for $\mathcal{U}_q(\mathfrak{sl}(3,\mathbb{C}))$, we find
some trihedron coefficients of the representation theory of
$\mathcal{U}_q(\mathfrak{sl}(3,\mathbb{C}))$. We study representation theory for
$\mathcal{U}_q(\mathfrak{sl}(4,\mathbb{C}))$.
We conjecture a complete set of relations for
$\mathcal{U}_q(\mathfrak{sl}(4,\mathbb{C}))$.


\newpage
\pagestyle{myheadings} \pagenumbering{arabic} \markright{  \rm
\normalsize CHAPTER 1. \hspace{0.5cm}
 Introduction }
\large
\chapter{Introduction}
\thispagestyle{myheadings}


There has been big progress in the theory bridging Lie algebras
and low-dimensional topology. These developments are based on
quantum groups, braided categories and new invariants of knots,
links and $3$-manifolds. After the discovery of the Jones
polynomial~\cite{Jones:poly}~\cite{Jones:braid}, Reshetikhin and
Turaev~\cite{RT:1}~\cite{RT:2} showed that braided categories
derived from quantum groups provide a natural generalization of
the Jones polynomial.

One of the developments is that a category of tangles with skein
relations leads to a braided category. If we decorate each
component of a tangle by a module over a simple Lie algebra, the
category becomes a ribbon category. Then we can get an invariant
of links, and sometimes 3-manifolds, from a functor constructed
in~\cite{Turaev:quantum}. To develop this theory further, we would
like to generalize the Jones-Wenzl projectors in the
Temperley-Lieb algebra to the quantization of other simple Lie
algebras. The $n$-th Temperley-Lieb algebra is realized as the
algebra of intertwining operators of the
$\mathcal{U}_q(\mathfrak{sl}(2,\mathbb{C}))$-module $V_1^{\otimes n}$, where
$V_1$ is the two-dimensional irreducible representation of
$\mathcal{U}_q(\mathfrak{sl}(2,\mathbb{C}))$. For each $n$, the algebra $T_n$
has an idempotent $f_n$ such that $f_nx=xf_n=\epsilon(x)f_n$ for
all $x\in T_n$ and $f_nf_n=f_n$, where $\epsilon$ is an
augmentation. These idempotents were first discovered by V.
Jones~\cite{Jones:subfactor} and H. Wenzl ~\cite{Wenzl:Proj}, and
they found a recursive formula:

$$f_{n}=f_{n-1}+\frac{[n-1]}{[n]}f_{n-1}e_{n-1}f_{n-1}.$$

So they
are named {\it Jones-Wenzl idempotents(Projectors)}.
Kuperberg~\cite{Kuperberg:spiders} defines a generalization of the
Temperley-Lieb category to the three rank two Lie algebras
$\mathfrak{sl}(3,\mathbb{C})$, $\mathfrak{sp}(4)$ and $G_2$. These
generalizations are called combinatorial rank two spiders. Also he
has proved that Jones-Wenzl projectors exist for simple Lie
algebras of rank $2$ and he called them {\it clasps}.
We will study how they can be expanded inductively in
Chapter $2$.

The skein module theory allows not only links but also
graphs. The invariants of the two simplest nontrivial trivalent graphs, the trihedron and
tetrahedron, have significant meaning and they are
called $3j$ and $6j$ symbols. So we can naturally ask how to compute
trihedron coefficients for $\mathcal{U}_q(\mathfrak{sl}(3,\mathbb{C}))$ as
suggested in ~\cite{Kuperberg:spiders}. In Chapter $3$, we will
apply our clasp expansions to find some trihedron coefficients for
$\mathcal{U}_q({\mathfrak sl}(3,\mathbb{C}))$.

Kuperberg's generalization of the Temperley-Lieb algebra is a set
of generators and relations for each rank $2$ Lie
algebra~\cite{Kuperberg:spiders}. The generators are easy to find,
but to get a complete set of relations is a challenging problem.
In Chapter 4, we follow Kuperberg's method to find some relations. We
conjecture a complete set of relations of
$\mathcal{U}_q(\mathfrak{sl}(4,\mathbb{C}))$.

\section{Preliminaries}

For simple terms, we refer to~\cite{Humphreys:gtm}~\cite{Kassel:gtm} and  \cite{KRT:knot}.

Quantum integers are defined as

\begin{align*}
[n]&= \frac{q^{n/2}-q^{-n/2}}{q^{1/2}-q^{-1/2}}\\
[0]&= 1\\
[n]!&= [n][n-1]\ldots [2][1]\\
\left[ \begin{matrix}  n\\
k \end{matrix} \right] &= \frac{[n]!}{[k]![n-k]!}
\end{align*}

Let $\mathfrak{sl}(n,\mathbb{C})$ be the Lie algebra of complex $n\times
n$-matrices with trace zero. Let $E_{i,j}$ be the elementary matrix
whose entries are all zero except $1$ in the $(i,j)$-th entry. Let
$E_i=E_{i,i+1}, F_i=E_{i+1,i}$ and $H_i=E_{i,i}-E_{i+1,i+1}$ where
$1\le i\le n-1$, then they generate  $\mathfrak{sl}(n,\mathbb{C})$ with
relations:

\begin{align*}
&[H_i,H_j]=0 &\mathrm{for}\hskip .1cm i, j=1, 2, \ldots, n-1\\
&[H_i,E_j]=\alpha_j(H_i)E_j &\mathrm{for}\hskip .1cm 1 \le i, j\le n-1\\
&[H_i,F_j]=-\alpha_j(H_i)F_j &\mathrm{for}\hskip .1cm 1 \le i, j\le n-1\\
&[E_i,F_j]=\delta_{ij}H_i &\mathrm{for}\hskip .1cm 1 \le i, j\le n-1\\
&[E_i,E_j]=0 &\mathrm{if}\hskip .1cm |i-j|\ge 2\\
&[F_i,F_j]=0 &\mathrm{if}\hskip .1cm |i-j|\ge 2\\
&[E_i,[E_i,E_j]]=0 &\mathrm{if}\hskip .1cm |i-j|=1\\
&[F_i,[F_i,F_j]]=0 &\mathrm{if}\hskip .1cm |i-j|=1
\end{align*}

where $\alpha_i$ is a linear form defined by
$$\alpha_{j}(H_i)=\begin{cases}2,&\mathrm{if}
\hskip .3cm i=j\\-1,&\mathrm{if}\hskip .3cm|i-j|=1\\
0,&\mathrm{Otherwise}\end{cases}$$

The quantum group $\mathcal{U}_q(\mathfrak{sl}(n,\mathbb{C}))$ is an
associative algebra over $\mathbb{C}(q)$ with generators, $E_i$, $
F_i$, $K_{i}^{\pm}$ with $1\le i\le n-1$, and relations:

\begin{align*}
&K_iK_i^{-1}=1=K_i^{-1}K_i &\mathrm{for}\hskip .1cm i=1, 2, \ldots, n-1\\
&K_iE_j=q^{\alpha_j(H_i)}E_jK_i &\mathrm{for}
\hskip .1cm i, j=1, 2, \ldots, n-1\\
&K_iF_j=q^{-\alpha_j(H_i)}F_jK_i &\mathrm{for}
\hskip .1cm i, j=1, 2, \ldots, n-1\\
&[E_i,F_j]=\delta_{ij}\frac{K_i-K_i^{-1}}{q-q^{-1}}
&\mathrm{for}\hskip .1cm i, j=1, 2, \ldots, n-1\\
&[E_i,E_j]=0 &\mathrm{if}\hskip .1cm |i-j|\ge 2\\
&[F_i,F_j]=0 &\mathrm{if}\hskip .1cm |i-j|\ge 2\\
&E_i^2E_j-(q+q^{-1})E_iE_jE_i+E_jE_i^2=0 &\mathrm{if}\hskip .1cm |i-j|=1\\
&F_i^2F_j-(q+q^{-1})F_iF_jF_i+F_jF_i^2=0 &\mathrm{if}\hskip .1cm
|i-j|=1
\end{align*}

Let $\mathfrak{h}$ be the Lie subalgebra of $\mathfrak{sl}(n,\mathbb{C})$
generated by $H_i$ and let $\Lambda \in\mathfrak{h}$ be the
integral lattice of linear forms on $H_i$ where $n-1\ge i\ge 1$.
Let $\lambda \in \Lambda, \epsilon = (\epsilon_1, \ldots,
\epsilon_{n-1})$, then there is a unique universal highest weight
module, a {\it Verma Module}, with highest weight $(\lambda,
\epsilon)$. $M(\lambda, \epsilon)$ has a unique simple quotient
$L(\lambda, \epsilon)$ which is highest module with highest weight
$(\lambda, \epsilon)$. Then $L(\lambda, \epsilon)$ is finite
dimensional if and only if $\lambda$ is dominant weight. One can
see that as $\mathcal{U}_q(\mathfrak{sl}(n,\mathbb{C}))$ module

$$L(\lambda, \epsilon) \cong L(\lambda, 0)\otimes L(0,\epsilon).$$

So we can study $L(\lambda, 0)$ which is denoted by $L(\lambda)$.
Then there is a theorem which connects studies of
$\mathcal{U}_q(\mathfrak{sl}(n,\mathbb{C}))$ modules and $\mathfrak{sl}(n,\mathbb{C})$ modules.

\begin{thm} ~\cite{KRT:knot}
\item{i)} Any finite dimensional simple $\mathcal{U}_q(\mathfrak{sl}(n,\mathbb{C}))$
module is of the form
$L(\lambda, \epsilon)$ where $\lambda$ is dominant weight and
$\epsilon \in (\mathbb{Z}/2\mathbb{Z})^{n-1}$.
\item{ii)} The character $ch(L(\lambda))$ is given by the same
formula as the character of simple  $\mathfrak{sl}(n)$ module
parameterized by the same highest weight.
\item{iii)} The multiplicity of a simple module $L(\nu)$
in the decomposition of the tensor product $L(\lambda) \otimes
L(\mu)$ of two simple modules is the same as for the decomposition
of the corresponding $\mathfrak{sl}(n,\mathbb{C})$ module.
\label{quantumtensor}
\end{thm}

\newpage
\pagestyle{myheadings}
\chapter{Single Clasp Expansions for Rank 2 Lie Algebras}
\thispagestyle{myheadings} \markright{  \rm \normalsize CHAPTER 2.
\hspace{0.5cm}
  Single Clasp Expansions for Rank 2 Lie Algebras }

\section{Introduction}
Let $T_n$ be the $n$-th Temperley-Lieb algebra with generators,
$1, e_1,e_2,\ldots , e_{n-1}$, and relations:

\begin{align*}
e_i^2&=-(q^{\frac{1}{2}}+q^{-\frac{1}{2}})e_i\\
e_ie_j&=e_je_i \hskip 2cm \mathrm{if}\hskip .2cm |i-j|\ge 2\\
e_i&=e_ie_{i\pm 1}e_i
\end{align*}

For each $n$, the algebra $T_n$ has an idempotent $f_n$ such that
$f_nx=xf_n=\epsilon(x)f_n$ for all $x\in T_n$, where $\epsilon$ is
an augmentation. These idempotents were first discovered by V.
Jones~\cite{Jones:subfactor} and H. Wenzl~\cite{Wenzl:Proj}. They
found a recursive formula:

$$f_{n}=f_{n-1}+\frac{[n-1]}{[n]}f_{n-1}e_{n-1}f_{n-1}$$

as in the following figure where we use a red box to represent
$f_n$:

\begin{eqnarray}
\pspicture[.4](-.1,-.3)(1.3,1.3)
\rput[t](.5,-.1){$n$}\qline(.5,0)(.5,.4)
\psframe[linecolor=darkred](0,.4)(1,.6)
\qline(.5,.6)(.5,1)\rput[b](.5,1.1){$n$}
\endpspicture
= \pspicture[.4](-.3,-.3)(1.5,1.3)
\rput[t](.5,-0.1){$n-1$}\qline(.5,0)(.5,.4)
\psframe[linecolor=darkred](0,.4)(1,.6)
\qline(.5,.6)(.5,1)\rput[b](.5,1.1){$n-1$} \qline(1.3,0)(1.3,1)
\endpspicture
+ \frac{[n-1]}{[n]} \pspicture[.4](-.3,-.3)(1.55,2.3)
\rput[t](.5,-.1){$n-1$} \qline(.5,0)(.5,.4)
\psframe[linecolor=darkred](0,.4)(1,.6)
\qline(.25,.6)(.25,1.4)\rput[l](.35,1){$n-2$}
\psframe[linecolor=darkred](0,1.4)(1,1.6) \qline(.5,1.6)(.5,2)
\rput[b](.5,2.1){$n-1$}
\psarc(1,.6){.2}{0}{180}\qline(1.2,0)(1.2,.6)
\psarc(1,1.4){.2}{180}{0}\qline(1.2,1.4)(1.2,2)
\endpspicture
\label{a12exp1}
\end{eqnarray}

So they are named {\it Jones-Wenzl idempotents(projectors)}. We
will recall an algebraic definition of Jones-Wenzl projectors in
section $1$. We refer to~\cite{Kuperberg:spiders} for definitions,
notation and simple calculations. We provide single clasp
expansions of generalized Jones-Wenzl projectors for
$\mathcal{U}_q({\mathfrak sl}(3,\mathbb{C}))$ in section $2$. In
section $3$ we study single clasp expansions of generalized
Jones-Wenzl projectors for $\mathcal{U}_q({\mathfrak sp}(4))$.

\section{Single Clasp Expansion for $\mathcal{U}_q({\mathfrak sl}(2,\mathbb{C}))$}

First we recall another definition of Jones-Wenzl projectors and
single clasp expansions for $\mathcal{U}_q({\mathfrak
sl}(2,\mathbb{C}))$. Then we use it to find trihedral coefficients
for $\mathcal{U}_q({\mathfrak sl}(2,\mathbb{C}))$.

\subsection{Jones-Wenzl Projector for $\mathcal{U}_q({\mathfrak sl}(2,\mathbb{C}))$}
Let us give a precise definition~\cite{khvanov:canonical} of a
{\it clasp} for ${\mathcal U}(\mathfrak{sl}(2,\mathbb{C}))$. Let
$V_i$ be an irreducible representation of highest weight $i$. Then
$i_n:V_n\rightarrow V_1^{\otimes n}$ is defined by
$$i_n(v^m)=\left[ \begin{matrix} n\\
\frac{n-m}{2} \end{matrix} \right]^{-1}\sum_{s, |s|=m}
q^{||s||_{-}}v^{s_1}\otimes\ldots\otimes v^{s_n}$$

and $\pi_n: V_1^{\otimes n}\rightarrow V_n$ is defined by

$$\pi_n(v^{s_1}\otimes v^{s_2}\otimes\ldots
\otimes v^{s_n})=q^{-||s||_{+}}v^{|s|}$$

where $s=(s_1,s_2, \ldots, s_n)$, $s_i=\pm 1$, $|s|=\sum s_i$ and
$||s||_{+}=\sum_{i<j} \{s_i>s_j\}$, $||s||_{-}=\sum_{i>j}
\{s_i>s_j\}$ and $\{a>b\}=1$ if $a>b$, and $0$ otherwise.

Then the composition $i_n\circ \pi_n$ is called a {\it Jones-Wenzl
projector}, denoted by $p_n$. It has the following properties 1)
it is an idempotent 2) $p_ne_i=0=e_ip_n$ where $e_i$ is a U-turn
from the $i$-th to the $i+1$-th string as in the following
figures.

$$
\pspicture[.4](-.3,0)(2.5,1)
\rput[r](-0.1,.5){$n$}\qline(0,.5)(.4,.5)
\psframe[linecolor=darkred](.4,0)(.6,1) \qline(.6,.5)(1.4,.5)
\rput[b](1,.6){$n$} \psframe[linecolor=darkred](1.4,0)(1.6,1)
\qline(1.6,.5)(2,.5)\rput[l](2.1,.5){$n$}
\endpspicture
= \pspicture[.4](-.5,0)(1.3,1)
\rput[r](-0.1,.5){$n$}\qline(0,.5)(.4,.5)
\psframe[linecolor=darkred](.4,0)(.6,1)
\qline(.6,.5)(1,.5)\rput[l](1.1,.5){$n$}
\endpspicture
\hskip 1cm ,\hskip 1cm \pspicture[.4](-.3,0)(3.3,1)
\qline(0,.5)(.4,.5)\rput[r](-.1,.5){$n$}
\psframe[linecolor=darkred](.4,0)(.6,1)
\qline(.6,.833)(1,.833)\rput[l](1.1,.833){$k$}
\psarc(.6,.5){.167}{-90}{90}
\qline(.6,.167)(1,.167)\rput[l](1.1,.167){$n-k-2$}
\endpspicture = 0
\label{claspaxiom}
$$

We can generalize the second property as follows: if we attach a
web with a cut path with less weight, then it is zero. Then we can
axiomatize these two properties to define generalized Jones-Wenzl
projectors for any simple Lie algebra.
Kuperberg~\cite{Kuperberg:spiders} proved that  Jones-Wenzl
projectors exist for simple Lie algebra of rank $2$ and he called
{\it clasps} (sometimes they are called {\it magic weaving
elements} or {\it boxes}). Here, we will call them clasps. For
$\mathcal{U}_q({\mathfrak sl}(2,\mathbb{C}))$, it is known that we
can inductively expand it as in equation~\ref{a12exp1}. For
advanced calculations, the single clasp expansion in
equation~\ref{a1exppic} is very useful and has been used
in~\cite{khvanov:canonical} for some beautiful results. By
symmetry, there are four different positions for the single clasp
expansion depending on where the clasp of weight $n-1$ is located.
For equation~\ref{a1exppic}, the clasp is located at the southwest
corner, which will be considered the standard expansion,
otherwise, we will state the location of the clasp.

\begin{eqnarray}
\pspicture[.45](-.1,-.3)(2.3,2.3)
\rput[t](1,-.1){$n$}\qline(1,0)(1,.4)
\psframe[linecolor=darkred](0,.4)(2,.6)
\qline(1,.6)(1,2)\rput[b](1,2.1){$n$}
\endpspicture
= \sum_{i=1}^{n} a_i \pspicture[.45](-.4,-.3)(2.35,2.5)
\rput[t](1,-.1){$n-1$} \qline(1,0)(1,.4)
\psframe[linecolor=darkred](0,.4)(2,.6) \qline(.2,.6)(.2,2)
\qline(.4,.6)(.4,2) \qline(.6,.6)(.6,2) \rput[l](0.75,2.3){$i$}
\rput[l](2.15,2.3){$1$} \qline(2.2,2)(2.2,1.5)
\qline(2.2,1.5)(1.5,1) \qline(1.5,1)(1.5,.6) \qline(2,2)(2,1.5)
\qline(2,1.5)(1.3,1) \qline(1.3,1)(1.3,.6) \qline(1.8,2)(1.8,1.5)
\qline(1.8,1.5)(1.1,1) \qline(1.1,1)(1.1,.6)
\psarc(1.9,.6){.3}{0}{180}\qline(2.2,0)(2.2,.6)
\psarc(1,1.7){.2}{180}{0}\qline(1.2,1.7)(1.2,2)
\qline(.8,1.7)(.8,2)
\endpspicture
\label{a1exppic}
\end{eqnarray}

\begin{prop}
The coefficients in equation~\ref{a1exppic} are
$$a_i=\frac{[n+1-i]}{[n]}.$$
\end{prop}
\begin{proof}
By attaching a $U$ turn at consecutive strings to the top, we have
the following $n-1$ equations.

$$a_{n-1}-[2]a_{n}=0.$$
For $i=1, 2, \ldots, n-2$,

$$a_i-[2]a_{i+1}+a_{i+2}=0.$$

One can see that these equations are independent. By attaching the
clasp of weight $n$ to the bottom of every web in
equation~\ref{a1exppic}, we get $a_1=1$ by the properties of a
clasp. This process is called a {\it normalization}. Then we check
the answer in the proposition satisfies these equations. Since
these webs in equation~\ref{a1exppic} form a basis, these
coefficients are unique.
\end{proof}

\subsection{Applications of Single Clasp Expansions
for $\mathcal{U}_q({\mathfrak sl}(2,\mathbb{C}))$}

We can easily prove the following propositions using the single
clasp expansion of $\mathcal{U}_q({\mathfrak sl}(2,\mathbb{C}))$.
Let $a+b=c+d$ and $b=$ min$\{a,b,c,d\}$.

\begin{eqnarray}
\pspicture[.45](-.5,-.3)(1.3,1.3)
\rput[t](0,-.1){$c$}\qline(0,0)(0,.4)
\rput[t](.7,-.1){$d$}\qline(.7,0)(.7,.4)
\psframe[linecolor=darkred](-.2,.4)(.9,.6)
\qline(0,.6)(0,1)\qline(.7,.6)(.7,1) \rput[b](0,1.1){$a$}
\rput[b](.7,1.1){$b$}
\endpspicture
= \sum_{k=0}^{b} a_k \pspicture[.45](-2.8,-2)(2.8,2)
\pcline(-1.3,-.9)(-1.3,.9)\Aput{$c-k$}
\pcline(1.3,.9)(1.3,-.9)\Aput{$b-k$}
\pccurve[angleA=270,angleB=270,ncurv=1](.7,.9)(-.7,.9)
\pccurve[angleA=90,angleB=90,ncurv=1](-.7,-.9)(.7,-.9)\rput(0,.8){$k$}
\rput(0,-.8){$k$} \qline(-1,-1.5)(-1,-1.1)\qline(-1,1.1)(-1,1.5)
\qline( 1,-1.5)( 1,-1.1)\qline( 1,1.1)( 1,1.5)
\psarc(0,.9){.9}{180}{270}\psarc(0,-.9){.9}{0}{90}
\rput[b](-1,1.7){$a$}\rput[b](1,1.7){$b$} \rput[t](-1,-1.7){$c$}
\rput[t](1,-1.7){$d$} \rput[br](.9,0){$d-b$}
\psframe[linecolor=darkred](-1.5,-1.1)(-.5,-.9)
\psframe[linecolor=darkred]( 1.5,-1.1)( .5,-.9)
\psframe[linecolor=darkred](-1.5, 1.1)(-.5, .9)
\psframe[linecolor=darkred]( 1.5, 1.1)( .5, .9)
\endpspicture
\label{a1abexp1}
\end{eqnarray}

\begin{prop}
The coefficients in equation~\ref{a1abexp1} are
$$a_k= \frac{[c]![b]![a+b-k]!}{[c-k]![b-k]![k]![a+b]!}.$$
\end{prop}
\begin{proof}
We induct on $a+b$. If $a+b=1$, it is clear. Without loss of
generality, we assume that $a\ge b$. Denote the diagram
corresponding to the coefficient $a_k$ in the right hand side of
equation~\ref{a1abexp1} by $D(k)$. By applying a single clasp
expansion for the clasp of weight $a+b$, then a single clasp
expansion of the clasp located at the northeast corner, we get

$$
\pspicture[.45](-.5,-.3)(1.3,1.3)
\rput[t](0,-.1){$c$}\qline(0,0)(0,.4)
\rput[t](.7,-.1){$d$}\qline(.7,0)(.7,.4)
\psframe[linecolor=darkred](-.2,.4)(.9,.6)
\qline(0,.6)(0,1)\qline(.7,.6)(.7,1) \rput[b](0,1.1){$a$}
\rput[b](.7,1.1){$b$}
\endpspicture
 = \pspicture[.45](-1.9,-2)(2,2)
\pcline(-1,.1)(-1,.9) \pcline(-1,-.9)(-1,-.1) \pcline(1,.9)(1,.1)
\pcline(1,-.1)(1,-.9) \pcline(1.4,.9)(1.4,-.9)
\pcline(-1,-1.5)(-1,-1.1) \pcline(-1,1.1)(-1,1.5)
\pcline(1,-1.1)(1,-1.5) \pcline(1,1.5)(1,1.1)
\rput[b](-1,1.7){$a$} \rput[b](1,1.7){$b$}
\rput[br](-1.2,.4){$a$}\rput[br](.8,.4){$b-1$}
\rput[br](-1.2,-.6){$c$}\rput[br](.8,-.6){$d-1$}
\rput[t](-1,-1.7){$c$}\rput[t](1,-1.7){$d$}
\psframe[linecolor=darkred](-1.5,-1.1)(-.5,-.9)
\psframe[linecolor=darkred]( 1.5,-1.1)( .5,-.9)
\psframe[linecolor=darkred](-1.5, 1.1)(-.5, .9)
\psframe[linecolor=darkred]( 1.5, 1.1)( .5, .9)
\psframe[linecolor=darkred](-1.3,-.1)(1.3, .1)
\endpspicture
-\frac{[a][c]}{[a+b][a+b-1]} \pspicture[.45](-2.4,-2)(2.2,2)
\pcline(-1,.1)(-1,.9)\Aput{$a-1$}
\pcline(-1,-.9)(-1,-.1)\Aput{$c-1$}
\pcline(1,.9)(1,.1)\Aput{$b-1$} \pcline(1,-.1)(1,-.9)\Aput{$d-1$}
\pccurve[angleA=135,angleB=45,ncurv=1](.7,-.9)(-.7,-.9)
\pccurve[angleA=225,angleB=315,ncurv=1](.7,.9)(-.7,.9)
\rput[b](-1,1.7){$a$}\rput[b](1,1.7){$b$}
\rput[t](-1,-1.7){$c$}\rput[t](1,-1.7){$d$}
\pcline(-1,-1.5)(-1,-1.1) \pcline(-1,1.1)(-1,1.5)
\pcline(1,-1.1)(1,-1.5) \pcline(1,1.5)(1,1.1)
\psframe[linecolor=darkred](-1.5,-1.1)(-.5,-.9)
\psframe[linecolor=darkred]( 1.5,-1.1)( .5,-.9)
\psframe[linecolor=darkred](-1.5, 1.1)(-.5, .9)
\psframe[linecolor=darkred]( 1.5, 1.1)( .5, .9)
\psframe[linecolor=darkred](-1.3,-.1)(1.3, .1)
\endpspicture
\label{a1abexp2}
$$

By induction, the right side equals

\begin{align*}
&\sum_{k=0}^{b-1}\frac{[c]![b-1]![a+b-1-k]!}
{[c-k]![b-1-k]![k]![a+b-1]!}D(k)\\
&+\frac{[a][c]}{[a+b][a+b-1]}\sum_{k=0}^{b-1}
\frac{[c-1]![b-1]![a+b-2-k]!}{[c-1-k]![b-1-k]![k]![a+b-2]!}D(k+1)\\
&=D(0)+\sum_{k=1}^{b-1}(\frac{[c]![b-1]![a+b-1-k]!}
{[c-k]![b-1-k]![k]![a+b-1]!}\\
&+\frac{[a][c]}{[a+b][a+b-1]}\frac{[c-1]![b-1]![a+b-2-(k-1)]!}
{[c-1-(k-1)]![b-1-(k-1)]![k-1]![a+b-2]!}) D(k)\\
&+\frac{[a][c]}{[a+b][a+b-1]}\frac{[c-1]![b-1]![a+b-2-(b-1)]!}
{[c-1-(b-1)]![b-1-(b-1)]![b-1]![a+b-2]!} D(b)\\
&=D(0)+\sum_{k=1}^{b-1}\frac{[c]![b]![a+b-k]!}
{[c-k]![b-k]![k]![a+b]!}
(\frac{[b-k][a+b]+[k][a]}{[b][a+b-k]})D(k)\\
& +\frac{[c]![b]![a]!}{[c-b]![0]![b]![a+b]!}D(b)
=\sum_{k=0}^{b}\frac{[c]![b-1]![a+b-k]!}{[c-k]![b-k]![k]![a+b]!}D(k).
\end{align*}
We use a well-known identity for quantum integers,
$$[m+r][n+r]=[m][n]+[m+n+r][r]$$
in the $6$-th line of the above equation with $n=-a, m=-k$ and
$r=a+b$.
\end{proof}

Next we look at the trihedron coefficient(or $3j$ symbol)
~\cite{Lickorish:gtm}~\cite{MV:3valent}~\cite{Turaev:quantum}.

\begin{prop}
The trihedron coefficient is
$$
\pspicture[.4](-2.2,-2.2)(2.2,2.2)
\psframe[linecolor=darkred](-1.5,-1.1)(-.5,-.9)
\psframe[linecolor=darkred]( 1.5,-1.1)( .5,-.9)
\psframe[linecolor=darkred](-1.5, 1.1)(-.5, .9)
\psframe[linecolor=darkred]( 1.5, 1.1)( .5, .9)
\psframe[linecolor=darkred](-.5,.1)(.5,-.1)
\pcline(-1.3,-.9)(-.3,-.1) \pcline(1.3,.9)(.3,.1)
\pcline(1.3,-.9)(.3,-.1) \pcline(-1.3,.9)(-.3,.1)
\pccurve[angleA=225,angleB=315,ncurv=1](.7,.9)(-.7,.9)
\pccurve[angleA=45,angleB=135,ncurv=1](-.7,-.9)(.7,-.9)\rput(0,.8){$j$}
\rput(0,-.8){$j$} \qline(-1,-1.5)(-1,-1.1)\qline(-1,1.1)(-1,1.5)
\qline(1,-1.5)(1,-1.1) \qline(1,1.1)(1,1.5) \qline(1,-1.5)(2,-1.5)
\qline(2,-1.5)(2,1.5) \qline(2,1.5)(1,1.5) \qline(-1,1.5)(-2,1.5)
\qline(-2,1.5)(-2,-1.5) \qline(-2,-1.5)(-1,-1.5)
 \rput[b](-1,.3){$i$}\rput[c](1,.3){$k$}
\endpspicture
=(-1)^{i+j+k}\frac{[i+j+k+1]![i]![j]![k]!}{[i+j]![j+k]![i+k]!}.
\label{a1triexp1}
$$
\end{prop}
\begin{proof}
The idea of the proof is identical to the previous proposition. We
induct on $i+j+k$. If $i+j+k=1$, it is just a circle, so its value
is $$-[2]=(-1)^1\frac{[2]![1]![0]![0]!}{[1]![1]![0]!}.$$ We apply
a single clasp expansion for the clasp of weight $i+k$, then
another single clasp expansion, for which the clasp is located at
the northeast corner. By induction, we have
\begin{align*}
&=-\frac{[j+k+1]}{[j+k]}(-1)^{i+j+k-1}
\frac{[i+j+k]![i]![j]![k-1]!}{[i+j]![i+k-1]![j+k-1]!} \\
&+\frac{[i][i]}{[i+k][i+k-1]}(-1)^{i+j+k-1}
\frac{[i+j+k]![i-1]![j+1]![k-1]!}{[i+j]![i+k-2]![j+k]!}\\
&=(-1)^{i+j+k}\frac{[i+j+k+1]![i]![j]![k]!}{[i+j]![i+k]![j+k]!}
(\frac{[i+k][j+k]}{[i+j+k+1][k]}-\frac{[i][j+1]}{[i+j+k+1][k]})\\
&=(-1)^{i+j+k}\frac{[i+j+k+1]![i]![j]![k]!}{[i+j]![i+k]![j+k]!}
\end{align*}
\end{proof}

\section{Single Clasp Expansion for $\mathcal{U}_q({\mathfrak sl}(3,\mathbb{C}))$}

A complete set of relations for skein theory of
$\mathcal{U}_q({\mathfrak sl}(3,\mathbb{C}))$ as given in
equations \ref{a2defr21} \ref{a2defr22} \ref{a2defr23} was found
in \cite{Kuperberg:spiders}. There is a relation for every {\it
elliptic face}, a face with less than 6 edges. We call the
relation~\ref{a2defr23} a {\it rectangular relation} and the
first(second) shape in the right side of the equality is called a
{\it horizontal(vertical}, respectively) {\it splitting}. For
several reasons, such as positivity and integrality, we use $-[2]$
in relation~\ref{a2defr22} but one can use a quantum integer $[2]$
and get an independent result. If one uses $[2]$, one can rewrite
all results in here by multiplying each trivalent vertex by the
complex number $i$.

\begin{eqnarray}
\pspicture[.45](-.6,-.5)(.6,.5)
\pscircle(0,0){.4}\psline[arrows=->,arrowscale=1.5](.1,.4)(.11,.4)
\endpspicture
& = & [3]  \label{a2defr21}
\end{eqnarray}

\begin{eqnarray}
 \pspicture[.45](-1.5,-.8)(1.5,.8)
\pnode(.4;180){a2} \pnode(.4;0){a3}
\qline(-1.2,0)(-.4,0)\psline[arrowscale=1.5]{->}(-.7,0)(-.9,0)
\nccurve[angleA=90,angleB=90,nodesep=1pt]{a2}{a3}\middlearrow
\nccurve[angleA=-90,angleB=-90,nodesep=1pt]{a2}{a3}\middlearrow
\qline(1.2,0)(.4,0)\psline[arrowscale=1.5]{->}(.9,0)(.7,0)
\endpspicture
& = & - [2] \pspicture[.45](-.8,-.8)(.8,.8) \qline(.6,0)(-.6,0)
\psline[arrowscale=1.5]{->}(.1,0)(-.1,0)
\endpspicture  \label{a2defr22}
\end{eqnarray}

\begin{eqnarray}
\pspicture[.45](-1.1,-1.1)(1.1,1.1) \qline(1,1)(.5,.5)
\qline(-1,1)(-.5,.5) \qline(-1,-1)(-.5,-.5) \qline(1,-1)(.5,-.5)
\qline(.5,.5)(-.5,.5) \qline(-.5,.5)(-.5,-.5)
\qline(-.5,-.5)(.5,-.5) \qline(.5,-.5)(.5,.5)
\psline[arrowscale=1.5]{<-}(.1,.5)(-.1,.5)
\psline[arrowscale=1.5]{<-}(-.1,-.5)(.1,-.5)
\psline[arrowscale=1.5]{<-}(.5,.1)(.5,-.1)
\psline[arrowscale=1.5]{<-}(-.5,-.1)(-.5,.1)
\psline[arrowscale=1.5]{->}(.85,.85)(.65,.65)
\psline[arrowscale=1.5]{<-}(-.85,.85)(-.65,.65)
\psline[arrowscale=1.5]{->}(-.85,-.85)(-.65,-.65)
\psline[arrowscale=1.5]{<-}(.85,-.85)(.65,-.65)
\endpspicture
&= &  \pspicture[.45](-1.1,-1.1)(1.1,1.1)
\pnode(1;45){a1}\pnode(1;135){a2}\pnode(1;225){a3}\pnode(1;315){a4}
\nccurve[angleA=225,angleB=315]{a1}{a2}\middlearrow
\nccurve[angleA=45,angleB=135]{a3}{a4}\middlearrow
\endpspicture + \pspicture[.45](-1.1,-1.1)(1.1,1.1)
\pnode(1;45){a1}\pnode(1;135){a2}\pnode(1;225){a3}\pnode(1;315){a4}
\nccurve[angleA=225,angleB=135]{a1}{a4}\middlearrow
\nccurve[angleA=45,angleB=315]{a3}{a2}\middlearrow
\endpspicture
\label{a2defr23}
\end{eqnarray}

A clasp for $\mathcal{U}_q({\mathfrak sl}(3,\mathbb{C}))$ can be
defined axiomatically: 1) it is an idempotent and 2) if we attach
a $U$ turn or a $Y$, it becomes zero. An explicit definition of
clasps for $\mathcal{U}_q({\mathfrak sl}(3,\mathbb{C}))$ can be
found in~\cite{Kuperberg:spiders}.

First we look at a single clasp expansion of the clasp of weight
$(a,0)$ where the weight $(a,b)$ stands for $a\lambda_1
+b\lambda_2$ and $\lambda_i$ is a fundamental dominant weight of
${\mathfrak sl}(3,\mathbb{C})$. Each directed edge represents
$V_{\lambda_i}$, the fundamental representation of the highest
weight $\lambda_i$. We might use the notation $+, -$ for
$V_{\lambda_1}, V_{\lambda_2}$ but it should be clear.

We recall the usual partial ordering of the weight lattice of
lattice of ${\mathfrak sl}(3,\mathbb{C})$ as
\begin{eqnarray*}
a\lambda_1 + b\lambda_2 & \succ & (a+1)\lambda_1 + (b-2)\lambda_2 \\
a\lambda_1 + b\lambda_2 & \succ & (a-2)\lambda_1 + (b+1)\lambda_2.
\end{eqnarray*}

A {\it cut path} is a path which is transverse to strings of web,
and the weight of a cut pass is $(a,b)$ if it passes $a$ strings
decorated by $V_{\lambda_1}$ and $b$ strings decorated by
$V_{\lambda_1}$.

\subsection{Single Clasp Expansions of a Clasp of Weight $(n,0)$ and $(0,n)$}

A basis for the single clasp expansion is given in
equation~\ref{a2exp}. If we attach a $Y$ on the top of webs in the
equation~\ref{a2exp}, there is at least one elliptic face on which
we can apply our relations. This process gives us exactly the same
equations we got for $\mathcal{U}_q({\mathfrak
sl}(2,\mathbb{C}))$. Thus, we can easily establish
proposition~\ref{a2exppro}. Moreover, this single clasp expansion
holds for any $\mathcal{U}_q({\mathfrak sl}(n,\mathbb{C}))$ where
$n\ge 4$ because $\mathcal{U}_q({\mathfrak sl}(3,\mathbb{C}))$ is
naturally embedded in $\mathcal{U}_q({\mathfrak
sl}(n,\mathbb{C}))$. Later we will mention the importance of this
fact. As same as for $\mathcal{U}_q({\mathfrak
sl}(2,\mathbb{C}))$, there are four different positions for the
single clasp expansions for $\mathcal{U}_q({\mathfrak
sl}(3,\mathbb{C}))$ so we use the same convention we used for
single clasp expansion for $\mathcal{U}_q({\mathfrak
sl}(2,\mathbb{C}))$.

\begin{eqnarray}
\pspicture[.45](-.1,-.3)(2.6,2.3)
\rput[t](1,-.1){$n$}\qline(1,0)(1,.4) \rput[c](1,1.2){$\ldots$}
\qline(.1,.6)(.1,2) \qline(.3,.6)(.3,2) \qline(.5,.6)(.5,2)
\qline(1.5,.6)(1.5,2) \qline(1.7,.6)(1.7,2) \qline(1.9,.6)(1.9,2)
\psline[arrowscale=1.5]{->}(.1,1)(.1,1.2)
\psline[arrowscale=1.5]{->}(.3,1)(.3,1.2)
\psline[arrowscale=1.5]{->}(.5,1)(.5,1.2)
\psline[arrowscale=1.5]{->}(1.5,1)(1.5,1.2)
\psline[arrowscale=1.5]{->}(1.7,1)(1.7,1.2)
\psline[arrowscale=1.5]{->}(1.9,1)(1.9,1.2) \rput[b](1,2.1){$n$}
\psframe[linecolor=darkred](0,.4)(2,.6)
\endpspicture
= \sum_{i=1}^{n} a_i \pspicture[.45](-.6,-.3)(3.35,2.5)
\rput[t](1,-.1){$n-1$} \rput[b](.4,1.8){$\ldots$}
\psline(1,0)(1,.4) \psline(.1,.6)(.1,2.2) \psline(.7,.6)(.7,2.2)
\psline(.9,1.8)(.9,2.2) \psline(.9,1.8)(1.1,1.6)
\psline(1.1,1.6)(1.1,1.4) \psline(1.1,1.4)(1.3,.9)
\psline(1.3,1.8)(1.1,1.6) \psline(1.3,1.8)(1.3,2.2)
\psline(1.3,.6)(1.3,.9) \psline(1.5,1.4)(1.3,.9)
\psline(1.5,1.4)(1.5,2.2) \psline(1.5,1.4)(1.7,.9)
\psline(2.2,1.4)(2,.9) \psline(2.2,1.4)(2.4,.9)
\psline(2.2,1.4)(2.2,2.2) \psline(2.4,.9)(2.6,1.4)
\psline(2.4,.6)(2.4,.9) \psline(2.6,1.4)(2.8,.9)
\psline(2.6,1.4)(2.6,2.2) \psline(2.8,.9)(3,1.4)
\psline(2.8,.9)(2.8,.6) \psline(3,1.4)(3.2,.9)
\psline(3.2,.9)(3.2,0) \psline[arrowscale=1.5]{->}(.1,1.3)(.1,1.5)
\psline[arrowscale=1.5]{->}(.7,1.3)(.7,1.5)
\psline[arrowscale=1.5]{->}(.9,1.9)(.9,2.1)
\psline[arrowscale=1.5]{->}(1.16,1.25)(1.24,1.05)
\psline[arrowscale=1.5]{->}(1.3,1.9)(1.3,2.1)
\psline[arrowscale=1.5]{->}(1.3,.65)(1.3,.85)
\psline[arrowscale=1.5]{->}(1.44,1.25)(1.36,1.05)
\psline[arrowscale=1.5]{->}(1.5,1.7)(1.5,1.9)
\psline[arrowscale=1.5]{->}(1.56,1.25)(1.64,1.05)
\psline[arrowscale=1.5]{->}(2.14,1.25)(2.06,1.05)
\psline[arrowscale=1.5]{->}(2.26,1.25)(2.34,1.05)
\psline[arrowscale=1.5]{->}(2.2,1.7)(2.2,1.9)
\psline[arrowscale=1.5]{<-}(2.46,1.05)(2.54,1.25)
\psline[arrowscale=1.5]{->}(2.4,.65)(2.4,.85)
\psline[arrowscale=1.5]{->}(2.64,1.25)(2.74,1.05)
\psline[arrowscale=1.5]{->}(2.6,1.7)(2.6,1.9)
\psline[arrowscale=1.5]{<-}(2.86,1.05)(2.94,1.25)
\psline[arrowscale=1.5]{->}(2.8,.65)(2.8,.85)
\psline[arrowscale=1.5]{<-}(3.2,.55)(3.2,.35)
\rput[b](.9,2.4){$i$} \rput[b](1.85,1.8){$\ldots$}
\rput[b](2.6,2.4){$1$} \psframe[linecolor=darkred](0,.4)(3,.6)
\endpspicture
\label{a2exp}
\end{eqnarray}

\begin{prop}
The coefficients in equation~\ref{a2exp} are
$$a_i=\frac{[n+1-i]}{[n]}.$$
\label{a2exppro}
\end{prop}

Also, we can easily find a single clasp expansion of the clasp of
weight $(0,b)$ by reversing arrows.

\subsection{Single Clasp Expansions of a Non-segregated Clasp of
Weight $(a,b)$}

The most interesting case is a single clasp expansion of the clasp
of weight $(a,b)$. First of all, we find the dimension of
$\mathrm{Inv}(V_{\lambda_1}^{\otimes a+1}\otimes
V_{\lambda_2}^{\otimes b}\otimes V_{(b-1)\lambda_1+a\lambda_2})$
in lemma~\ref{lema2dim}.

\begin{lem}
$$\mathrm{dim}(\mathrm{Inv}(V_{\lambda_1}^{\otimes a+1}\otimes
V_{\lambda_2}^{\otimes b}\otimes
V_{(b-1)\lambda_1+a\lambda_2}))=\begin{cases}a &\mathrm{if} \hskip 1cm b=0\\
(a+1)b &\mathrm{if}\hskip 1cm b>0\end{cases}$$ \label{lema2dim}
\end{lem}

\begin{proof}
For $b=0$, the result follows from proposition~\ref{a2exppro}. We
induct on $b$ to show the multiplicity of
$V_{a\lambda_{1}+(b-1)\lambda_2}$, $V_{(a+1)\lambda_1+b\lambda_2}$
and $V_{(a-1)\lambda_1+(b+1)\lambda_2}$ in the decomposition of
$V_{1,0}^{\otimes a+1}\otimes V_{0,1}^{\otimes b}$ into
irreducible representations are $ (a+1)b$, $1$ and $a$
respectively. By a simple application of Schur's Lemma, we find
that the dimension of $\mathrm{Inv}(V_{\lambda_1}^{\otimes
a+1}\otimes V_{\lambda_2}^{\otimes b}\otimes
V_{(b-1)\lambda_1+a\lambda_2})$ is equal to the multiplicity of
$V_{a\lambda_{1}+(b-1)\lambda_2}$ in the decomposition of
$V_{\lambda_1}^{\otimes a+1}\otimes V_{\lambda_2}^{\otimes b}$
into irreducible representations.
\end{proof}

Lemma~\ref{lema2dim} work for any $A_n$ where $n\ge 2$ by
replacing $\lambda_2$ by $\lambda_{n}$. Moreover the single clasp
expansion of the clasp of weight $a\lambda_1+b\lambda_n$ as in
equation ~\ref{a2exp} is also true for any
$\mathcal{U}_q({\mathfrak sl}(n,\mathbb{C}))$ where $n\ge 3$.

We need a set of basis webs with nice rectangular order, but we
can not find one in the general case. Even if one find a basis,
they have many hexagonal faces which make it very difficult to get
numerical relations. So we start from an alternative, {\it
non-segregated} clasp. A non-segregated clasp is obtained from the
segregated clasp by attaching a sequence of $H$'s until we get the
desired shape of edge orientations. Fortunately, there is a
canonical way to find by putting $H$ from the leftmost string of
weight $\lambda_2$ or $-$ until it reach to the desired position.
In the following lemma, we will show the non-segregated clasp is
well defined. The following figure is an example of a
non-segregated clasp of weight $(2,3)$ and how to obtain it from a
segregated clasp of weight $(2,3)$.

$$
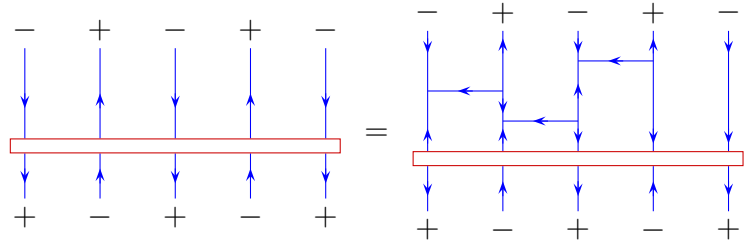
\pspicture[.42](-1.4,-1.4)(3.4,1.4) \uput{3pt}[90](0,1){$+$}
\uput{3pt}[90](-1,1){$-$} \uput{3pt}[270](-1,-1){$+$}
\uput{3pt}[90](1,1){$-$} \uput{3pt}[90](2,1){$+$}
\uput{3pt}[90](3,1){$-$} \uput{3pt}[270](0,-1){$-$}
\uput{3pt}[270](1,-1){$+$} \uput{3pt}[270](2,-1){$-$}
\uput{3pt}[270](3,-1){$+$}
\psline(0,-.2)(0,1)\psline[arrowscale=1.5]{->}(0,.2)(0,.4)
\psline(1,-.2)(1,1)\psline[arrowscale=1.5]{<-}(1,.2)(1,.4)
\psline(2,-.2)(2,1)\psline[arrowscale=1.5]{->}(2,.2)(2,.4)
\psline(3,-.2)(3,1)\psline[arrowscale=1.5]{<-}(3,.2)(3,.4)
\psline(-1,-.2)(-1,1)\psline[arrowscale=1.5]{<-}(-1,.2)(-1,.4)
\psline(-1,-.4)(-1,-1)\psline[arrowscale=1.5]{->}(-1,-.6)(-1,-.8)
\psline(0,-.4)(0,-1)\psline[arrowscale=1.5]{<-}(0,-.6)(0,-.8)
\psline(1,-.4)(1,-1)\psline[arrowscale=1.5]{->}(1,-.6)(1,-.8)
\psline(2,-.4)(2,-1)\psline[arrowscale=1.5]{<-}(2,-.6)(2,-.8)
\psline(3,-.4)(3,-1)\psline[arrowscale=1.5]{->}(3,-.6)(3,-.8)
\psframe[linecolor=darkred](-1.2,-.4)(3.2,-.2)
\endpspicture
= \pspicture[.42](-1.4,-1.4)(3.4,1.8) \uput{3pt}[90](0,1.4){$+$}
\uput{3pt}[90](-1,1.4){$-$} \uput{3pt}[270](-1,-1){$+$}
\uput{3pt}[90](1,1.4){$-$} \uput{3pt}[90](2,1.4){$+$}
\uput{3pt}[90](3,1.4){$-$} \uput{3pt}[270](0,-1){$-$}
\uput{3pt}[270](1,-1){$+$} \uput{3pt}[270](2,-1){$-$}
\uput{3pt}[270](3,-1){$+$}
\psline(0,-.2)(0,1.4)\psline[arrowscale=1.5]{->}(0,1.1)(0,1.3)
\psline(1,-.2)(1,1.4)\psline[arrowscale=1.5]{<-}(1,1.1)(1,1.3)
\psline(2,-.2)(2,1.4)\psline[arrowscale=1.5]{->}(2,1.1)(2,1.3)
\psline(3,-.2)(3,1.4)\psline[arrowscale=1.5]{<-}(3,1.1)(3,1.3)
\psline(-1,-.2)(-1,1.4)\psline[arrowscale=1.5]{<-}(-1,1.1)(-1,1.3)
\psline(-1,-.4)(-1,-1)\psline[arrowscale=1.5]{->}(-1,-.6)(-1,-.8)
\psline(0,-.4)(0,-1)\psline[arrowscale=1.5]{<-}(0,-.6)(0,-.8)
\psline(1,-.4)(1,-1)\psline[arrowscale=1.5]{->}(1,-.6)(1,-.8)
\psline(2,-.4)(2,-1)\psline[arrowscale=1.5]{<-}(2,-.6)(2,-.8)
\psline(3,-.4)(3,-1)\psline[arrowscale=1.5]{->}(3,-.6)(3,-.8)
\psline(-1,.6)(0,.6)\psline[arrowscale=1.5]{->}(-.4,.6)(-.6,.6)
\psline(1,.2)(0,.2)\psline[arrowscale=1.5]{->}(.6,.2)(.4,.2)
\psline(1,1)(2,1)\psline[arrowscale=1.5]{->}(1.6,1)(1.4,1)
\psline[arrowscale=1.5]{->}(-1,-.1)(-1,.1)\psline[arrowscale=1.5]{->}(0,-.1)(0,.1)
\psline[arrowscale=1.5]{<-}(1,-.1)(1,.1)\psline[arrowscale=1.5]{<-}(2,-.1)(2,.1)
\psline[arrowscale=1.5]{<-}(3,-.1)(3,.1)\psline[arrowscale=1.5]{<-}(0,.3)(0,.5)
\psline[arrowscale=1.5]{->}(1,.5)(1,.7)
\psframe[linecolor=darkred](-1.2,-.4)(3.2,-.2)
\endpspicture \label{nsclaspex}
$$

\begin{lem}
Non-segregated clasps are well-defined.
\end{lem}

\begin{proof}
Let $\alpha$ be a sequence of $H$'s which induce the same
non-segregated clasp. We find the first string from the leftmost
of sign $-$ which does move to right. If there is no such a
string, then the clasp is canonical. Otherwise, there exist two
consecutive $H$'s which can be removed by the horizontal splitting
because it had to reach the desired position. We induct on the
length of the sequence and it completes the proof.
\end{proof}

We also find that non-segregated clasps satisfy two properties of
segregated clasps.

\begin{lem}
1) Two consecutive non-segregated clasps is equal to a
non-segregated clasp.

2) If we attach a web to a non-segregated clasp and if it has a
cut path of which weight is less than the weight of the clasp,
then it is zero.
\end{lem}
\begin{proof} Since there is two consecutive non-segregated clasps,
we see that the bottom end of the upper non-segregated clasp and
the top end of the lower non-segregated clasp are the same
non-segregated. Since non-segregated clasp does not depend on the
choice of order of attaching $H$'s, we fix the canonical one for
the both side then we an see that one is the other's inverse, the
inverse is just the horizontal reflection. Therefore, all the
$H$'s in the middle cancel out and standard clasps are
idempotents.

The second part is obvious because adding $H$'s does not change
the weight of the minimal path so we can change to a segregated
clasp by adding more $H$'s. Then the clasp becomes zero.
\end{proof}

The following equation~\ref{a2abexp} is a single clasp expansion
of a non-segregated clasp of weight $(a,b)$. Let us denote the web
corresponding to the coefficient $a_{i,j}$ by $D_{i,j}$. These
webs form a basis because there are no applicable relations.

\begin{eqnarray}
\pspicture[.4](-.3,-.7)(2.5,2.5) \psline(.1,.1)(.1,2)
\psline(.1,-.1)(.1,-.6) \psline(.3,.1)(.3,2)
\psline(.3,-.1)(.3,-.6) \psline(.5,.1)(.5,2)
\psline(.5,-.1)(.5,-.6) \psline(1.5,.1)(1.5,2)
\psline(1.5,-.1)(1.5,-.6) \psline(1.7,.1)(1.7,2)
\psline(1.7,-.1)(1.7,-.6) \psline(1.9,.1)(1.9,2)
\psline(1.9,-.1)(1.9,-.5) \psline(1.9,-.5)(2.1,-.5)
\psline(2.1,-.5)(2.1,2) \uput{3pt}[90](1,.7){$\ldots$}
\uput{3pt}[90](.3,2.1){$b$} \uput{3pt}[90](1.7,2.1){$a+1$}
\uput{3pt}[270](1,-.3){$\ldots$}
\psline[arrowscale=1.5]{<-}(.1,.8)(.1,1)
\psline[arrowscale=1.5]{<-}(.3,.8)(.3,1)
\psline[arrowscale=1.5]{<-}(.5,.8)(.5,2)
\psline[arrowscale=1.5]{->}(1.5,.8)(1.5,1)
\psline[arrowscale=1.5]{->}(1.7,.8)(1.7,1)
\psline[arrowscale=1.5]{->}(1.9,.8)(1.9,1)
\psline[arrowscale=1.5]{->}(2.1,.8)(2.1,1)
\psframe[linecolor=darkred](0,-.1)(2,.1)
\endpspicture
= \sum_{i=1}^{b}\sum_{j=0}^{a} \hskip .3cm a_{ij}
\pspicture[.4](-.4,-.8)(7,2.5) \psline(.1,.1)(.1,2)
\psline(.1,-.1)(.1,-.6) \psline(.3,.1)(.3,2)
\psline(.3,-.1)(.3,-.6) \psline(1.3,.1)(1.3,2)
\psline(1.3,-.1)(1.3,-.6) \psline(1.5,1)(1.5,2)
\psline(1.5,1)(1.8,1) \psline(2.2,2)(2.2,1.3)
\psline(2.2,.7)(2.2,.1) \psline(2.2,-.1)(2.2,-.6)
\psline(3,1)(3.8,1) \psline(4.6,1.3)(4.6,2)
\psline(4.6,-.1)(4.6,-.6) \psline(4.6,.7)(4.6,.1)
\psline(5,1)(5.3,1) \psline(5.3,1)(5.3,2)
\psline(5.5,-.1)(5.5,-.6) \psline(5.5,.1)(5.5,2)
\psline(6.5,-.1)(6.5,-.6) \psline(6.5,.1)(6.5,2)
\psline(6.7,-.1)(6.7,-.6) \psline(6.7,.1)(6.7,2)
\uput{3pt}[90](1.5,2.1){$1$} \uput{3pt}[90](2.2,2.1){$i-1$}
\uput{3pt}[90](4.6,2.1){$j$} \uput{3pt}[90](5.3,2.1){$1$}
\uput{3pt}[0](2.2,0.3){$i-1$} \uput{3pt}[0](4.6,0.3){$j$}
\uput{3pt}[270](.8,1.3){$\ldots$} \uput{3pt}[270](6,1.3){$\ldots$}
\uput{3pt}[270](.8,-.3){$\ldots$} \uput{3pt}[270](6,-.3){$\ldots$}
\psline[arrowscale=1.5]{<-}(.1,.8)(.1,1)
\psline[arrowscale=1.5]{<-}(.3,.8)(.3,1)
\psline[arrowscale=1.5]{<-}(1.3,.8)(1.3,1)
\psline[arrowscale=1.5]{->}(1.5,1)(1.8,1)
\psline[arrowscale=1.5]{->}(2.2,1.8)(2.2,1.6)
\psline[arrowscale=1.5]{->}(2.2,.5)(2.2,.3)
\psline[arrowscale=1.5]{->}(3.1,1)(3.3,1)
\psline[arrowscale=1.5]{->}(3.5,1)(3.7,1)
\psline[arrowscale=1.5]{->}(4.6,1.6)(4.6,1.8)
\psline[arrowscale=1.5]{<-}(4.6,.5)(4.6,.3)
\psline[arrowscale=1.5]{->}(5,1)(5.3,1)
\psline[arrowscale=1.5]{->}(5.5,.8)(5.5,1)
\psline[arrowscale=1.5]{->}(6.5,.8)(6.5,1)
\psline[arrowscale=1.5]{->}(6.7,.8)(6.7,1)
\psframe[linecolor=darkred](0,-.1)(6.8,.1)
\psframe[linecolor=emgreen](1.8,.7)(3,1.3)
\psframe[linecolor=emgreen](3.8,.7)(5,1.3)
\psline[linecolor=darkred,linestyle=dashed](3.4,2.2)(3.4,-.7)
\uput{1pt}[90](2.4,.7){$(1)$} \uput{1pt}[90](4.4,.7){$(2)$}
\endpspicture
\label{a2abexp}
\end{eqnarray}

Kuperberg~\cite{Kuperberg:spiders}  showed that for a fixed
boundary, interior can be filled by a cut out from the hexagonal
tiling of the plane with the given boundary. For our cases, there
are two possible fillings but we use the maximal cut out of the
hexagonal tiling. We draw examples of the case $i=6, j=5$ and the
first one in equation~\ref{hexaex} is not maximal cut out and the
second one is the maximal cut out which fits to the left rectangle
and the last one is the maximal cut out which fits to the right
rectangle as the number indicates in equation~\ref{a2abexp}.

\begin{eqnarray}
\pspicture[.5](-.6,-2.1)(3.9,1.5)
\psline(-.5,0)(0,0)\psline[arrowscale=1.5]{->}(-.35,0)(-.15,0)
\psline(3.3,0)(3.8,0)\psline[arrowscale=1.5]{->}(3.45,0)(3.65,0)
\multips(.3,.3)(.6,0){5}{\psline(0,0)(0,1)\psline[arrowscale=1.5]{->}(0,.6)(0,.4)}
\multips(.6,-1.3)(.6,0){5}{\psline(0,0)(0,1)\psline[arrowscale=1.5]{->}(0,.6)(0,.4)}
\multips(.3,.3)(.6,0){4}{\psline(0,0)(0.9,-.6)}
\psline(0,0)(.6,-.3)\psline(3.3,0)(2.7,.3)
\psframe[linecolor=emgreen](0,-1)(3.3,1)
\endpspicture
, \pspicture[.5](-.6,-2.1)(3.9,1.5)
\psline(-.5,0)(0,0)\psline[arrowscale=1.5]{->}(-.35,0)(-.15,0)
\psline(3.3,0)(3.8,0)\psline[arrowscale=1.5]{->}(3.45,0)(3.65,0)
\multips(.3,.3)(.6,0){5}{\psline(0,0)(0,1)\psline[arrowscale=1.5]{->}(0,.6)(0,.4)}
\multips(.6,-1.3)(.6,0){5}{\psline(0,0)(0,1)\psline[arrowscale=1.5]{->}(0,.6)(0,.4)}
\multips(.3,.3)(.6,0){5}{\psline(0,0)(0.3,-.6)\psline[arrowscale=1.5]{->}(0.1,-.2)(0.2,-.4)}
\multips(.6,-.3)(.6,0){4}{\psline(0,0)(0.3,.6)\psline[arrowscale=1.5]{->}(0.1,.2)(0.2,.4)}
\psline(0,0)(.3,.3)\psline(3.3,0)(3,-.3)
\psframe[linecolor=emgreen](0,-1)(3.3,1)
\uput{3pt}[90](1.8,-2.1){$(1)$}
\endpspicture
, \pspicture[.5](-.6,-2.1)(3.9,1.5)
\psline(-.5,0)(0,0)\psline[arrowscale=1.5]{->}(-.35,0)(-.15,0)
\psline(3.3,0)(3.8,0)\psline[arrowscale=1.5]{->}(3.45,0)(3.65,0)
\multips(.6,.3)(.6,0){5}{\psline(0,0)(0,1)\psline[arrowscale=1.5]{<-}(0,.6)(0,.4)}
\multips(.3,-1.3)(.6,0){5}{\psline(0,0)(0,1)\psline[arrowscale=1.5]{<-}(0,.6)(0,.4)}
\multips(.6,.3)(.6,0){4}{\psline(0,0)(0.3,-.6)\psline[arrowscale=1.5]{<-}(0.1,-.2)(0.2,-.4)}
\multips(.3,-.3)(.6,0){5}{\psline(0,0)(0.3,.6)\psline[arrowscale=1.5]{<-}(0.1,.2)(0.2,.4)}
\psline(0,0)(.3,-.3)\psline(3.3,0)(3,.3)
\psframe[linecolor=emgreen](0,-1)(3.3,1)
\uput{3pt}[90](1.8,-2.1){$(2)$}
\endpspicture\label{hexaex}
\end{eqnarray}

\begin{thm}
The coefficients in equation~\ref{a2abexp} are
$$a_{i,j}=\frac{[b-i+1]}{[b]}\frac{[b+j+1]}{[a+b+1]}.$$
\label{a2abexpthm1}
\end{thm}
\begin{proof}
As usual, we attach a $Y$ or a $U$ turn to find one exceptional
and three types of equations as follow.

$$[3]a_{1,0}-[2]a_{1,1}-[2]a_{2,0}+a_{2,1}=0.$$

Type I : For $j=0,1,\ldots, a$,

$$a_{b-1,j}-[2]a_{b,j}=0.$$

Type II : For $i=1,2, \ldots, b-2$ and $j=0,1,\ldots, a$.

$$a_{i,j}-[2]a_{i+1,j}+a_{i+2,j}=0.$$

Type III : For $i=1,2, \ldots, b$ and $j=0,1,\ldots, a-2$.

$$a_{i,j}-[2]a_{i,j+1}+a_{i,j+2}=0.$$

We establish the following lemma~\ref{a2abexplem1} first.
\begin{lem}
Let $a_{1,0}=x$, then the coefficients in the
equation~\ref{a2abexp} is
$$a_{i,j}=\frac{[b-i+1]}{[b]}\frac{[b+j+1]}{[b+1]}x.$$
\label{a2abexplem1}
\end{lem}
\begin{proof}
First we can see that the right side of equation~\ref{a2abexp} is
a basis for a single clasp expansion because its number of the
webs is equal to the dimension as in lemma~\ref{lema2dim} and none
of these webs has any faces. Second we find that these equations
has at least $(a+1)b$ independent equations. Then we plug in these
coefficients to equations to check that they are the right
coefficients.
\end{proof}

Usually we normalize one basis web in the expansion to get a known
value. But we can not normalize for this expansion yet because it
is not a segregated clasp. Thus we use lemma~\ref{a2abexplem3} to
find that the coefficient of $a_{1,a}$ is $1$. Then, we get
$x=\frac{[b+1]}{[a+b+1]}$ and it completes the proof of the
theorem.
\end{proof}

\subsection{Double Clasp Expansion of a Segregated Clasp of
Weight $(a,b)$}

Now we study a single and a double clasps expansion of a
segregated clasp of weight $(a,b)$. We will start with an example,
a single clasp expansion of a clasp of weight $(2,2)$. To apply
theorem~\ref{a2abexpthm1}, we add some H's to change the
segregated clasp to a non-segregated clasp.

$$
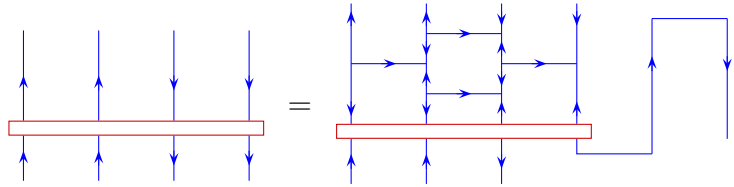
\pspicture[.42](-1.4,-1.2)(2.4,1.4)
\psline(-1,-.2)(-1,1)\psline[arrowscale=1.5]{->}(-1,.2)(-1,.4)
\psline(0,-.2)(0,1)\psline[arrowscale=1.5]{->}(0,.2)(0,.4)
\psline(1,-.2)(1,1)\psline[arrowscale=1.5]{<-}(1,.2)(1,.4)
\psline(2,-.2)(2,1)\psline[arrowscale=1.5]{<-}(2,.2)(2,.4)
\psline(-1,-.4)(-1,-1)\psline[arrowscale=1.5]{<-}(-1,-.6)(-1,-.8)
\psline(0,-.4)(0,-1)\psline[arrowscale=1.5]{<-}(0,-.6)(0,-.8)
\psline(1,-.4)(1,-1)\psline[arrowscale=1.5]{->}(1,-.6)(1,-.8)
\psline(2,-.4)(2,-1)\psline[arrowscale=1.5]{->}(2,-.6)(2,-.8)
\psframe[linecolor=darkred](-1.2,-.4)(2.2,-.2)
\endpspicture
= \pspicture[.42](-1.4,-1.2)(4.2,1.5)
\psline(-1,-.2)(-1,1.4)\psline[arrowscale=1.5]{->}(-1,1.1)(-1,1.3)
\psline(0,-.2)(0,1.4)\psline[arrowscale=1.5]{->}(0,1.1)(0,1.3)
\psline(1,-.2)(1,1.4)\psline[arrowscale=1.5]{<-}(1,1.1)(1,1.3)
\psline(2,-.2)(2,1.4)\psline[arrowscale=1.5]{<-}(2,1.1)(2,1.3)
\psline(3,-.6)(3,1.2)\psline[arrowscale=1.5]{->}(3,.5)(3,.7)
\psline(2,-.4)(2,-.6)(3,-.6) \psline(3,1.2)(4,1.2)
\psline(-1,-.4)(-1,-1)\psline[arrowscale=1.5]{<-}(-1,-.6)(-1,-.8)
\psline(0,-.4)(0,-1)\psline[arrowscale=1.5]{<-}(0,-.6)(0,-.8)
\psline(1,-.4)(1,-1)\psline[arrowscale=1.5]{->}(1,-.6)(1,-.8)
\psline(4,-.4)(4,1.2)\psline[arrowscale=1.5]{->}(4,.7)(4,.5)
\psline(0,1)(1,1)\psline[arrowscale=1.5]{->}(.4,1)(.6,1)
\psline(-1,.6)(0,.6)\psline[arrowscale=1.5]{<-}(-.4,.6)(-.6,.6)
\psline(0,.2)(1,.2)\psline[arrowscale=1.5]{<-}(.6,.2)(.4,.2)
\psline(1,.6)(2,.6)\psline[arrowscale=1.5]{<-}(1.6,.6)(1.4,.6)
\psline[arrowscale=1.5]{<-}(-1,-.1)(-1,.1)\psline[arrowscale=1.5]{<-}(0,-.1)(0,.1)
\psline[arrowscale=1.5]{->}(1,-.1)(1,.1)\psline[arrowscale=1.5]{->}(2,-.1)(2,.1)
\psline[arrowscale=1.5]{->}(0,.3)(0,.5)\psline[arrowscale=1.5]{<-}(0,.7)(0,.9)
\psline[arrowscale=1.5]{->}(1,.5)(1,.3)\psline[arrowscale=1.5]{->}(1,.7)(1,.9)
\psframe[linecolor=darkred](-1.2,-.4)(2.2,-.2)
\endpspicture \label{nsclaspex22}
$$

Then we expand the clasp. In the following
equation~\ref{a2doubleex1}, we will omit the direction of edges
unless there is an ambiguity.

\begin{eqnarray}\nonumber  \frac{[3]}{[5]}
\pspicture[.42](-.4,-.6)(2.7,2)
\psframe[linecolor=darkred](-.2,-.2)(2.2,0) \psline(0,-.7)(0,-.2)
\psline(0,0)(0,1.7) \psline(0,1.1)(0.5,1.1)
\psline(0.5,1.7)(0.5,0.5) \psline(0.5,0.5)(1,0.5)
\psline(1,0.5)(1,1.7) \psline(0.5,0.9)(1,0.9)
\psline(0.5,1.3)(1,1.3) \psline(1,1.1)(1.5,1.1)
\psline(1.5,1.7)(1.5,0.0) \psline(2,0.0)(2,1.7)
\psline(2,1.7)(2.5,1.7) \psline(2.5,1.7)(2.5,-0.7)
\psline(0,-.7)(0,-.2) \psline(1.5,-.7)(1.5,-.2)
\psline(2,-0.7)(2,-0.2)
\endpspicture + \frac{[4]}{[5]}
\pspicture[.42](-.4,-.6)(2.7,2)
\psframe[linecolor=darkred](-.2,-.2)(2.2,0) \psline(0,-.7)(0,-.2)
\psline(1,-.7)(1,-.2) \psline(2,-0.7)(2,-0.2) \psline(0,0)(0,1.7)
\psline(0.5,1.7)(0.5,0.5) \psline(1,0.0)(1,1.7)
\psline(1.5,1.7)(1.5,0.7) \psline(2,0.0)(2,1.7)
\psline(2.5,1.7)(2.5,-0.7) \psline(0,1.1)(0.5,1.1)
\psline(0.5,0.5)(1,0.5) \psline(0.5,0.9)(1,0.9)
\psline(0.5,1.3)(1,1.3) \psline(1,1.1)(1.5,1.1)
\psline(1,0.7)(1.5,0.7) \psline(2,1.7)(2.5,1.7)
\endpspicture + \frac{[5]}{[5]}
\pspicture[.42](-.4,-.6)(2.7,2)
\psframe[linecolor=darkred](-.2,-.2)(2.2,0) \psline(0,-.7)(0,-.2)
\psline(1.5,-.7)(1.5,-.2) \psline(1,-.7)(1,-.2)
\psline(0,0)(0,1.7) \psline(0.5,1.7)(0.5,0.5)
\psline(1,0.0)(1,1.7) \psline(1.5,1.7)(1.5,0.0)
\psline(2,0.9)(2,1.7) \psline(2.5,1.7)(2.5,-0.7)
\psline(0,1.1)(0.5,1.1) \psline(0.5,0.5)(1,0.5)
\psline(0.5,0.9)(1,0.9) \psline(0.5,1.3)(1,1.3)
\psline(1,1.1)(1.5,1.1) \psline(1,0.7)(1.5,0.7)
\psline(1.5,0.9)(2,0.9) \psline(2,1.7)(2.5,1.7)
\endpspicture \\ + \frac{[1][2]}{[3][5]}
\pspicture[.42](-.4,-.6)(2.7,2)
\psframe[linecolor=darkred](-.2,-.2)(2.2,0)
\psline(0.5,-0.7)(0.5,-0.2) \psline(1.5,-.7)(1.5,-.2)
\psline(2,-0.7)(2,-0.2) \psline(0,.7)(0,1.7)
\psline(0.5,1.7)(0.5,0.0) \psline(1,0.5)(1,1.7)
\psline(1.5,1.7)(1.5,0.0) \psline(2,0.0)(2,1.7)
\psline(2.5,1.7)(2.5,-0.7) \psline(0,1.1)(0.5,1.1)
\psline(0,.7)(0.5,0.7) \psline(0.5,0.5)(1,0.5)
\psline(0.5,0.9)(1,0.9) \psline(0.5,1.3)(1,1.3)
\psline(1,1.1)(1.5,1.1) \psline(2,1.7)(2.5,1.7)
\endpspicture + \frac{[1][2]}{[4][5]}
\pspicture[.42](-.4,-.6)(2.7,2)
\psframe[linecolor=darkred](-.2,-.2)(2.2,0)
\psline(0.5,-0.7)(0.5,-0.2) \psline(1,-.7)(1,-.2)
\psline(2,-0.7)(2,-0.2) \psline(0,.7)(0,1.7)
\psline(0.5,1.7)(0.5,0.0) \psline(1,0.0)(1,1.7)
\psline(1.5,1.7)(1.5,0.7) \psline(2,0.0)(2,1.7)
\psline(2.5,1.7)(2.5,-0.7) \psline(0,1.1)(0.5,1.1)
\psline(0,.7)(0.5,0.7) \psline(0.5,0.5)(1,0.5)
\psline(0.5,0.9)(1,0.9) \psline(0.5,1.3)(1,1.3)
\psline(1,1.1)(1.5,1.1) \psline(1,0.7)(1.5,0.7)
\psline(2,1.7)(2.5,1.7)
\endpspicture
+ \frac{[1][2]}{[5][5]} \pspicture[.42](-.4,-.6)(2.7,2)
\psframe[linecolor=darkred](-.2,-.2)(2.2,0)
\psline(0.5,-0.7)(0.5,-0.2) \psline(1,-.7)(1,-.2)
\psline(1.5,-.7)(1.5,-.2) \psline(0,.7)(0,1.7)
\psline(0.5,1.7)(0.5,0.0) \psline(1,0.0)(1,1.7)
\psline(1.5,1.7)(1.5,0.0) \psline(2,0.9)(2,1.7)
\psline(2.5,1.7)(2.5,-0.7) \psline(0,1.1)(0.5,1.1)
\psline(0,.7)(0.5,0.7) \psline(0.5,0.5)(1,0.5)
\psline(0.5,0.9)(1,0.9) \psline(0.5,1.3)(1,1.3)
\psline(1,1.1)(1.5,1.1) \psline(1,0.7)(1.5,0.7)
\psline(1.5,0.9)(2,0.9) \psline(2,1.7)(2.5,1.7)
\endpspicture
\label{a2doubleex1}
\end{eqnarray}

These webs can be expanded using relations. For example, the first
one can be expanded by equation~\ref{a2defr23} as

$$
\pspicture[.42](-.4,-.6)(2.7,1.8)
\psframe[linecolor=darkred](-.2,-.2)(2.2,0) \psline(0,0)(0,1.7)
\psline(0,1.1)(0.5,1.1) \psline(0.5,1.7)(0.5,0.5)
\psline(0.5,0.5)(1,0.5) \psline(1,0.5)(1,1.7)
\psline(0.5,0.9)(1,0.9) \psline(0.5,1.3)(1,1.3)
\psline(1,1.1)(1.5,1.1) \psline(1.5,1.7)(1.5,0.0)
\psline(2,0.0)(2,1.7) \psline(2,1.7)(2.5,1.7)
\psline(2.5,1.7)(2.5,-0.7) \psline(0,-.7)(0,-.2)
\psline(1.5,-.7)(1.5,-.2) \psline(2,-0.7)(2,-0.2)
\endpspicture
=\pspicture[.42](-.4,-.6)(2.7,1.8)
\psframe[linecolor=darkred](-.2,-.2)(2.2,0) \psline(0,0)(0,1.7)
\psline(0.5,1.7)(0.5,1.1) \psline(1,1.1)(1,1.7)
\psline(1.5,1.7)(1.5,0.0) \psline(0,1.1)(0.5,1.1)
\psline(1,1.1)(1.5,1.1) \psline(2,0.0)(2,1.7)
\psline(2,1.7)(2.5,1.7) \psline(2.5,1.7)(2.5,-0.7)
\psline(0,-.7)(0,-.2) \psline(1.5,-.7)(1.5,-.2)
\psline(2,-0.7)(2,-0.2)
\pcline[linestyle=dotted,dotsep=1.3pt,linewidth=.7pt](0.5,1.3)(1,1.3)
\pcline[linestyle=dotted,dotsep=1.3pt,linewidth=.7pt](0.5,1.1)(0.5,0.5)
\pcline[linestyle=dotted,dotsep=1.3pt,linewidth=.7pt](1,1.1)(1,0.5)
\pcline[linestyle=dotted,dotsep=1.3pt,linewidth=.7pt](1,0.9)(0.5,0.9)
\pcline[linestyle=dotted,dotsep=1.3pt,linewidth=.7pt](1,0.5)(0.5,0.5)
\endpspicture -[2]
\pspicture[.42](-.4,-.6)(2.7,1.8)
\psframe[linecolor=darkred](-.2,-.2)(2.2,0) \psline(0,0)(0,1.7)
\psline(0.5,1.7)(0.5,1.3) \psline(1,1.3)(1,1.7)
\psline(1.5,1.7)(1.5,0.0) \psline(0,1.1)(0.5,1.1)
\psline(0.5,1.1)(0.5,0.9) \psline(0.5,0.9)(1,0.9)
\psline(1,0.9)(1,1.1) \psline(1,1.1)(1.5,1.1)
\psline(0.5,1.3)(1,1.3) \psline(2,0.0)(2,1.7)
\psline(2,1.7)(2.5,1.7) \psline(2.5,1.7)(2.5,-0.7)
\psline(0,-.7)(0,-.2) \psline(1.5,-.7)(1.5,-.2)
\psline(2,-0.7)(2,-0.2)
\pcline[linestyle=dotted,dotsep=1.3pt,linewidth=.7pt](0.5,1.1)(0.5,1.3)
\pcline[linestyle=dotted,dotsep=1.3pt,linewidth=.7pt](1,1.1)(1,1.3)
\pcline[linestyle=dotted,dotsep=1.3pt,linewidth=.7pt](0.5,0.9)(0.5,0.5)
\pcline[linestyle=dotted,dotsep=1.3pt,linewidth=.7pt](0.5,0.5)(1,0.5)
\pcline[linestyle=dotted,dotsep=1.3pt,linewidth=.7pt](1,0.5)(1,0.9)
\endpspicture
$$

For some small cases, we can expand this way but it will be
difficult to manage all possible expansions. An other way to look
at this expansion is to use paths : since there are five points on
the top and three points right above the clasp and these three
points have to be connected to points on the top (otherwise, we
have a cut path with weight less than $(2,1)$ which makes the web
zero), we have two $Y$'s or one $U$ turn. We first find all
possible disjoint, monotone(except at Y's) paths connecting these
points. For the web on above example, there are eight
possibilities as follows.

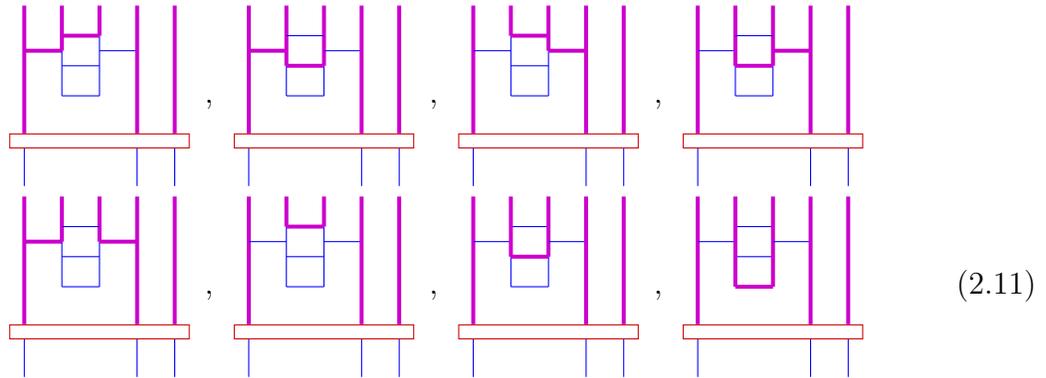
\begin{eqnarray}\nonumber
\pspicture[.42](-.4,-.6)(2.4,1.8)
\psframe[linecolor=darkred](-.2,-.2)(2.2,0)
\psline[linewidth=1.5pt,linecolor=pup](0,0)(0,1.7)
\psline[linewidth=1.5pt,linecolor=pup](0,1.1)(0.5,1.1)
\psline(0.5,1.1)(0.5,0.5)
\psline[linewidth=1.5pt,linecolor=pup](0.5,1.1)(0.5,1.7)
\psline(0.5,0.5)(1,0.5) \psline(1,0.5)(1,1.3)
\psline[linewidth=1.5pt,linecolor=pup](1,1.3)(1,1.7)
\psline(0.5,0.9)(1,0.9)
\psline[linewidth=1.5pt,linecolor=pup](0.5,1.3)(1,1.3)
\psline(1,1.1)(1.5,1.1)
\psline[linewidth=1.5pt,linecolor=pup](1.5,1.7)(1.5,0.0)
\psline[linewidth=1.5pt,linecolor=pup](2,0.0)(2,1.7)
\psline(0,-.7)(0,-.2) \psline(1.5,-.7)(1.5,-.2)
\psline(2,-0.7)(2,-0.2)
\endpspicture ,
\pspicture[.42](-.4,-.6)(2.4,1.8)
\psframe[linecolor=darkred](-.2,-.2)(2.2,0)
\psline[linewidth=1.5pt,linecolor=pup](0,0)(0,1.7)
\psline[linewidth=1.5pt,linecolor=pup](0,1.1)(0.5,1.1)
\psline(0.5,0.9)(0.5,0.5)
\psline[linewidth=1.5pt,linecolor=pup](0.5,0.9)(0.5,1.7)
\psline(0.5,0.5)(1,0.5) \psline(1,0.5)(1,0.9)
\psline[linewidth=1.5pt,linecolor=pup](1,0.9)(1,1.7)
\psline(0.5,1.3)(1,1.3)
\psline[linewidth=1.5pt,linecolor=pup](0.5,0.9)(1,0.9)
\psline(1,1.1)(1.5,1.1)
\psline[linewidth=1.5pt,linecolor=pup](1.5,1.7)(1.5,0.0)
\psline[linewidth=1.5pt,linecolor=pup](2,0.0)(2,1.7)
\psline(0,-.7)(0,-.2) \psline(1.5,-.7)(1.5,-.2)
\psline(2,-0.7)(2,-0.2)
\endpspicture ,
\pspicture[.42](-.4,-.6)(2.4,1.8)
\psframe[linecolor=darkred](-.2,-.2)(2.2,0)
\psline[linewidth=1.5pt,linecolor=pup](0,0)(0,1.7)
\psline(0,1.1)(0.5,1.1)
\psline[linewidth=1.5pt,linecolor=pup](0.5,1.7)(0.5,1.3)
\psline(0.5,1.3)(0.5,0.5) \psline(0.5,0.5)(1,0.5)
\psline(1,0.5)(1,1.1)
\psline[linewidth=1.5pt,linecolor=pup](1,1.1)(1,1.7)
\psline(0.5,0.9)(1,0.9)
\psline[linewidth=1.5pt,linecolor=pup](0.5,1.3)(1,1.3)
\psline[linewidth=1.5pt,linecolor=pup](1,1.1)(1.5,1.1)
\psline[linewidth=1.5pt,linecolor=pup](1.5,1.7)(1.5,0.0)
\psline[linewidth=1.5pt,linecolor=pup](2,0.0)(2,1.7)
\psline(0,-.7)(0,-.2) \psline(1.5,-.7)(1.5,-.2)
\psline(2,-0.7)(2,-0.2)
\endpspicture ,
\pspicture[.42](-.4,-.6)(2.4,1.8)
\psframe[linecolor=darkred](-.2,-.2)(2.2,0)
\psline[linewidth=1.5pt,linecolor=pup](0,0)(0,1.7)
\psline(0,1.1)(0.5,1.1)
\psline[linewidth=1.5pt,linecolor=pup](0.5,1.7)(0.5,0.9)
\psline(0.5,0.9)(0.5,0.5) \psline(0.5,0.5)(1,0.5)
\psline(1,0.5)(1,0.9)
\psline[linewidth=1.5pt,linecolor=pup](1,0.9)(1,1.7)
\psline(0.5,1.3)(1,1.3)
\psline[linewidth=1.5pt,linecolor=pup](0.5,0.9)(1,0.9)
\psline[linewidth=1.5pt,linecolor=pup](1,1.1)(1.5,1.1)
\psline[linewidth=1.5pt,linecolor=pup](1.5,1.7)(1.5,0.0)
\psline[linewidth=1.5pt,linecolor=pup](2,0.0)(2,1.7)
\psline(0,-.7)(0,-.2) \psline(1.5,-.7)(1.5,-.2)
\psline(2,-0.7)(2,-0.2)
\endpspicture \\
\pspicture[.42](-.4,-.6)(2.4,1.8)
\psframe[linecolor=darkred](-.2,-.2)(2.2,0)
\psline[linewidth=1.5pt,linecolor=pup](0,0)(0,1.7)
\psline[linewidth=1.5pt,linecolor=pup](0,1.1)(0.5,1.1)
\psline(0.5,1.1)(0.5,0.5)
\psline[linewidth=1.5pt,linecolor=pup](0.5,1.1)(0.5,1.7)
\psline(0.5,0.5)(1,0.5) \psline(1,0.5)(1,1.1)
\psline[linewidth=1.5pt,linecolor=pup](1,1.1)(1,1.7)
\psline(0.5,0.9)(1,0.9) \psline(0.5,1.3)(1,1.3)
\psline[linewidth=1.5pt,linecolor=pup](1,1.1)(1.5,1.1)
\psline[linewidth=1.5pt,linecolor=pup](1.5,1.7)(1.5,0.0)
\psline[linewidth=1.5pt,linecolor=pup](2,0.0)(2,1.7)
\psline(0,-.7)(0,-.2) \psline(1.5,-.7)(1.5,-.2)
\psline(2,-0.7)(2,-0.2)
\endpspicture,
\pspicture[.42](-.4,-.6)(2.4,1.8)
\psframe[linecolor=darkred](-.2,-.2)(2.2,0)
\psline[linewidth=1.5pt,linecolor=pup](0,0)(0,1.7)
\psline(0,1.1)(0.5,1.1) \psline(0.5,1.3)(0.5,0.5)
\psline(0.5,0.5)(1,0.5) \psline(1,0.5)(1,1.3)
\psline(0.5,0.9)(1,0.9)\psline[linewidth=1.5pt,linecolor=pup](1,1.3)(1,1.7)
\psline[linewidth=1.5pt,linecolor=pup](0.5,1.3)(0.5,1.7)
\psline[linewidth=1.5pt,linecolor=pup](0.5,1.3)(1,1.3)
\psline(1,1.1)(1.5,1.1)
\psline[linewidth=1.5pt,linecolor=pup](1.5,1.7)(1.5,0.0)
\psline[linewidth=1.5pt,linecolor=pup](2,0.0)(2,1.7)
\psline(0,-.7)(0,-.2) \psline(1.5,-.7)(1.5,-.2)
\psline(2,-0.7)(2,-0.2)
\endpspicture,
\pspicture[.42](-.4,-.6)(2.4,1.8)
\psframe[linecolor=darkred](-.2,-.2)(2.2,0)
\psline[linewidth=1.5pt,linecolor=pup](0,0)(0,1.7)
\psline(0,1.1)(0.5,1.1) \psline(0.5,1.3)(0.5,0.5)
\psline(0.5,0.5)(1,0.5) \psline(1,0.5)(1,0.9)
\psline(0.5,1.3)(1,1.3)\psline[linewidth=1.5pt,linecolor=pup](1,0.9)(1,1.7)
\psline[linewidth=1.5pt,linecolor=pup](0.5,0.9)(0.5,1.7)
\psline[linewidth=1.5pt,linecolor=pup](0.5,0.9)(1,0.9)
\psline(1,1.1)(1.5,1.1)
\psline[linewidth=1.5pt,linecolor=pup](1.5,1.7)(1.5,0.0)
\psline[linewidth=1.5pt,linecolor=pup](2,0.0)(2,1.7)
\psline(0,-.7)(0,-.2) \psline(1.5,-.7)(1.5,-.2)
\psline(2,-0.7)(2,-0.2)
\endpspicture,
\pspicture[.42](-.4,-.6)(2.4,1.8)
\psframe[linecolor=darkred](-.2,-.2)(2.2,0)
\psline[linewidth=1.5pt,linecolor=pup](0,0)(0,1.7)
\psline(0,1.1)(0.5,1.1) \psline(0.5,1.3)(1,1.3)
\psline(0.5,0.9)(1,0.9)\psline[linewidth=1.5pt,linecolor=pup](1,0.5)(1,1.7)
\psline[linewidth=1.5pt,linecolor=pup](0.5,0.5)(0.5,1.7)
\psline[linewidth=1.5pt,linecolor=pup](0.5,0.5)(1,0.5)
\psline(1,1.1)(1.5,1.1)
\psline[linewidth=1.5pt,linecolor=pup](1.5,1.7)(1.5,0.0)
\psline[linewidth=1.5pt,linecolor=pup](2,0.0)(2,1.7)
\psline(0,-.7)(0,-.2) \psline(1.5,-.7)(1.5,-.2)
\psline(2,-0.7)(2,-0.2)
\endpspicture
\end{eqnarray}

If we examine them to determine whether it will appeared in the
actual expansion, the first two in the second row appear but the
rest of them do not. One can see that $U$ turn can appear only
once at the very top.

\begin{eqnarray}\nonumber & \pspicture[.42](-.7,-1.2)(2.2,1.5)
\psline(-.5,-.2)(-.5,1.4)\psline[arrowscale=1.5]{->}(-.5,1.1)(-.5,1.3)
\psline(0,-.2)(0,1.4)\psline[arrowscale=1.5]{->}(0,1.1)(0,1.3)
\psline(.5,-.2)(.5,1.4)\psline[arrowscale=1.5]{<-}(.5,1.1)(.5,1.3)
\psline(1,-.2)(1,1.4)\psline[arrowscale=1.5]{<-}(1,1.1)(1,1.3)
\psline(1.5,-.6)(1.5,1.2)\psline[arrowscale=1.5]{->}(1.5,.5)(1.5,.7)
\psline(1,-.4)(1,-.6)(1.5,-.6) \psline(1.5,1.2)(2,1.2)
\psline(-.5,-.4)(-.5,-1)\psline[arrowscale=1.5]{<-}(-.5,-.6)(-.5,-.8)
\psline(0,-.4)(0,-1)\psline[arrowscale=1.5]{<-}(0,-.6)(0,-.8)
\psline(.5,-.4)(.5,-1)\psline[arrowscale=1.5]{->}(.5,-.6)(.5,-.8)
\psline(2,-.4)(2,1.2)\psline[arrowscale=1.5]{->}(2,.7)(2,.5)
\psline(0,1)(.5,1)\psline[arrowscale=1.5]{->}(.15,1)(.35,1)
\psline(-.5,.6)(0,.6)\psline[arrowscale=1.5]{<-}(-.15,.6)(-.35,.6)
\psline(0,.2)(.5,.2)\psline[arrowscale=1.5]{<-}(.35,.2)(.15,.2)
\psline(.5,.6)(1,.6)\psline[arrowscale=1.5]{<-}(.65,.6)(.85,.6)
\psline[arrowscale=1.5]{<-}(-.5,-.1)(-.5,.1)\psline[arrowscale=1.5]{<-}(0,-.1)(0,.1)
\psline[arrowscale=1.5]{->}(.5,-.1)(.5,.1)\psline[arrowscale=1.5]{->}(1,-.1)(1,.1)
\psline[arrowscale=1.5]{->}(0,.3)(0,.5)\psline[arrowscale=1.5]{<-}(0,.7)(0,.9)
\psline[arrowscale=1.5]{->}(.5,.5)(.5,.3)\psline[arrowscale=1.5]{->}(.5,.7)(.5,.9)
\psframe[linecolor=darkred](-.7,-.4)(1.2,-.2)
\endpspicture = \frac{[3]}{[5]}\left[ -[2]
\pspicture[.42](-.4,-.6)(2.7,2)
\psframe[linecolor=darkred](-.2,-.2)(2.2,0) \psline(0,0)(0,1.7)
\psline(0.5,1.7)(0.5,1.1) \psline(1,1.1)(1,1.7)
\psline(1.5,1.7)(1.5,0.0) \psline(0,1.1)(0.5,1.1)
\psline(1,1.1)(1.5,1.1) \psline(2,0.0)(2,1.7)
\psline(2,1.7)(2.5,1.7) \psline(2.5,1.7)(2.5,-0.7)
\psline(0,-.7)(0,-.2) \psline(1.5,-.7)(1.5,-.2)
\psline(2,-0.7)(2,-0.2)
\endpspicture -[2]
\pspicture[.42](-.4,-.6)(2.7,2)
\psframe[linecolor=darkred](-.2,-.2)(2.2,0) \psline(0,0)(0,1.7)
\psline(0.5,1.7)(0.5,1.3) \psline(1,1.3)(1,1.7)
\psline(1.5,1.7)(1.5,0.0) \psline(0,1.1)(0.5,1.1)
\psline(0.5,1.1)(0.5,0.9) \psline(0.5,0.9)(1,0.9)
\psline(1,0.9)(1,1.1) \psline(1,1.1)(1.5,1.1)
\psline(0.5,1.3)(1,1.3) \psline(2,0.0)(2,1.7)
\psline(2,1.7)(2.5,1.7) \psline(2.5,1.7)(2.5,-0.7)
\psline(0,-.7)(0,-.2) \psline(1.5,-.7)(1.5,-.2)
\psline(2,-0.7)(2,-0.2)
\endpspicture \right] \\ \nonumber
&+ \frac{[4]}{[5]} \left[ -[2] \pspicture[.42](-.4,-.6)(2.7,2)
\psframe[linecolor=darkred](-.2,-.2)(2.2,0) \psline(0,-.7)(0,-.2)
\psline(1,-.7)(1,-.2) \psline(2,-0.7)(2,-0.2) \psline(0,0)(0,1.7)
\psline(0.5,1.7)(0.5,0.5) \psline(1,1.1)(1,1.7)
\psline(1,0.0)(1,0.5) \psline(1.5,1.7)(1.5,1.1)
\psline(2,0.0)(2,1.7) \psline(2.5,1.7)(2.5,-0.7)
\psline(0,1.1)(0.5,1.1) \psline(0.5,0.5)(1,0.5)
\psline(0.5,1.3)(1,1.3) \psline(1,1.1)(1.5,1.1)
\psline(2,1.7)(2.5,1.7)
\endpspicture + \pspicture[.42](-.4,-.6)(2.7,2)
\psframe[linecolor=darkred](-.2,-.2)(2.2,0) \psline(0,-.7)(0,-.2)
\psline(1,-.7)(1,-.2) \psline(2,-0.7)(2,-0.2) \psline(0,0)(0,1.7)
\psline(0.5,1.7)(0.5,1.1) \psline(1,1.1)(1,1.7)
\psline(1,0.0)(1,0.7) \psline(1.5,1.7)(1.5,0.7)
\psline(2,0.0)(2,1.7) \psline(2.5,1.7)(2.5,-0.7)
\psline(0,1.1)(0.5,1.1) \psline(1,1.1)(1.5,1.1)
\psline(1,0.7)(1.5,0.7) \psline(2,1.7)(2.5,1.7)
\endpspicture
+ \pspicture[.42](-.4,-.6)(2.7,2)
\psframe[linecolor=darkred](-.2,-.2)(2.2,0) \psline(0,-.7)(0,-.2)
\psline(1,-.7)(1,-.2) \psline(2,-0.7)(2,-0.2) \psline(0,0)(0,1.7)
\psline(0.5,1.7)(0.5,1.3) \psline(1,1.3)(1,1.7)
\psline(1,0.0)(1,0.7) \psline(1.5,1.7)(1.5,0.7)
\psline(2,0.0)(2,1.7) \psline(2.5,1.7)(2.5,-0.7)
\psline(0,1.1)(0.5,1.1) \psline(0.5,1.1)(0.5,0.9)
\psline(0.5,0.9)(1,0.9) \psline(1,0.9)(1,1.1)
\psline(1,1.1)(1.5,1.1) \psline(0.5,1.3)(1,1.3)
\psline(1,0.7)(1.5,0.7) \psline(2,1.7)(2.5,1.7)
\endpspicture\right]
\\
\nonumber &+ \frac{[5]}{[5]} \left[
\pspicture[.42](-.4,-.6)(2.7,2)
\psframe[linecolor=darkred](-.2,-.2)(2.2,0) \psline(0,-.7)(0,-.2)
\psline(1.5,-.7)(1.5,-.2) \psline(1,-.7)(1,-.2)
\psline(0,0)(0,1.7) \psline(0.5,1.7)(0.5,0.5)
\psline(1,0.0)(1,0.5) \psline(1,0.7)(1,1.7)
\psline(1.5,0.7)(1.5,0.0) \psline(1.5,0.9)(1.5,1.7)
\psline(2,0.9)(2,1.7) \psline(2.5,1.7)(2.5,-0.7)
\psline(0,1.1)(0.5,1.1) \psline(0.5,0.5)(1,0.5)
\psline(0.5,1.3)(1,1.3) \psline(1,0.7)(1.5,0.7)
\psline(1.5,0.9)(2,0.9) \psline(2,1.7)(2.5,1.7)
\endpspicture + \pspicture[.42](-.4,-.6)(2.7,2)
\psframe[linecolor=darkred](-.2,-.2)(2.2,0) \psline(0,-.7)(0,-.2)
\psline(1.5,-.7)(1.5,-.2) \psline(1,-.7)(1,-.2)
\psline(0,0)(0,1.7) \psline(0.5,1.7)(0.5,0.5)
\psline(1,1.1)(1,1.7) \psline(1,0.0)(1,0.5)
\psline(1.5,1.7)(1.5,1.1) \psline(1.5,0.0)(1.5,0.9)
\psline(2,0.9)(2,1.7) \psline(2.5,1.7)(2.5,-0.7)
\psline(0,1.1)(0.5,1.1) \psline(0.5,0.5)(1,0.5)
\psline(0.5,1.3)(1,1.3) \psline(1,1.1)(1.5,1.1)
\psline(1.5,0.9)(2,0.9) \psline(2,1.7)(2.5,1.7)
\endpspicture\right]
\\
\nonumber & + \frac{[1][2]}{[3][5]} \left[
\pspicture[.42](-.4,-.6)(2.7,2)
\psframe[linecolor=darkred](-.2,-.2)(2.2,0)
\psline(0.5,-0.7)(0.5,-0.2) \psline(1.5,-.7)(1.5,-.2)
\psline(2,-0.7)(2,-0.2) \psline(0,.7)(0,1.7)
\psline(0.5,1.7)(0.5,1.3) \psline(0.5,0.0)(0.5,0.7)
\psline(1,1.3)(1,1.7) \psline(1.5,1.7)(1.5,0.0)
\psline(2,0.0)(2,1.7) \psline(2.5,1.7)(2.5,-0.7)
\psline(0,.7)(0.5,0.7) \psline(0,1.1)(0.5,1.1)
\psline(0.5,1.1)(0.5,0.9) \psline(0.5,0.9)(1,0.9)
\psline(1,0.9)(1,1.1) \psline(1,1.1)(1.5,1.1)
\psline(0.5,1.3)(1,1.3) \psline(2,1.7)(2.5,1.7)
\endpspicture
+ \pspicture[.42](-.4,-.6)(2.7,2)
\psframe[linecolor=darkred](-.2,-.2)(2.2,0)
\psline(0.5,-0.7)(0.5,-0.2) \psline(1.5,-.7)(1.5,-.2)
\psline(2,-0.7)(2,-0.2) \psline(0,.7)(0,1.7)
\psline(0.5,0.0)(0.5,0.7) \psline(0.5,1.7)(0.5,1.1)
\psline(1,1.1)(1,1.7) \psline(1.5,1.7)(1.5,0.0)
\psline(2,0.0)(2,1.7) \psline(2.5,1.7)(2.5,-0.7)
\psline(0,.7)(0.5,0.7) \psline(0,1.1)(0.5,1.1)
\psline(1,1.1)(1.5,1.1)
\endpspicture + \pspicture[.42](-.4,-.6)(2.7,2) \
\psframe[linecolor=darkred](-.2,-.2)(2.2,0)
\psline(0.5,-0.7)(0.5,-0.2) \psline(1.5,-.7)(1.5,-.2)
\psline(2,-0.7)(2,-0.2) \psline(0,.7)(0,1.7)
\psline(0.5,1.7)(0.5,0.7) \psline(0.5,0.0)(0.5,0.5)
\psline(1,0.5)(1,1.7) \psline(1.5,1.7)(1.5,0.0)
\psline(2,0.0)(2,1.7) \psline(2.5,1.7)(2.5,-0.7)
\psline(0,1.1)(0.5,1.1) \psline(0,.7)(0.5,0.7)
\psline(0.5,0.5)(1,0.5) \psline(0.5,1.3)(1,1.3)
\psline(1,1.1)(1.5,1.1) \psline(2,1.7)(2.5,1.7)
\endpspicture\right]
\\ \nonumber &+ \frac{[1][2]}{[4][5]}\left[
\pspicture[.42](-.4,-.6)(2.7,2)
\psframe[linecolor=darkred](-.2,-.2)(2.2,0)
\psline(0.5,-0.7)(0.5,-0.2) \psline(1,-.7)(1,-.2)
\psline(2,-0.7)(2,-0.2) \psline(0,1.1)(0,1.7)
\psline(0.5,0.9)(0.5,0.0) \psline(0.5,1.1)(0.5,1.7)
\psline(1,0.0)(1,1.7) \psline(1.5,1.7)(1.5,0.7)
\psline(2,0.0)(2,1.7) \psline(2.5,1.7)(2.5,-0.7)
\psline(0,1.1)(0.5,1.1) \psline(0.5,0.5)(1,0.5)
\psline(0.5,0.9)(1,0.9) \psline(0.5,1.3)(1,1.3)
\psline(1,1.1)(1.5,1.1) \psline(1,0.7)(1.5,0.7)
\psline(2,1.7)(2.5,1.7)
\endpspicture
+ \pspicture[.42](-.4,-.6)(2.7,2)
\psframe[linecolor=darkred](-.2,-.2)(2.2,0)
\psline(0.5,-0.7)(0.5,-0.2) \psline(1,-.7)(1,-.2)
\psline(2,-0.7)(2,-0.2) \psline(0,.7)(0,1.7)
\psline(0.5,0.7)(0.5,0.0) \psline(0.5,0.9)(0.5,1.7)
\psline(1,0.0)(1,0.9) \psline(1,1.1)(1,1.7)
\psline(1.5,1.7)(1.5,1.1) \psline(2,0.0)(2,1.7)
\psline(2.5,1.7)(2.5,-0.7) \psline(0,.7)(0.5,0.7)
\psline(0.5,0.5)(1,0.5) \psline(0.5,0.9)(1,0.9)
\psline(0.5,1.3)(1,1.3) \psline(1,1.1)(1.5,1.1)
\psline(2,1.7)(2.5,1.7)
\endpspicture + \pspicture[.42](-.4,-.6)(2.7,2)
\psframe[linecolor=darkred](-.2,-.2)(2.2,0)
\psline(0.5,-0.7)(0.5,-0.2) \psline(1,-.7)(1,-.2)
\psline(2,-0.7)(2,-0.2) \psline(0,.7)(0,1.7)
\psline(0.5,0.7)(0.5,0.0) \psline(0.5,0.9)(0.5,1.7)
\psline(1,0.0)(1,0.7) \psline(1,0.9)(1,1.7)
\psline(1.5,1.7)(1.5,0.7) \psline(2,0.0)(2,1.7)
\psline(2.5,1.7)(2.5,-0.7) \psline(0,.7)(0.5,0.7)
\psline(0.5,0.5)(1,0.5) \psline(0.5,0.9)(1,0.9)
\psline(0.5,1.3)(1,1.3) \psline(1,0.7)(1.5,0.7)
\psline(2,1.7)(2.5,1.7)
\endpspicture\right]
\\
&+ \frac{[1][2]}{[5][5]} \left[ \pspicture[.42](-.4,-.6)(2.7,2)
\psframe[linecolor=darkred](-.2,-.2)(2.2,0)
\psline(0.5,-0.7)(0.5,-0.2) \psline(1,-.7)(1,-.2)
\psline(1.5,-.7)(1.5,-.2) \psline(0,.7)(0,1.7)
\psline(0.5,0.0)(0.5,0.7) \psline(0.5,1.7)(0.5,0.9)
\psline(1,0.0)(1,0.9) \psline(1,1.1)(1,1.7)
\psline(1.5,1.7)(1.5,0.0) \psline(2,0.9)(2,1.7)
\psline(2.5,1.7)(2.5,-0.7) \psline(0,.7)(0.5,0.7)
\psline(0.5,0.5)(1,0.5) \psline(0.5,0.9)(1,0.9)
\psline(1,1.1)(1.5,1.1) \psline(1,0.7)(1.5,0.7)
\psline(1.5,0.9)(2,0.9) \psline(2,1.7)(2.5,1.7)
\endpspicture + \pspicture[.42](-.4,-.6)(2.7,2)
\psframe[linecolor=darkred](-.2,-.2)(2.2,0)
\psline(0.5,-0.7)(0.5,-0.2) \psline(1,-.7)(1,-.2)
\psline(1.5,-.7)(1.5,-.2) \psline(0,.7)(0,1.7)
\psline(0.5,0.7)(0.5,0.0) \psline(0.5,1.3)(0.5,1.7)
\psline(1,0.0)(1,1.1) \psline(1,1.3)(1,1.7)
\psline(1.5,0.9)(1.5,0.0) \psline(1.5,1.1)(1.5,1.7)
\psline(2,0.9)(2,1.7) \psline(2.5,1.7)(2.5,-0.7)
\psline(0,.7)(0.5,0.7) \psline(0.5,0.5)(1,0.5)
\psline(0.5,1.3)(1,1.3) \psline(1,1.1)(1.5,1.1)
\psline(1.5,0.9)(2,0.9)
\psline(2,1.7)(2.5,1.7)\endpspicture\right] \label{a2doubleex2}
\end{eqnarray}

Thus, we get a single clasp expansion of a clasp of weight $(2,2)$
as follows.

\begin{eqnarray}\nonumber
\pspicture[.42](-.8,-1.2)(1.3,1.2)
\psline(-.5,-.2)(-.5,1)\psline[arrowscale=1.5]{->}(-.5,.2)(-.5,.4)
\psline(0,-.2)(0,1)\psline[arrowscale=1.5]{->}(0,.2)(0,.4)
\psline(.5,-.2)(.5,1)\psline[arrowscale=1.5]{<-}(.5,.2)(.5,.4)
\psline(1,-.2)(1,1)\psline[arrowscale=1.5]{<-}(1,.2)(1,.4)
\psline(-.5,-.4)(-.5,-1)\psline[arrowscale=1.5]{<-}(-.5,-.6)(-.5,-.8)
\psline(0,-.4)(0,-1)\psline[arrowscale=1.5]{<-}(0,-.6)(0,-.8)
\psline(.5,-.4)(.5,-1)\psline[arrowscale=1.5]{->}(.5,-.6)(.5,-.8)
\psline(1,-.4)(1,-1)\psline[arrowscale=1.5]{->}(1,-.6)(1,-.8)
\psframe[linecolor=darkred](-.7,-.4)(1.2,-.2)
\endpspicture
&= \pspicture[.42](-.8,-1.2)(1.3,1.2) \psline(-.5,-.5)(-.5,1)
\psline(0,-.5)(0,1) \psline(.5,-.5)(.5,1) \psline(-.5,-.7)(-.5,-1)
\psline(0,-.7)(0,-1) \psline(.5,-.7)(.5,-1) \psline(1,1)(1,-1)
\psframe[linecolor=darkred](-.7,-.7)(.7,-.5)
\endpspicture
+ \frac{[1]}{[2]} \pspicture[.42](-.8,-1.2)(1.3,1.2)
\psline(-.5,-.5)(-.5,1) \psline(0,-.5)(0,1)
\psline(.5,-.5)(.5,-.2)(.75,0.1)(.75,.4)(.5,.7)(.5,1)
\psline(1,-1)(1,-.2)(.75,0.1) \psline(.75,.4)(1,.7)(1,1)
\psline(-.5,-.7)(-.5,-1) \psline(0,-.7)(0,-1)
\psline(.5,-.7)(.5,-1)
\psframe[linecolor=darkred](-.7,-.7)(.7,-.5)
\endpspicture
-\frac{[1]}{[2][5]} \pspicture[.42](-.8,-1.2)(1.3,1.2)
\psline(-.5,-.5)(-.5,-.2)(-.25,0.1)(-.25,.4)(-.5,.7)(-.5,1)
\psline(0,-.5)(0,-.2)(.25,0.1)(.75,.4)(.5,.7)(.5,1)
\psline(0,-.2)(-.25,0.1)(-.25,.4)(0,.7)(0,1)
\psline(.25,.1)(.5,-.2)(.5,-.5)(.5,-.2)(.75,0.1)(1,-.2)(1,-1)
\psline(.75,.4)(1,.7)(1,1) \psline(-.5,-.7)(-.5,-1)
\psline(0,-.7)(0,-1) \psline(.5,-.7)(.5,-1)
\psframe[linecolor=darkred](-.7,-.7)(.7,-.5)
\endpspicture\\
&-\frac{[2]}{[5]} \pspicture[.42](-.8,-1.2)(1.3,1.2)
\psline(-.5,-.5)(-.5,1) \psline(0,-.5)(0,-.2)(.25,0.1)(1,.7)(1,1)
\psline(0,1)(0,0.7)(.25,.4)(.5,.7)(.5,1)
\psline(.5,-.5)(.5,-.2)(.25,.1)(.5,-.2)(.75,.1)(1,-.2)(1,-1)
\psline(-.5,-.7)(-.5,-1) \psline(0,-.7)(0,-1)
\psline(.5,-.7)(.5,-1)
\psframe[linecolor=darkred](-.7,-.7)(.7,-.5)
\endpspicture -\frac{[1]}{[5]} \pspicture[.42](-.8,-1.2)(1.3,1.2)
\psline(-.5,-.5)(-.5,1)
\psline(0,-.5)(0,-.2)(.25,0.1)(.25,.4)(0,.7)(0,1)
\psline(.5,-.5)(.5,-.2)(.25,.1)(.25,.4)(.5,.7)(.5,1)
\psline(.5,.7)(.75,.4)(1,.7)(1,1)
\psline(.5,-.2)(.75,.1)(1,-.2)(1,-1) \psline(-.5,-.7)(-.5,-1)
\psline(0,-.7)(0,-1) \psline(.5,-.7)(.5,-1)
\psframe[linecolor=darkred](-.7,-.7)(.7,-.5)
\endpspicture
-\frac{[1]}{[5]} \pspicture[.42](-.8,-1.2)(1.3,1.2)
\psline(-.5,-.5)(-.5,-.2)(-.25,.1)(-.25,.4)(-.5,.7)(-.5,1)
\psline(-.25,.4)(0,.7)(0,1) \psline(-.25,.1)(.5,.7)(.5,1)
\psline(0,-.5)(0,-.2)(.25,.1)(1,.7)(1,1)
\psline(.25,.1)(.5,-.2)(.75,.1)(1,-.2)(1,-1)
\psline(-.5,-.7)(-.5,-1) \psline(0,-.7)(0,-1)
\psline(.5,-.7)(.5,-1)\psline(.5,-.5)(.5,-.2)
\psframe[linecolor=darkred](-.7,-.7)(.7,-.5)
\endpspicture
\label{a2single22}
\end{eqnarray}

By attaching $(2,1)$ clasps on the left top of every web in the
right side of equation~\ref{a2single22}, we get the following
double claps expansion of the clasp of weight $(2,2)$.

\begin{eqnarray}
\pspicture[.42](-.8,-1.2)(1.3,1.2)
\psline(-.5,-.2)(-.5,1)\psline[arrowscale=1.5]{->}(-.5,.2)(-.5,.4)
\psline(0,-.2)(0,1)\psline[arrowscale=1.5]{->}(0,.2)(0,.4)
\psline(.5,-.2)(.5,1)\psline[arrowscale=1.5]{<-}(.5,.2)(.5,.4)
\psline(1,-.2)(1,1)\psline[arrowscale=1.5]{<-}(1,.2)(1,.4)
\psline(-.5,-.4)(-.5,-1)\psline[arrowscale=1.5]{<-}(-.5,-.6)(-.5,-.8)
\psline(0,-.4)(0,-1)\psline[arrowscale=1.5]{<-}(0,-.6)(0,-.8)
\psline(.5,-.4)(.5,-1)\psline[arrowscale=1.5]{->}(.5,-.6)(.5,-.8)
\psline(1,-.4)(1,-1)\psline[arrowscale=1.5]{->}(1,-.6)(1,-.8)
\psframe[linecolor=darkred](-.7,-.4)(1.2,-.2)
\endpspicture
&= \pspicture[.42](-.8,-1.2)(1.3,1.7) \psline(-.5,-.5)(-.5,1.5)
\psline(0,-.5)(0,1.5) \psline(.5,-.5)(.5,1.5)
\psline(-.5,-.7)(-.5,-1) \psline(0,-.7)(0,-1)
\psline(.5,-.7)(.5,-1) \psline(1,1.5)(1,-1)
\psframe[linecolor=darkred](-.7,-.7)(.7,-.5)
\endpspicture
+\frac{[1]}{[2]} \pspicture[.42](-.8,-1.2)(1.3,1.7)
\psline(-.5,-.5)(-.5,1) \psline(0,-.5)(0,1)
\psline(.5,-.5)(.5,-.2)(.75,0.1)(.75,.4)(.5,.7)(.5,1)
\psline(1,-1)(1,-.2)(.75,0.1) \psline(.75,.4)(1,.7)(1,1.5)
\psline(-.5,-.7)(-.5,-1) \psline(0,-.7)(0,-1)
\psline(.5,-.7)(.5,-1) \psline(-.5,1.5)(-.5,1.2)
\psline(0,1.5)(0,1.2) \psline(.5,1.5)(.5,1.2)
\psframe[linecolor=darkred](-.7,-.7)(.7,-.5)
\psframe[linecolor=darkred](-.7,1)(.7,1.2)
\endpspicture  -\frac{[1]}{[5]} \pspicture[.42](-.8,-1.2)(1.3,1.7)
\psline(-.5,-.5)(-.5,1)
\psline(0,-.5)(0,-.2)(.25,0.1)(.25,.4)(0,.7)(0,1)
\psline(.5,-.5)(.5,-.2)(.25,.1)(.25,.4)(.5,.7)(.5,1)
\psline(.5,.7)(.75,.4)(1,.7)(1,1.5) \psline(-.5,1.5)(-.5,1.2)
\psline(0,1.5)(0,1.2) \psline(.5,1.5)(.5,1.2)
\psline(.5,-.2)(.75,.1)(1,-.2)(1,-1) \psline(-.5,-.7)(-.5,-1)
\psline(0,-.7)(0,-1) \psline(.5,-.7)(.5,-1)
\psframe[linecolor=darkred](-.7,-.7)(.7,-.5)
\psframe[linecolor=darkred](-.7,1)(.7,1.2)
\endpspicture
\label{a2double22}
\end{eqnarray}

Unfortunately, there is no particular order we can put for these
basis webs for single clasp expansions. But for the double clasp
expansion, we can generalize the euqation~\ref{a2double22} as
follow. In equation~\ref{a2ab2exp}, the green box between two
clasps is the unique cut out from the hexagonal tiling with the
given boundary as we have seen in Figure~\ref{hexaex}. For
equation~\ref{a2ab2exp} we assume that $a\ge b\ge 1$.

\begin{eqnarray}
\pspicture[.4](-.5,-1)(1.2,2.5) \rput[t](0,-.5){$a$}
\rput[b](0,1.9){$a$} \rput[b](.7,1.9){$b$} \rput[t](.7,-.5){$b$}
\psline(0,.6)(0,1.8)\psline[arrowscale=1.5]{->}(0,.9)(0,1.1)
\psline(0,-.4)(0,.4)\psline[arrowscale=1.5]{<-}(0,.1)(0,-.1)
\psline(.7,1.8)(.7,.6)\psline[arrowscale=1.5]{<-}(.7,.9)(.7,1.1)
\psline(.7,.4)(.7,-.4)\psline[arrowscale=1.5]{->}(.7,.1)(.7,-.1)
\psframe[linecolor=darkred](-.2,.4)(.9,.6)
\endpspicture
= \pspicture[.4](-.5,-1)(2.4,2.5) \rput[t](0,-.5){$a$}
\rput[b](0,1.9){$a$} \rput[b](1.2,1.9){$b-1$}
\rput[t](1.2,-.5){$b-1$}
\psline(0,.6)(0,1.8)\psline[arrowscale=1.5]{->}(0,.9)(0,1.1)
\psline(0,-.4)(0,.4)\psline[arrowscale=1.5]{<-}(0,.1)(0,-.1)
\psline(1.2,1.8)(1.2,.6)\psline[arrowscale=1.5]{<-}(1.2,.9)(1.2,1.1)
\psline(1.2,.4)(1.2,-.4)\psline[arrowscale=1.5]{->}(1.2,.1)(1.2,-.1)
\psline(2.2,1.8)(2.2,-.4)\psline[arrowscale=1.5]{->}(2.2,.1)(2.2,-.1)
\psframe[linecolor=darkred](-.2,.4)(1.9,.6)
\endpspicture
+\alpha \pspicture[.4](-.8,-1)(2.5,2.5) \rput[b](0,1.9){$a$}
\rput[bl](.4,1.9){$b-1$}\rput[t](0,-.5){$a$}\rput[tl](.4,-.5){$b-1$}
\rput[bl](0.5,.4){$b-2$}\rput[br](-.1,.4){$a$}
\psline(0,-.4)(0,0)\psline[arrowscale=1.5]{->}(0,-.3)(0,-.1)
\psline(0,.2)(0,1)\psline[arrowscale=1.5]{->}(0,.5)(0,.7)
\psline(0,1.2)(0,1.8)\psline[arrowscale=1.5]{->}(0,1.4)(0,1.6)
\psline(.4,0)(.4,-0.4)\psline[arrowscale=1.5]{->}(.4,-.1)(.4,-.3)
\psline(.4,1)(.4,.2)\psline[arrowscale=1.5]{->}(.4,.6)(.4,.4)
\psline(.4,1.8)(.4,1.2)\psline[arrowscale=1.5]{->}(.4,1.6)(.4,1.4)
\psline(2.25,1.8)(2.25,1)\psline[arrowscale=1.5]{->}(2.25,1.6)(2.25,1.4)
\psline(2.25,.2)(2.25,-.4)\psline[arrowscale=1.5]{->}(2.25,0)(2.25,-.2)
\psline(2,.45)(2,.75)\psline[arrowscale=1.5]{->}(2,.5)(2,.7)
\pccurve[angleA=270,angleB=180,ncurv=1](1.75,1)(2,.75)
\pccurve[angleA=0,angleB=270,ncurv=1](2,.75)(2.25,1)
\pccurve[angleA=90,angleB=180,ncurv=1](1.75,.2)(2,.45)
\pccurve[angleA=0,angleB=90,ncurv=1](2,.45)(2.25,.2)
\psframe[linecolor=darkred](-.2,0)(1.9,.2)
\psframe[linecolor=darkred](-.2,1)(1.9,1.2)\endpspicture +\beta
\pspicture[.4](-1.3,-1)(2.7,2.5) \rput[br](0,1.9){$a-1$}
\rput[br](1.75,1.9){$b-1$} \rput[tr](0,-.5){$a-1$}
\rput[tr](1.75,-.5){$b-1$} \rput[b](0.5,1.9){$1$}
\rput[b](2.25,1.9){$1$} \rput[t](0.5,-.5){$1$}
\rput[t](2.25,-.5){$1$} \rput[br](-.1,.4){$a-1$}
\psline(0,-.5)(0,-.1)\psline[arrowscale=1.5]{->}(0,-.4)(0,-.4)
\psline(0,.1)(0,1.1)\psline[arrowscale=1.5]{->}(0,.4)(0,.6)
\psline(0,1.3)(0,1.8)\psline[arrowscale=1.5]{->}(0,1.45)(0,1.65)
\psline(.5,-.5)(.5,-.1)\psline[arrowscale=1.5]{->}(.5,-.4)(.5,-.2)
\psline(.5,1.3)(.5,1.8)\psline[arrowscale=1.5]{->}(.5,1.45)(.5,1.65)
\psline(1.5,-.1)(1.5,-0.5)\psline[arrowscale=1.5]{->}(1.5,-.2)(1.5,-.4)
\psline(1.5,1.8)(1.5,1.3)\psline[arrowscale=1.5]{->}(1.5,1.65)(1.5,1.45)
\psline(2.25,1.8)(2.25,1.1) \psline(2.25,.1)(2.25,-.5)
\pcline(.65,.75)(.5,1.1)\middlearrow
\pcline(.5,.1)(.65,.45)\middlearrow
\pcline(1.5,1.1)(1.25,.75)\middlearrow
\pcline(1.25,.45)(1.5,.1)\middlearrow
\pcline(2.25,1.1)(1.75,.65)\middlearrow
\pcline(1.75,.55)(2.25,.1)\middlearrow
\psline[linecolor=emgreen](.5,.6)(.75,.85)
\psline[linecolor=emgreen](.5,.6)(.75,.35)
\psline[linecolor=emgreen](2,.6)(.75,.85)
\psline[linecolor=emgreen](2,.6)(.75,.35)
\psframe[linecolor=darkred](-.2,-.1)(1.9,.1)
\psframe[linecolor=darkred](-.2,1.1)(1.9,1.3)\endpspicture
\label{a2ab2exp}
\end{eqnarray}

\begin{thm}
The coefficients in equation~\ref{a2ab2exp} are
$\alpha=\frac{[b-1]}{[b]}$, $\beta = -\frac{[a]}{[b][a+b+1]}$. We
assume that $[0]=0$ for $\alpha$. \label{a2ab2expthm}
\end{thm}
\begin{proof}
It follows from lemma~\ref{a2abexplem2} and
lemma~\ref{a2abexplem3} that $\alpha=a_{2,a}$ and
$\beta=a_{2,a-1}-[2]a_{1,a-1}+a_{1,a}$.
\end{proof}

To prove two key lemmas, we generalize the idea of paths in the
first example. First, we attach $H$'s as in figure~\ref{a2abexp4}
to all basis webs in equation~\ref{a2abexp} each of basis web is
denoted by $D_{ij}$. After attaching $H$'s as in
figure~\ref{a2abexp4}, the resulting web is denoted by $\tilde
D_{ij}$.

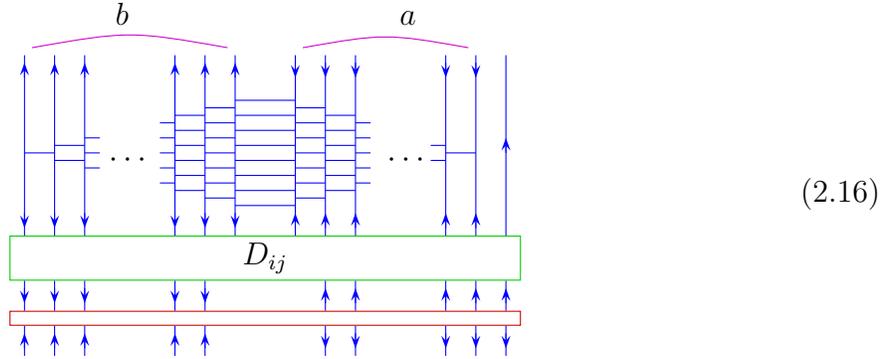
\begin{eqnarray}
\pspicture[.45](-.5,-1.2)(7.5,3.8)
\pccurve[angleA=10,angleB=170,ncurv=1,linecolor=pup](.2,3.1)(2.8,3.1)
\pccurve[angleA=10,angleB=170,ncurv=1,linecolor=pup](3.8,3.1)(6,3.1)
\rput[b](3.3,0.1){$D_{ij}$} \rput[b](1.4,3.4){$b$}
\rput[b](5.2,3.4){$a$}
\psline(.1,.6)(.1,3)\psline[arrowscale=1.5]{->}(.1,.9)(.1,.7)
\psline[arrowscale=1.5]{->}(.1,2.7)(.1,2.9)
\psline(.1,0)(.1,-.4)\psline[arrowscale=1.5]{->}(.1,-.1)(.1,-.3)
\psline(.1,-.6)(.1,-1)\psline[arrowscale=1.5]{->}(.1,-.9)(.1,-.7)
\psline(.5,.6)(.5,3)\psline[arrowscale=1.5]{->}(.5,.9)(.5,.7)
\psline[arrowscale=1.5]{->}(.5,2.7)(.5,2.9)
\psline(.5,0)(.5,-.4)\psline[arrowscale=1.5]{->}(.5,-.1)(.5,-.3)
\psline(.5,-.6)(.5,-1)\psline[arrowscale=1.5]{->}(.5,-.9)(.5,-.7)
\psline(.9,.6)(.9,3)\psline[arrowscale=1.5]{->}(.9,.9)(.9,.7)
\psline[arrowscale=1.5]{->}(.9,2.7)(.9,2.9)
\psline(.9,0)(.9,-.4)\psline[arrowscale=1.5]{->}(.9,-.1)(.9,-.3)
\psline(.9,-.6)(.9,-1)\psline[arrowscale=1.5]{->}(.9,-.9)(.9,-.7)
\psline(2.1,.6)(2.1,3)\psline[arrowscale=1.5]{->}(2.1,.9)(2.1,.7)
\psline[arrowscale=1.5]{->}(2.1,2.7)(2.1,2.9)
\psline(2.1,0)(2.1,-.4)\psline[arrowscale=1.5]{->}(2.1,-.1)(2.1,-.3)
\psline(2.1,-.6)(2.1,-1)\psline[arrowscale=1.5]{->}(2.1,-.9)(2.1,-.7)
\psline(2.5,.6)(2.5,3)\psline[arrowscale=1.5]{->}(2.5,.9)(2.5,.7)
\psline[arrowscale=1.5]{->}(2.5,2.7)(2.5,2.9)
\psline(2.5,0)(2.5,-.4)\psline[arrowscale=1.5]{->}(2.5,-.1)(2.5,-.3)
\psline(2.5,-.6)(2.5,-1)\psline[arrowscale=1.5]{->}(2.5,-.9)(2.5,-.7)
\psline(2.9,.6)(2.9,3)\psline[arrowscale=1.5]{->}(2.9,.9)(2.9,.7)
\psline[arrowscale=1.5]{->}(2.9,2.7)(2.9,2.9)
\psline(3.7,.6)(3.7,3)\psline[arrowscale=1.5]{<-}(3.7,.9)(3.7,.7)
\psline[arrowscale=1.5]{<-}(3.7,2.7)(3.7,2.9)
\psline(4.1,.6)(4.1,3)\psline[arrowscale=1.5]{<-}(4.1,.9)(4.1,.7)
\psline[arrowscale=1.5]{<-}(4.1,2.7)(4.1,2.9)
\psline(4.1,0)(4.1,-.4)\psline[arrowscale=1.5]{<-}(4.1,-.1)(4.1,-.3)
\psline(4.1,-.6)(4.1,-1)\psline[arrowscale=1.5]{<-}(4.1,-.9)(4.1,-.7)
\psline(4.5,.6)(4.5,3)\psline[arrowscale=1.5]{<-}(4.5,.9)(4.5,.7)
\psline[arrowscale=1.5]{<-}(4.5,2.7)(4.5,2.9)
\psline(4.5,0)(4.5,-.4)\psline[arrowscale=1.5]{<-}(4.5,-.1)(4.5,-.3)
\psline(4.5,-.6)(4.5,-1)\psline[arrowscale=1.5]{<-}(4.5,-.9)(4.5,-.7)
\psline(5.7,.6)(5.7,3)\psline[arrowscale=1.5]{<-}(5.7,.9)(5.7,.7)
\psline[arrowscale=1.5]{<-}(5.7,2.7)(5.7,2.9)
\psline(5.7,0)(5.7,-.4)\psline[arrowscale=1.5]{<-}(5.7,-.1)(5.7,-.3)
\psline(5.7,-.6)(5.7,-1)\psline[arrowscale=1.5]{<-}(5.7,-.9)(5.7,-.7)
\psline(6.1,.6)(6.1,3)\psline[arrowscale=1.5]{<-}(6.1,.9)(6.1,.7)
\psline[arrowscale=1.5]{<-}(6.1,2.7)(6.1,2.9)
\psline(6.1,0)(6.1,-.4)\psline[arrowscale=1.5]{<-}(6.1,-.1)(6.1,-.3)
\psline(6.1,-.6)(6.1,-1)\psline[arrowscale=1.5]{<-}(6.1,-.9)(6.1,-.7)
\psline(6.5,.6)(6.5,3)\psline[arrowscale=1.5]{->}(6.5,1.7)(6.5,1.9)
\psline(6.5,0)(6.5,-.4)\psline[arrowscale=1.5]{<-}(6.5,-.1)(6.5,-.3)
\psline(6.5,-.6)(6.5,-1)\psline[arrowscale=1.5]{<-}(6.5,-.9)(6.5,-.7)
\qline(.1,1.7)(.5,1.7) \qline(.5,1.8)(.9,1.8)
\qline(.5,1.6)(.9,1.6) \qline(.9,1.9)(1.1,1.9)
\qline(.9,1.7)(1.1,1.7) \qline(.9,1.5)(1.1,1.5)
\qline(1.9,2.1)(2.1,2.1) \qline(1.9,1.9)(2.1,1.9)
\qline(1.9,1.7)(2.1,1.7) \qline(1.9,1.5)(2.1,1.5)
\qline(1.9,1.3)(2.1,1.3) \qline(2.1,2.2)(2.5,2.2)
\qline(2.1,2)(2.5,2) \qline(2.1,1.8)(2.5,1.8)
\qline(2.1,1.6)(2.5,1.6) \qline(2.1,1.4)(2.5,1.4)
\qline(2.1,1.2)(2.5,1.2) \qline(2.5,2.3)(2.9,2.3)
\qline(2.5,2.1)(2.9,2.1) \qline(2.5,1.9)(2.9,1.9)
\qline(2.5,1.7)(2.9,1.7) \qline(2.5,1.5)(2.9,1.5)
\qline(2.5,1.3)(2.9,1.3) \qline(2.5,1.1)(2.9,1.1)
\qline(4.1,2.3)(3.7,2.3) \qline(4.1,2.1)(3.7,2.1)
\qline(4.1,1.9)(3.7,1.9) \qline(4.1,1.7)(3.7,1.7)
\qline(4.1,1.5)(3.7,1.5) \qline(4.1,1.3)(3.7,1.3)
\qline(4.1,1.1)(3.7,1.1) \qline(4.1,2.2)(4.5,2.2)
\qline(4.1,2)(4.5,2) \qline(4.1,1.8)(4.5,1.8)
\qline(4.1,1.6)(4.5,1.6) \qline(4.1,1.4)(4.5,1.4)
\qline(4.1,1.2)(4.5,1.2) \qline(4.5,2.1)(4.7,2.1)
\qline(4.5,1.9)(4.7,1.9) \qline(4.5,1.7)(4.7,1.7)
\qline(4.5,1.5)(4.7,1.5) \qline(4.5,1.3)(4.7,1.3)
\qline(5.7,1.7)(6.1,1.7) \qline(5.5,1.8)(5.7,1.8)
\qline(5.5,1.6)(5.7,1.6) \qline(2.9,1)(3.7,1)
\qline(2.9,1.2)(3.7,1.2) \qline(2.9,1.4)(3.7,1.4)
\qline(2.9,1.6)(3.7,1.6) \qline(2.9,1.8)(3.7,1.8)
\qline(2.9,2)(3.7,2) \qline(2.9,2.2)(3.7,2.2)
\qline(2.9,2.4)(3.7,2.4)\rput[t](1.5,1.7){$\cdots$}\rput[t](5.2,1.7){$\cdots$}
\psframe[linecolor=darkred](-0.1,-.6)(6.7,-0.4)
\psframe[linecolor=emgreen](-0.1,0)(6.7,0.6)
\endpspicture
\label{a2abexp4}
\end{eqnarray}

As we have seen in the example, $\tilde D_{i,j}$ is not a basis
web because it contains some elliptic faces. If we decompose each
$\tilde D_{i,j}$ into a linear combination of some webs which have
no elliptic faces, then the union of all these resulting webs
actually forms a basis. Let us prove that these webs actually form
a basis which will be denoted by $D'_{i',j'}$. As vector spaces,
this change, adding $H$'s, induces an isomorphism between two web
spaces. Its matrix representation with respect to these web bases
$\{D_{i,j}\}$ and $\{ D'_{i',j'}\}$ is an $(a+1)b\times (a+1)b$
matrix whose entries are $0,1$ or $-[2]$. In general, we will not
be able to write this matrix because there are many nonzero
entries in every columns and rows. But we know that the
determinant of this matrix is $\pm [2]^{ab}$ because each one $H$
contributes $\pm [2]$ depending on the directions.

Since $\tilde D_{i,j}$ is not a basis web, to find a single clasp
expansion, we might have to use relations to find its linear
expansion into a new web basis $D'_{i',j'}$. In general this might
not be done. If we just limit ourself to a double clasp expansion,
We could use the idea of paths as we demonstrated in the example.
Let us formally define it, a {\it stem} of a web. Geometrically it
is transversal to cut paths. From $\tilde D_{i,j}$, we see that
there are $a+b+1$ points on top but only $a+b-1$ lines right above
the clasp. Because of one of properties of the clasp of weight
$(a,b-1)$: if we have a cut path of weight which is less than
$(a,b-1)$, then the web becomes zero, we must have $a+b-1$
vertical lines which connect top $a+b+1$ nodes to the clasp of
weight $(a,b-1)$ for non-vanishing webs after applying relations.
It is clear that these connecting lines should be mutually
disjoint, otherwise, we will have a cut path with weight less than
$(a,b-1)$. A {\it stem} of a web is a disjoint union of lines as
we described. Unfortunately some of stems do not arise all cases
because it may not be obtained by removing elliptic faces. If a
stem appears, we call it an {\it admissible stem}. For single
clasp expansion, finding all these stems will be more difficult
than an expansion by relations but for the double clasp expansion
of segregated clasps, there are only few possible admissible stems
whose coefficient is nonzero.

\begin{lem}
After attaching a clasp of weight $(a,b-1)$ to top of webs $\tilde
D_{i,j}$ from equation~\ref{a2abexp4}, the only $3$ non-vanishing
shapes are those in Figure~\ref{a2ab2exppo}.

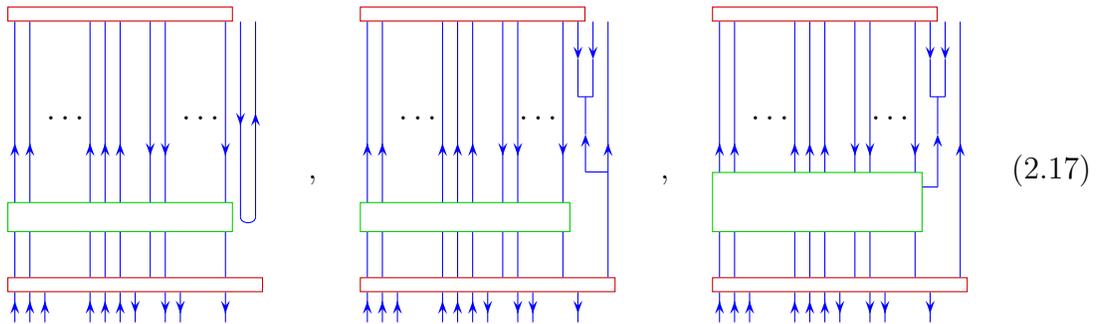
\begin{eqnarray}
\pspicture[.45](-.5,-1.2)(4,3.5) \rput[t](.8,1.8){$\cdots$}
\rput[t](2.6,1.8){$\cdots$} \pcline(.1,-0.4)(.1,3)\middlearrow
\pcline(.3,-0.4)(.3,3)\middlearrow
\pcline(1.1,-0.4)(1.1,3)\middlearrow
\pcline(1.3,-0.4)(1.3,3)\middlearrow
\pcline(1.5,-0.4)(1.5,3)\middlearrow
\pcline(1.9,3)(1.9,-0.4)\middlearrow
\pcline(2.1,3)(2.1,-0.4)\middlearrow
\pcline(2.9,3)(2.9,-0.4)\middlearrow
\pcline(3.1,3)(3.1,0.4)\middlearrow
\pcline(3.3,0.4)(3.3,3)\middlearrow
\pccurve[angleA=-90,angleB=-90,ncurv=1](3.1,0.4)(3.3,0.4)
\pcline(0.1,-1)(0.1,-.6)\middlearrow
\pcline(0.3,-1)(0.3,-.6)\middlearrow
\pcline(0.5,-1)(0.5,-.6)\middlearrow
\pcline(1.1,-1)(1.1,-.6)\middlearrow
\pcline(1.3,-1)(1.3,-.6)\middlearrow
\pcline(1.5,-1)(1.5,-.6)\middlearrow
\pcline(1.7,-.6)(1.7,-1)\middlearrow
\pcline(2.1,-.6)(2.1,-1)\middlearrow
\pcline(2.3,-.6)(2.3,-1)\middlearrow
\pcline(2.9,-.6)(2.9,-1)\middlearrow
\psframe[linecolor=darkred](0,3)(3,3.2)
\psframe[linecolor=darkred](0,-.6)(3.4,-0.4)
\psframe[linecolor=emgreen,fillstyle=solid,fillcolor=white](0,.2)(3,0.6)
\endpspicture
, \pspicture[.45](-.5,-1.2)(4,3.5) \rput[t](.8,1.8){$\cdots$}
\rput[t](2.4,1.8){$\cdots$} \pcline(.1,-0.4)(.1,3)\middlearrow
\pcline(.3,-0.4)(.3,3)\middlearrow
\pcline(1.1,-0.4)(1.1,3)\middlearrow
\pcline(1.3,-0.4)(1.3,3)\middlearrow
\pcline(1.5,-0.4)(1.5,3)\middlearrow
\pcline(1.9,3)(1.9,-0.4)\middlearrow
\pcline(2.1,3)(2.1,-0.4)\middlearrow
\pcline(2.7,3)(2.7,-0.4)\middlearrow
\pcline(3.3,-0.4)(3.3,3)\middlearrow
\pcline(0.1,-1)(0.1,-.6)\middlearrow
\pcline(0.3,-1)(0.3,-.6)\middlearrow
\pcline(0.5,-1)(0.5,-.6)\middlearrow
\pcline(1.1,-1)(1.1,-.6)\middlearrow
\pcline(1.3,-1)(1.3,-.6)\middlearrow
\pcline(1.5,-1)(1.5,-.6)\middlearrow
\pcline(1.7,-.6)(1.7,-1)\middlearrow
\pcline(2.1,-.6)(2.1,-1)\middlearrow
\pcline(2.3,-.6)(2.3,-1)\middlearrow
\pcline(2.9,-.6)(2.9,-1)\middlearrow \qline(2.9,2)(3.1,2)
\psline[arrows=->,arrowscale=1.5](2.9,3)(2.9,2.5)\psline(2.9,2.5)(2.9,2)
\psline[arrows=->,arrowscale=1.5](3.1,3)(3.1,2.5)\psline(3.1,2.5)(3.1,2)
\psline[arrows=->,arrowscale=1.5](3,1)(3,1.5)\psline(3,1.5)(3,2)
\qline(3,1)(3.3,1) \psframe[linecolor=darkred](0,-.6)(3.4,-0.4)
\psframe[linecolor=darkred](0,3)(3,3.2)
\psframe[linecolor=emgreen,fillstyle=solid,fillcolor=white](0,.2)(2.8,0.6)
\endpspicture ,
\pspicture[.45](-.5,-1.2)(4,3.5) \rput[t](.8,1.8){$\cdots$}
\rput[t](2.4,1.8){$\cdots$} \pcline(.1,-0.4)(.1,3)\middlearrow
\pcline(.3,-0.4)(.3,3)\middlearrow
\pcline(1.1,-0.4)(1.1,3)\middlearrow
\pcline(1.3,-0.4)(1.3,3)\middlearrow
\pcline(1.5,-0.4)(1.5,3)\middlearrow
\pcline(1.9,3)(1.9,-0.4)\middlearrow
\pcline(2.1,3)(2.1,-0.4)\middlearrow
\pcline(2.7,3)(2.7,-0.4)\middlearrow
\pcline(3.3,-0.4)(3.3,3)\middlearrow
\pcline(0.1,-1)(0.1,-.6)\middlearrow
\pcline(0.3,-1)(0.3,-.6)\middlearrow
\pcline(0.5,-1)(0.5,-.6)\middlearrow
\pcline(1.1,-1)(1.1,-.6)\middlearrow
\pcline(1.3,-1)(1.3,-.6)\middlearrow
\pcline(1.5,-1)(1.5,-.6)\middlearrow
\pcline(1.7,-.6)(1.7,-1)\middlearrow
\pcline(2.1,-.6)(2.1,-1)\middlearrow
\pcline(2.3,-.6)(2.3,-1)\middlearrow
\pcline(2.9,-.6)(2.9,-1)\middlearrow \qline(2.9,2)(3.1,2)
\psline[arrows=->,arrowscale=1.5](2.9,3)(2.9,2.5)\psline(2.9,2.5)(2.9,2)
\psline[arrows=->,arrowscale=1.5](3.1,3)(3.1,2.5)\psline(3.1,2.5)(3.1,2)
\psline[arrows=->,arrowscale=1.5](3,.8)(3,1.5)\psline(3,1.5)(3,2)
\qline(3,.8)(2.7,.8) \psframe[linecolor=darkred](0,-.6)(3.4,-0.4)
\psframe[linecolor=darkred](0,3)(3,3.2)
\psframe[linecolor=emgreen,fillstyle=solid,fillcolor=white](0,.2)(2.8,1)
\endpspicture\label{a2ab2exppo}\end{eqnarray}
\label{a2abexplem2}
\end{lem}

\begin{proof}
From $\tilde D_{i,j}$ we see that there are $a+b+1$ lines on top
and $a+b-1$ lines right above the clasp. If we repeatedly use the
rectangular relation as in equation~\ref{a2defr23}, we can push up
the $Y$'s so that there are either two $Y$'s or one $U$ shape at
the top. It is possible to have two adjacent $Y$'s which appear in
the second and third figures in Figure~\ref{a2ab2exppo} but a $U$
turn can appear in only two places because of the orientation of
edges. If we attach the $(a,b-1)$ clasp to the top of the
resulting web from the left and $U$ or $Y$ shape appear just below
it, the web becomes zero. Therefore only these three webs do not
vanish.
\end{proof}

For the next lemma, we will find all $\tilde D_{i,j}$'s which can
be transformed to each of the figures in Figure~\ref{a2ab2exppo}.

\begin{lem}
Only $\tilde D_{1,a}(\tilde D_{2,a})$ can be transformed to the
first(second, respectively) shape in Figure~\ref{a2ab2exppo}. Only
the three webs, $\tilde D_{1,a-1}, \tilde D_{1,a}$ and $\tilde
D_{2,a-1}$ can be transformed to the last shape. Moreover, all of
these transformations use only rectangular relations as in
equation~\ref{a2defr23} except for the transformation from $\tilde
D_{1,a-1}$ to the third figure uses one loop relation in
equation~\ref{a2defr22}. \label{a2abexplem3}
\end{lem}
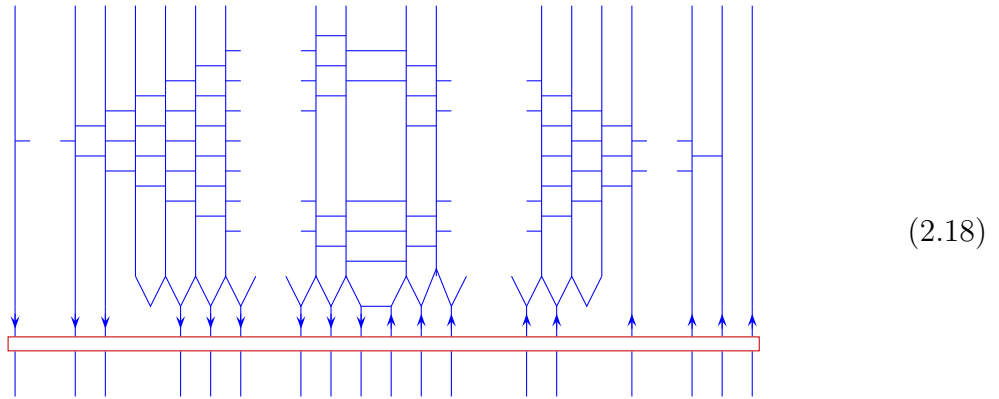
\begin{proof}

\begin{eqnarray}
\pspicture[.4](-.2,-.3)(10.1,5.1)
\qline(0.0,0.6)(0.0,1.0)\psline[arrowscale=1.5]{->}(0.0,.9)(0.0,.7)
\qline(0.8,0.6)(0.8,1.0)\psline[arrowscale=1.5]{->}(0.8,.9)(0.8,.7)
\qline(1.2,0.6)(1.2,1.0)\psline[arrowscale=1.5]{->}(1.2,.9)(1.2,.7)
\qline(2.2,0.6)(2.2,1.0)\psline[arrowscale=1.5]{->}(2.2,.9)(2.2,.7)
\qline(2.6,0.6)(2.6,1.0)\psline[arrowscale=1.5]{->}(2.6,.9)(2.6,.7)
\qline(3.0,0.6)(3.0,1.0)\psline[arrowscale=1.5]{->}(3.0,.9)(3.0,.7)
\qline(3.8,0.6)(3.8,1.0)\psline[arrowscale=1.5]{->}(3.8,.9)(3.8,.7)
\qline(4.2,0.6)(4.2,1.0)\psline[arrowscale=1.5]{->}(4.2,.9)(4.2,.7)
\qline(4.6,0.6)(4.6,1.0)\psline[arrowscale=1.5]{->}(4.6,.9)(4.6,.7)
\qline(5.0,0.6)(5.0,1.0)\psline[arrowscale=1.5]{<-}(5.0,.9)(5.0,.7)
\qline(5.4,0.6)(5.4,1.0)\psline[arrowscale=1.5]{<-}(5.4,.9)(5.4,.7)
\qline(5.8,0.6)(5.8,1.0)\psline[arrowscale=1.5]{<-}(5.8,.9)(5.8,.7)
\qline(6.8,0.6)(6.8,1.0)\psline[arrowscale=1.5]{<-}(6.8,.9)(6.8,.7)
\qline(7.2,0.6)(7.2,1.0)\psline[arrowscale=1.5]{<-}(7.2,.9)(7.2,.7)
\qline(8.2,0.6)(8.2,1.0)\psline[arrowscale=1.5]{<-}(8.2,.9)(8.2,.7)
\qline(9.0,0.6)(9.0,1.0)\psline[arrowscale=1.5]{<-}(9.0,.9)(9.0,.7)
\qline(9.4,0.6)(9.4,1.0)\psline[arrowscale=1.5]{<-}(9.4,.9)(9.4,.7)
\qline(9.8,0.6)(9.8,1.0)\psline[arrowscale=1.5]{<-}(9.8,.9)(9.8,.7)
\qline(0.0,0.4)(0.0,-0.2) \qline(0.8,0.4)(0.8,-0.2)
\qline(1.2,0.4)(1.2,-0.2) \qline(2.2,0.4)(2.2,-0.2)
\qline(2.6,0.4)(2.6,-0.2) \qline(3.0,0.4)(3.0,-0.2)
\qline(3.8,0.4)(3.8,-0.2) \qline(4.2,0.4)(4.2,-0.2)
\qline(4.6,0.4)(4.6,-0.2) \qline(5.0,0.4)(5.0,-0.2)
\qline(5.4,0.4)(5.4,-0.2) \qline(5.8,0.4)(5.8,-0.2)
\qline(6.8,0.4)(6.8,-0.2) \qline(7.2,0.4)(7.2,-0.2)
\qline(8.2,0.4)(8.2,-0.2) \qline(9.0,0.4)(9.0,-0.2)
\qline(9.4,0.4)(9.4,-0.2) \qline(9.8,0.4)(9.8,-0.2)
\qline(0.0,1.0)(0.,5.0)\qline(0.8,1.0)(0.8,5.0)\qline(1.2,1.0)(1.2,5.0)
\qline(1.6,1.4)(1.6,5.0)\qline(2.0,1.4)(2.0,5.0)\qline(2.4,1.4)(2.4,5.0)
\qline(2.8,1.4)(2.8,5.0)\qline(4.0,1.4)(4.0,5.0)\qline(4.4,1.4)(4.4,5.0)
\qline(5.2,1.4)(5.2,5.0)\qline(5.6,1.4)(5.6,5.0)\qline(7.0,1.4)(7.0,5.0)
\qline(7.4,1.4)(7.4,5.0)\qline(7.8,1.4)(7.8,5.0)\qline(8.2,1.0)(8.2,5.0)
\qline(9.0,1.0)(9.0,5.0)\qline(9.4,1.0)(9.4,5.0)\qline(9.8,1.0)(9.8,5.0)
\psline(1.6,1.4)(1.8,1.0)(2.0,1.4)(2.2,1.0)(2.4,1.4)(2.6,1.0)(2.8,1.4)(3.0,1.0)(3.2,1.4)
\psline(3.6,1.4)(3.8,1.0)(4.0,1.4)(4.2,1.0)(4.4,1.4)(4.6,1.0)(5.0,1.0)(5.2,1.4)(5.4,1.0)(5.6,1.5)(5.8,1.0)(6.0,1.4)
\psline(6.6,1.4)(6.8,1.0)(7.0,1.4)(7.2,1.0)(7.4,1.4)(7.6,1.0)(7.8,1.4)
\qline(0,3.2)(.2,3.2)
\qline(0.6,3.2)(0.8,3.2)\qline(0.8,3.0)(1.2,3.0)\qline(0.8,3.4)(1.2,3.4)
\qline(1.2,2.8)(1.6,2.8)\qline(1.2,3.2)(1.6,3.2)\qline(1.2,3.6)(1.6,3.6)
\qline(1.6,2.6)(2.0,2.6)\qline(1.6,3.0)(2.0,3.0)\qline(1.6,3.4)(2.0,3.4)
\qline(1.6,3.8)(2.0,3.8)\qline(2.0,2.4)(2.4,2.4)\qline(2.0,2.8)(2.4,2.8)
\qline(2.0,3.2)(2.4,3.2)\qline(2.0,3.6)(2.4,3.6)\qline(2.0,4.0)(2.4,4.0)
\qline(2.4,2.2)(2.8,2.2)\qline(2.4,2.6)(2.8,2.6)\qline(2.4,3.0)(2.8,3.0)
\qline(2.4,3.4)(2.8,3.4)\qline(2.4,3.8)(2.8,3.8)\qline(2.4,4.2)(2.8,4.2)
\qline(2.8,2.0)(3.0,2.0)\qline(2.8,2.4)(3.0,2.4)\qline(2.8,2.8)(3.0,2.8)
\qline(2.8,3.2)(3.0,3.2)\qline(2.8,3.6)(3.0,3.6)\qline(2.8,4.0)(3.0,4.0)
\qline(2.8,4.4)(3.0,4.4)\qline(3.8,2.0)(4.0,2.0)\qline(3.8,2.4)(4.0,2.4)
\qline(3.8,3.6)(4.0,3.6)\qline(3.8,4.0)(4.0,4.0)\qline(3.8,4.4)(4.0,4.4)
\qline(4.0,1.8)(4.4,1.8)\qline(4.0,2.2)(4.4,2.2)\qline(4.0,3.8)(4.4,3.8)
\qline(4.0,4.2)(4.4,4.2)\qline(4.0,4.6)(4.4,4.6)\qline(4.4,1.6)(5.2,1.6)
\qline(5.2,2.0)(4.4,2.0)\qline(5.2,2.4)(4.4,2.4)\qline(5.2,4.0)(4.4,4.0)
\qline(5.2,4.4)(4.4,4.4)\qline(5.2,1.8)(5.6,1.8)\qline(5.2,2.2)(5.6,2.2)
\qline(5.2,3.8)(5.6,3.8)\qline(5.2,4.2)(5.6,4.2)\qline(5.2,3.4)(5.6,3.4)
\qline(5.6,2.0)(5.8,2.0)\qline(5.6,2.4)(5.8,2.4)\qline(5.6,4.0)(5.8,4.0)
\qline(5.6,3.6)(5.8,3.6)\qline(7.0,2.0)(6.8,2.0)\qline(6.8,2.4)(7.0,2.4)
\qline(7.0,4.0)(6.8,4.0)\qline(6.8,3.6)(7.0,3.6)\qline(7.0,2.2)(7.4,2.2)
\qline(7.0,2.6)(7.4,2.6)\qline(7.0,3.0)(7.4,3.0)\qline(7.0,3.4)(7.4,3.4)
\qline(7.0,3.8)(7.4,3.8)\qline(7.4,2.4)(7.8,2.4)\qline(7.4,2.8)(7.8,2.8)
\qline(7.4,3.2)(7.8,3.2)\qline(7.4,3.6)(7.8,3.6)\qline(7.8,2.6)(8.2,2.6)
\qline(7.8,3.0)(8.2,3.0)\qline(7.8,3.4)(8.2,3.4)\qline(8.2,2.8)(8.4,2.8)
\qline(8.2,3.2)(8.4,3.2)\qline(8.8,2.8)(9.0,2.8)\qline(8.8,3.2)(9.0,3.2)
\qline(9.0,3.0)(9.4,3.0)
\psframe[linecolor=darkred](-.1,.4)(9.9,.6)
\endpspicture
\label{a2abexp5}\end{eqnarray}

For the first shape in figure~\ref{a2ab2exppo}, it is easy to see
that we need to look at $\tilde D_{i,a}$, for $i=1, 2, \ldots, b$,
otherwise the last two strings can not be changed to the shape
with a $U$ turn as shown in figure~\ref{a2abexp5}.

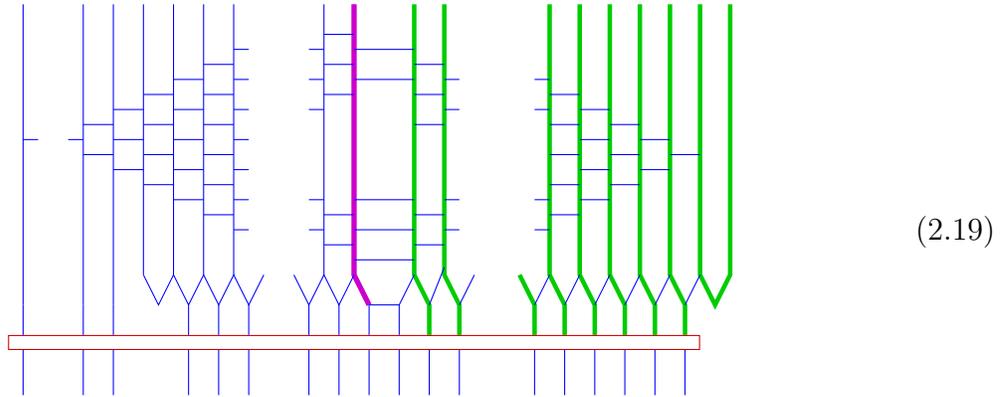
\begin{eqnarray}
\pspicture[.4](-.2,-.3)(10.1,5.1)
\qline(0.0,0.6)(0.0,1.0)\qline(0.8,0.6)(0.8,1.0)\qline(1.2,0.6)(1.2,1.0)
\qline(2.2,0.6)(2.2,1.0)\qline(2.6,0.6)(2.6,1.0)\qline(3.0,0.6)(3.0,1.0)
\qline(3.8,0.6)(3.8,1.0)\qline(4.2,0.6)(4.2,1.0)\qline(4.6,0.6)(4.6,1.0)
\qline(5.0,0.6)(5.0,1.0)\qline(0.0,0.4)(0.0,-0.2)\qline(0.8,0.4)(0.8,-0.2)
\qline(1.2,0.4)(1.2,-0.2)\qline(2.2,0.4)(2.2,-0.2)\qline(2.6,0.4)(2.6,-0.2)
\qline(3.0,0.4)(3.0,-0.2)\qline(3.8,0.4)(3.8,-0.2)\qline(4.2,0.4)(4.2,-0.2)
\qline(4.6,0.4)(4.6,-0.2)\qline(5.0,0.4)(5.0,-0.2)\qline(5.4,0.4)(5.4,-0.2)
\qline(5.8,0.4)(5.8,-0.2)\qline(6.8,0.4)(6.8,-0.2)\qline(7.2,0.4)(7.2,-0.2)
\qline(7.6,0.4)(7.6,-0.2)\qline(8.0,0.4)(8.0,-0.2)\qline(8.4,0.4)(8.4,-0.2)
\qline(8.8,0.4)(8.8,-0.2)\qline(0.0,1.0)(0.,5.0)\qline(0.8,1.0)(0.8,5.0)
\qline(1.2,1.0)(1.2,5.0)\qline(1.6,1.4)(1.6,5.0)\qline(2.0,1.4)(2.0,5.0)
\qline(2.4,1.4)(2.4,5.0)\qline(2.8,1.4)(2.8,5.0)\qline(4.0,1.4)(4.0,5.0)
\qline(4.4,1.4)(4.4,5.0)
\psline[linecolor=emgreen,linewidth=1.7pt](6.6,1.4)(6.8,1)(6.8,.6)
\psline[linecolor=emgreen,linewidth=1.7pt](5.4,.6)(5.4,1)(5.2,1.4)(5.2,5.0)
\psline[linecolor=emgreen,linewidth=1.7pt](5.8,.6)(5.8,1)(5.6,1.4)(5.6,5.0)
\psline[linecolor=emgreen,linewidth=1.7pt](7.2,.6)(7.2,1)(7.0,1.4)(7.0,5.0)
\psline[linecolor=emgreen,linewidth=1.7pt](7.6,.6)(7.6,1)(7.4,1.4)(7.4,5.0)
\psline[linecolor=emgreen,linewidth=1.7pt](8,.6)(8,1)(7.8,1.4)(7.8,5.0)
\psline[linecolor=emgreen,linewidth=1.7pt](8.4,.6)(8.4,1)(8.2,1.4)(8.2,5.0)
\psline[linecolor=emgreen,linewidth=1.7pt](8.8,.6)(8.8,1)(8.6,1.4)(8.6,5.0)
\psline[linecolor=emgreen,linewidth=1.7pt](9.4,5)(9.4,1.4)(9.2,1)(9.0,1.4)(9.0,5.0)
\psline[linecolor=pup,linewidth=2pt](4.4,5)(4.4,1.4)(4.6,1)
\psline(1.6,1.4)(1.8,1)(2,1.4)(2.2,1)(2.4,1.4)(2.6,1)(2.8,1.4)(3,1)(3.2,1.4)
\psline(3.6,1.4)(3.8,1)(4,1.4)(4.2,1)(4.4,1.4)\psline(4.6,1)(5,1)(5.2,1.4)
\qline(5.4,1)(5.6,1.5)\qline(5.8,1)(6,1.4)
\qline(6.8,1)(7,1.4)\qline(7.2,1)(7.4,1.4)\qline(7.6,1)(7.8,1.4)
\qline(8,1)(8.2,1.4)\qline(8.4,1)(8.6,1.4)\qline(8.8,1)(9,1.4)
\qline(0,3.2)(.2,3.2)
\qline(0.6,3.2)(0.8,3.2)\qline(0.8,3.0)(1.2,3.0)\qline(0.8,3.4)(1.2,3.4)
\qline(1.2,2.8)(1.6,2.8)\qline(1.2,3.2)(1.6,3.2)\qline(1.2,3.6)(1.6,3.6)
\qline(1.6,2.6)(2.0,2.6)\qline(1.6,3.0)(2.0,3.0)\qline(1.6,3.4)(2.0,3.4)
\qline(1.6,3.8)(2.0,3.8)\qline(2.0,2.4)(2.4,2.4)\qline(2.0,2.8)(2.4,2.8)
\qline(2.0,3.2)(2.4,3.2)\qline(2.0,3.6)(2.4,3.6)\qline(2.0,4.0)(2.4,4.0)
\qline(2.4,2.2)(2.8,2.2)\qline(2.4,2.6)(2.8,2.6)\qline(2.4,3.0)(2.8,3.0)
\qline(2.4,3.4)(2.8,3.4)\qline(2.4,3.8)(2.8,3.8)\qline(2.4,4.2)(2.8,4.2)
\qline(2.8,2.0)(3.0,2.0)\qline(2.8,2.4)(3.0,2.4)\qline(2.8,2.8)(3.0,2.8)
\qline(2.8,3.2)(3.0,3.2)\qline(2.8,3.6)(3.0,3.6)\qline(2.8,4.0)(3.0,4.0)
\qline(2.8,4.4)(3.0,4.4)\qline(3.8,2.0)(4.0,2.0)\qline(3.8,2.4)(4.0,2.4)
\qline(3.8,3.6)(4.0,3.6)\qline(3.8,4.0)(4.0,4.0)\qline(3.8,4.4)(4.0,4.4)
\qline(4.0,1.8)(4.4,1.8)\qline(4.0,2.2)(4.4,2.2)\qline(4.0,3.8)(4.4,3.8)
\qline(4.0,4.2)(4.4,4.2)\qline(4.0,4.6)(4.4,4.6)\qline(4.4,1.6)(5.2,1.6)
\qline(5.2,2.0)(4.4,2.0)\qline(5.2,2.4)(4.4,2.4)\qline(5.2,4.0)(4.4,4.0)
\qline(5.2,4.4)(4.4,4.4)\qline(5.2,1.8)(5.6,1.8)\qline(5.2,2.2)(5.6,2.2)
\qline(5.2,3.8)(5.6,3.8)\qline(5.2,4.2)(5.6,4.2)\qline(5.2,3.4)(5.6,3.4)
\qline(5.6,2.0)(5.8,2.0)\qline(5.6,2.4)(5.8,2.4)\qline(5.6,4.0)(5.8,4.0)
\qline(5.6,3.6)(5.8,3.6)\qline(7.0,2.0)(6.8,2.0)\qline(6.8,2.4)(7.0,2.4)
\qline(7.0,4.0)(6.8,4.0)\qline(6.8,3.6)(7.0,3.6)\qline(7.0,2.2)(7.4,2.2)
\qline(7.0,2.6)(7.4,2.6)\qline(7.0,3.0)(7.4,3.0)\qline(7.0,3.4)(7.4,3.4)
\qline(7.0,3.8)(7.4,3.8)\qline(7.4,2.4)(7.8,2.4)\qline(7.4,2.8)(7.8,2.8)
\qline(7.4,3.2)(7.8,3.2)\qline(7.4,3.6)(7.8,3.6)\qline(7.8,2.6)(8.2,2.6)
\qline(7.8,3.0)(8.2,3.0)\qline(7.8,3.4)(8.2,3.4)\qline(8.2,2.8)(8.6,2.8)
\qline(8.2,3.2)(8.6,3.2)\qline(8.6,3)(9.0,3)
\psframe[linecolor=darkred](-.2,.4)(9,.6)
\endpspicture
\label{a2abexp5-1}\end{eqnarray}

Now we look at the $D_{i,a}$ where $i>1$ as in
figure~\ref{a2abexp5-1}. Since we picked where the $U$ turn
appears already, one can find a candidate for a stem as thick and
shaded (green in color) line from the right hand side but we can
not finish because the purple string can not be join to the bottom
clasp without being zero(it will force to have a generator caps
off). So only nonzero admissible stems should be obtained from
$\tilde D_{1,a}$. We split the rectangle(only one in the middle)
vertically(horizontal splitting vanishes immediately) and it
creates another rectangle at right top side of previous place. We
have to split vertically except in the last step, for this
rectangle, as in the figure~\ref{a2abexp6}, both splits do not
vanish. The vertical split gives us the first shape
figure~\ref{a2ab2exppo} and the horizontal split gives the third
shape in figure ~\ref{a2ab2exppo}.

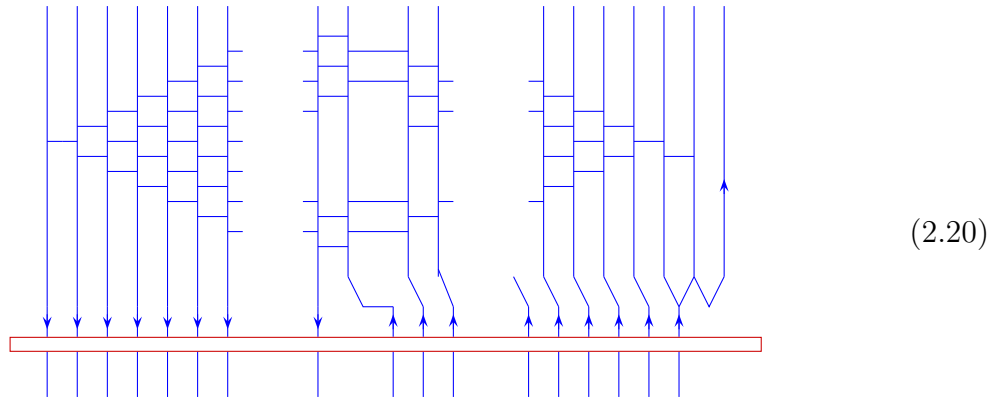
\begin{eqnarray}
\pspicture[.4](-.2,-.3)(10.1,5.1)
\qline(0.4,0.6)(0.4,1.0)\psline[arrowscale=1.5]{->}(0.4,.9)(0.4,.7)
\qline(0.8,0.6)(0.8,1.0)\psline[arrowscale=1.5]{->}(0.8,.9)(0.8,.7)
\qline(1.2,0.6)(1.2,1.0)\psline[arrowscale=1.5]{->}(1.2,.9)(1.2,.7)
\qline(1.6,0.6)(1.6,1.0)\psline[arrowscale=1.5]{->}(1.6,.9)(1.6,.7)
\qline(2,0.6)(2,1.0)\psline[arrowscale=1.5]{->}(2,.9)(2,.7)
\qline(2.4,0.6)(2.4,1.0)\psline[arrowscale=1.5]{->}(2.4,.9)(2.4,.7)
\qline(2.8,0.6)(2.8,1.0)\psline[arrowscale=1.5]{->}(2.8,.9)(2.8,.7)
\qline(4.0,0.6)(4.0,1.0)\psline[arrowscale=1.5]{->}(4.0,.9)(4.0,.7)
\qline(5.0,0.6)(5.0,1.0)\psline[arrowscale=1.5]{<-}(5.0,.9)(5.0,.7)
\qline(5.4,0.6)(5.4,1.0)\psline[arrowscale=1.5]{<-}(5.4,.9)(5.4,.7)
\qline(5.8,0.6)(5.8,1.0)\psline[arrowscale=1.5]{<-}(5.8,.9)(5.8,.7)
\qline(6.8,0.6)(6.8,1.0)\psline[arrowscale=1.5]{<-}(6.8,.9)(6.8,.7)
\qline(7.2,0.6)(7.2,1.0)\psline[arrowscale=1.5]{<-}(7.2,.9)(7.2,.7)
\qline(7.6,0.6)(7.6,1.0)\psline[arrowscale=1.5]{<-}(7.6,.9)(7.6,.7)
\qline(8.0,0.6)(8.0,1.0)\psline[arrowscale=1.5]{<-}(8,.9)(8,.7)
\qline(8.4,0.6)(8.4,1.0)\psline[arrowscale=1.5]{<-}(8.4,.9)(8.4,.7)
\qline(8.8,0.6)(8.8,1.0)\psline[arrowscale=1.5]{<-}(8.8,.9)(8.8,.7)
\qline(0.4,0.4)(0.4,-0.2) \qline(0.8,0.4)(0.8,-0.2)
\qline(1.2,0.4)(1.2,-0.2) \qline(1.6,0.4)(1.6,-0.2)
\qline(2,0.4)(2,-0.2) \qline(2.4,0.4)(2.4,-0.2)
\qline(2.8,0.4)(2.8,-0.2) \qline(4,0.4)(4,-0.2)
\qline(5.0,0.4)(5.0,-0.2) \qline(5.4,0.4)(5.4,-0.2)
\qline(5.8,0.4)(5.8,-0.2) \qline(6.8,0.4)(6.8,-0.2)
\qline(7.2,0.4)(7.2,-0.2) \qline(7.6,0.4)(7.6,-0.2)
\qline(8.0,0.4)(8.0,-0.2) \qline(8.4,0.4)(8.4,-0.2)
\qline(8.8,0.4)(8.8,-0.2)
\qline(0.4,1.0)(0.4,5.0)\qline(0.8,1.0)(0.8,5.0)\qline(1.2,1.0)(1.2,5.0)
\qline(1.6,1.0)(1.6,5.0)\qline(2.0,1.0)(2.0,5.0)\qline(2.4,1.0)(2.4,5.0)
\qline(2.8,1.0)(2.8,5.0)\qline(4.0,1.0)(4.0,5.0)\qline(4.4,1.4)(4.4,5.0)
\qline(5.2,1.4)(5.2,5.0)\qline(5.6,1.4)(5.6,5.0)\qline(7.0,1.4)(7.0,5.0)
\qline(7.4,1.4)(7.4,5.0)\qline(7.8,1.4)(7.8,5.0)\qline(8.2,1.4)(8.2,5.0)
\qline(9.0,1.4)(9.0,5.0)\qline(8.6,1.4)(8.6,5.0)
\psline(4.4,1.4)(4.6,1.0)(5.0,1.0)\qline(5.2,1.4)(5.4,1.0)
\qline(5.6,1.5)(5.8,1.0)\qline(6.6,1.4)(6.8,1.0)\qline(7.0,1.4)(7.2,1.0)
\qline(7.4,1.4)(7.6,1.0)\qline(7.8,1.4)(8,1)\qline(8.2,1.4)(8.4,1.0)
\qline(8.6,1.4)(8.8,1)\psline(9,1.4)(9.2,1)(9.4,1.4)\qline(0.4,3.2)(.6,3.2)
\qline(0.6,3.2)(0.8,3.2)\qline(0.8,3.0)(1.2,3.0)\qline(0.8,3.4)(1.2,3.4)
\qline(1.2,2.8)(1.6,2.8)\qline(1.2,3.2)(1.6,3.2)\qline(1.2,3.6)(1.6,3.6)
\qline(1.6,2.6)(2.0,2.6)\qline(1.6,3.0)(2.0,3.0)\qline(1.6,3.4)(2.0,3.4)
\qline(1.6,3.8)(2.0,3.8)\qline(2.0,2.4)(2.4,2.4)\qline(2.0,2.8)(2.4,2.8)
\qline(2.0,3.2)(2.4,3.2)\qline(2.0,3.6)(2.4,3.6)\qline(2.0,4.0)(2.4,4.0)
\qline(2.4,2.2)(2.8,2.2)\qline(2.4,2.6)(2.8,2.6)\qline(2.4,3.0)(2.8,3.0)
\qline(2.4,3.4)(2.8,3.4)\qline(2.4,3.8)(2.8,3.8)\qline(2.4,4.2)(2.8,4.2)
\qline(2.8,2.0)(3.0,2.0)\qline(2.8,2.4)(3.0,2.4)\qline(2.8,2.8)(3.0,2.8)
\qline(2.8,3.2)(3.0,3.2)\qline(2.8,3.6)(3.0,3.6)\qline(2.8,4.0)(3.0,4.0)
\qline(2.8,4.4)(3.0,4.4)\qline(3.8,2.0)(4.0,2.0)\qline(3.8,2.4)(4.0,2.4)
\qline(3.8,3.6)(4.0,3.6)\qline(3.8,4.0)(4.0,4.0)\qline(3.8,4.4)(4.0,4.4)
\qline(4.0,1.8)(4.4,1.8)\qline(4.0,2.2)(4.4,2.2)\qline(4.0,3.8)(4.4,3.8)
\qline(4.0,4.2)(4.4,4.2)\qline(4.0,4.6)(4.4,4.6)\qline(5.2,2.0)(4.4,2.0)
\qline(5.2,2.4)(4.4,2.4)\qline(5.2,4.0)(4.4,4.0)\qline(5.2,4.4)(4.4,4.4)
\qline(5.2,2.2)(5.6,2.2)\qline(5.2,3.8)(5.6,3.8)\qline(5.2,4.2)(5.6,4.2)
\qline(5.2,3.4)(5.6,3.4)\qline(5.6,2.4)(5.8,2.4)\qline(5.6,4.0)(5.8,4.0)
\qline(5.6,3.6)(5.8,3.6)\qline(6.8,2.4)(7.0,2.4)\qline(7.0,4.0)(6.8,4.0)
\qline(6.8,3.6)(7.0,3.6)\qline(7.0,2.6)(7.4,2.6)\qline(7.0,3.0)(7.4,3.0)
\qline(7.0,3.4)(7.4,3.4)\qline(7.0,3.8)(7.4,3.8)\qline(7.4,2.8)(7.8,2.8)
\qline(7.4,3.2)(7.8,3.2)\qline(7.4,3.6)(7.8,3.6)\qline(7.8,3.0)(8.2,3.0)
\qline(7.8,3.4)(8.2,3.4)\qline(8.2,3.2)(8.6,3.2)\qline(9.4,1.4)(9.4,5)
\qline(8.6,3)(9,3)\qline(8.8,1)(9,1.4)
\psline[arrowscale=1.5]{->}(9.4,2.5)(9.4,2.7)
\psframe[linecolor=darkred](-.1,.4)(9.9,.6)
\endpspicture
\label{a2abexp6}\end{eqnarray}

A similar argument works for the second one in
figure~\ref{a2ab2exppo}. The third figure in~\ref{a2ab2exppo} is a
little subtle. First one can see that none of the $\tilde D_{i,j}$
work if either $i>2$ or $j<a-1$. Thus, we only need to check
$\tilde D_{1,a-1}, \tilde D_{1,a}, \tilde D_{2,a-1}$ and $\tilde
D_{2,a}$ but we already know about $ \tilde D_{1,a}, \tilde
D_{2,a}$. The following figure~\ref{a2abexp7} shows the nonzero
admissible stem for $\tilde D_{1,a-1}$. As usual, we draw a stem
as a union of thick and purple lines.

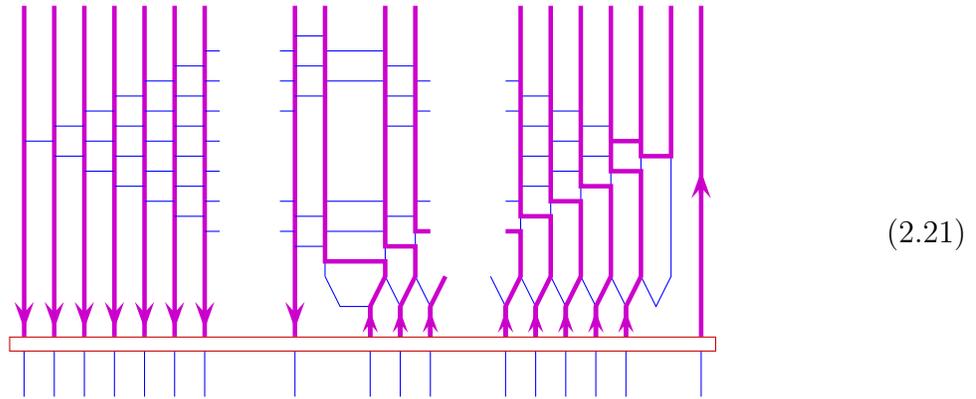
\begin{eqnarray}
\pspicture[.4](-.2,-.3)(10.1,5.1)
\psline[linecolor=pup,linewidth=1.7pt](0.4,0.6)(0.4,5.0)
\psline[arrowscale=1.5,linecolor=pup,linewidth=1.7pt]{->}(0.4,.9)(0.4,.7)
\psline[linecolor=pup,linewidth=1.7pt](0.8,0.6)(0.8,5.0)
\psline[arrowscale=1.5,linecolor=pup,linewidth=1.7pt]{->}(0.8,.9)(0.8,.7)
\psline[linecolor=pup,linewidth=1.7pt](1.2,0.6)(1.2,5.0)
\psline[arrowscale=1.5,linecolor=pup,linewidth=1.7pt]{->}(1.2,.9)(1.2,.7)
\psline[linecolor=pup,linewidth=1.7pt](1.6,0.6)(1.6,5.0)
\psline[arrowscale=1.5,linecolor=pup,linewidth=1.7pt]{->}(1.6,.9)(1.6,.7)
\psline[linecolor=pup,linewidth=1.7pt](2,0.6)(2,5.0)
\psline[arrowscale=1.5,linecolor=pup,linewidth=1.7pt]{->}(2,.9)(2,.7)
\psline[linecolor=pup,linewidth=1.7pt](2.4,0.6)(2.4,5.0)
\psline[arrowscale=1.5,linecolor=pup,linewidth=1.7pt]{->}(2.4,.9)(2.4,.7)
\psline[linecolor=pup,linewidth=1.7pt](2.8,0.6)(2.8,5.0)
\psline[arrowscale=1.5,linecolor=pup,linewidth=1.7pt]{->}(2.8,.9)(2.8,.7)
\psline[linecolor=pup,linewidth=1.7pt](4.0,0.6)(4.0,5.0)
\psline[arrowscale=1.5,linecolor=pup,linewidth=1.7pt]{->}(4.0,.9)(4.0,.7)
\psline[linecolor=pup,linewidth=1.7pt](5.0,0.6)(5.0,1.0)
\psline[arrowscale=1.2,linecolor=pup,linewidth=1.7pt]{<-}(5.0,.94)(5.0,.8)
\psline[linecolor=pup,linewidth=1.7pt](5.4,0.6)(5.4,1.0)
\psline[arrowscale=1.2,linecolor=pup,linewidth=1.7pt]{<-}(5.4,.94)(5.4,.8)
\psline[linecolor=pup,linewidth=1.7pt](5.8,0.6)(5.8,1.0)
\psline[arrowscale=1.2,linecolor=pup,linewidth=1.7pt]{<-}(5.8,.94)(5.8,.8)
\psline[linecolor=pup,linewidth=1.7pt](6.8,0.6)(6.8,1.0)
\psline[arrowscale=1.2,linecolor=pup,linewidth=1.7pt]{<-}(6.8,.94)(6.8,.8)
\psline[linecolor=pup,linewidth=1.7pt](7.2,0.6)(7.2,1.0)
\psline[arrowscale=1.2,linecolor=pup,linewidth=1.7pt]{<-}(7.2,.94)(7.2,.8)
\psline[linecolor=pup,linewidth=1.7pt](7.6,0.6)(7.6,1.0)
\psline[arrowscale=1.2,linecolor=pup,linewidth=1.7pt]{<-}(7.6,.94)(7.6,.8)
\psline[linecolor=pup,linewidth=1.7pt](8.0,0.6)(8.0,1.0)
\psline[arrowscale=1.2,linecolor=pup,linewidth=1.7pt]{<-}(8,.94)(8,.8)
\psline[linecolor=pup,linewidth=1.7pt](8.4,0.6)(8.4,1.0)
\psline[arrowscale=1.2,linecolor=pup,linewidth=1.7pt]{<-}(8.4,.94)(8.4,.8)
\psline[linecolor=pup,linewidth=1.7pt](9.4,0.6)(9.4,5.0)
\psline[arrowscale=1.5,linecolor=pup,linewidth=1.7pt]{<-}(9.4,2.8)(9.4,2.5)
\qline(0.4,0.4)(0.4,-0.2) \qline(0.8,0.4)(0.8,-0.2)
\qline(1.2,0.4)(1.2,-0.2) \qline(1.6,0.4)(1.6,-0.2)
\qline(2,0.4)(2,-0.2) \qline(2.4,0.4)(2.4,-0.2)
\qline(2.8,0.4)(2.8,-0.2) \qline(4,0.4)(4,-0.2)
\qline(5.0,0.4)(5.0,-0.2) \qline(5.4,0.4)(5.4,-0.2)
\qline(5.8,0.4)(5.8,-0.2) \qline(6.8,0.4)(6.8,-0.2)
\qline(7.2,0.4)(7.2,-0.2) \qline(7.6,0.4)(7.6,-0.2)
\qline(8.0,0.4)(8.0,-0.2) \qline(8.4,0.4)(8.4,-0.2)
\qline(9.4,0.4)(9.4,-0.2) \qline(4.4,1.4)(4.4,1.6)
\qline(5.2,1.6)(5.2,1.8)\qline(5.6,1.8)(5.6,2.0)\qline(7.0,2)(7.0,2.2)
\qline(7.4,2.4)(7.4,2.2)\qline(7.8,2.4)(7.8,2.6)\qline(8.2,2.6)(8.2,2.8)
\qline(8.6,2.8)(8.6,3.0)\psline(9,3)(9,1.4)(8.8,1)(8.6,1.4)
\psline(4.4,1.4)(4.6,1.0)(5.0,1.0)\qline(5.2,1.4)(5.4,1.0)\qline(5.6,1.4)(5.8,1.0)
\qline(6.6,1.4)(6.8,1)\qline(7,1.4)(7.2,1)\qline(7.4,1.4)(7.6,1)\qline(7.8,1.4)(8,1)
\qline(8.2,1.4)(8.4,1)
\qline(0.4,3.2)(0.8,3.2)\qline(0.8,3.0)(1.2,3.0)\qline(0.8,3.4)(1.2,3.4)
\qline(1.2,2.8)(1.6,2.8)\qline(1.2,3.2)(1.6,3.2)\qline(1.2,3.6)(1.6,3.6)
\qline(1.6,2.6)(2.0,2.6)\qline(1.6,3.0)(2.0,3.0)\qline(1.6,3.4)(2.0,3.4)
\qline(1.6,3.8)(2.0,3.8)\qline(2.0,2.4)(2.4,2.4)\qline(2.0,2.8)(2.4,2.8)
\qline(2.0,3.2)(2.4,3.2)\qline(2.0,3.6)(2.4,3.6)\qline(2.0,4.0)(2.4,4.0)
\qline(2.4,2.2)(2.8,2.2)\qline(2.4,2.6)(2.8,2.6)\qline(2.4,3.0)(2.8,3.0)
\qline(2.4,3.4)(2.8,3.4)\qline(2.4,3.8)(2.8,3.8)\qline(2.4,4.2)(2.8,4.2)
\qline(2.8,2.0)(3.0,2.0)\qline(2.8,2.4)(3.0,2.4)\qline(2.8,2.8)(3.0,2.8)
\qline(2.8,3.2)(3.0,3.2)\qline(2.8,3.6)(3.0,3.6)\qline(2.8,4.0)(3.0,4.0)
\qline(2.8,4.4)(3.0,4.4)\qline(3.8,2.0)(4.0,2.0)\qline(3.8,2.4)(4.0,2.4)
\qline(3.8,3.6)(4.0,3.6)\qline(3.8,4.0)(4.0,4.0)\qline(3.8,4.4)(4.0,4.4)
\qline(4.0,1.8)(4.4,1.8)\qline(4.0,2.2)(4.4,2.2)\qline(4.0,3.8)(4.4,3.8)
\qline(4.0,4.2)(4.4,4.2)\qline(4.0,4.6)(4.4,4.6)\qline(5.2,2.0)(4.4,2.0)
\qline(5.2,2.4)(4.4,2.4)\qline(5.2,4.0)(4.4,4.0)\qline(5.2,4.4)(4.4,4.4)
\qline(5.2,2.2)(5.6,2.2)\qline(5.2,3.8)(5.6,3.8)\qline(5.2,4.2)(5.6,4.2)
\qline(5.2,3.4)(5.6,3.4)\qline(5.6,2.4)(5.8,2.4)\qline(5.6,4.0)(5.8,4.0)
\qline(5.6,3.6)(5.8,3.6)\qline(6.8,2.4)(7.0,2.4)\qline(7.0,4.0)(6.8,4.0)
\qline(6.8,3.6)(7.0,3.6)\qline(7.0,2.6)(7.4,2.6)\qline(7.0,3.0)(7.4,3.0)
\qline(7.0,3.4)(7.4,3.4)\qline(7.0,3.8)(7.4,3.8)\qline(7.4,2.8)(7.8,2.8)
\qline(7.4,3.2)(7.8,3.2)\qline(7.4,3.6)(7.8,3.6)\qline(7.8,3.0)(8.2,3.0)
\qline(7.8,3.4)(8.2,3.4)
\psline[linecolor=pup,linewidth=1.7pt](4.4,5)(4.4,1.6)(5.2,1.6)(5.2,1.4)(5.0,1)(5.0,.6)
\psline[linecolor=pup,linewidth=1.7pt](5.2,5)(5.2,1.8)(5.6,1.8)(5.6,1.4)(5.4,1)(5.4,.6)
\psline[linecolor=pup,linewidth=1.7pt](5.6,5)(5.6,2.0)(5.8,2.0)
\psline[linecolor=pup,linewidth=1.7pt](6,1.4)(5.8,1)(5.8,.6)
\psline[linecolor=pup,linewidth=1.7pt](6.8,2)(7,2)(7,1.4)(6.8,1)(6.8,.6)
\psline[linecolor=pup,linewidth=1.7pt](7,5)(7,2.2)(7.4,2.2)(7.4,1.4)(7.2,1)(7.2,.6)
\psline[linecolor=pup,linewidth=1.7pt](7.4,5)(7.4,2.4)(7.8,2.4)(7.8,1.4)(7.6,1)(7.6,.6)
\psline[linecolor=pup,linewidth=1.7pt](7.8,5)(7.8,2.6)(8.2,2.6)(8.2,1.4)(8,1)(8,.6)
\psline[linecolor=pup,linewidth=1.7pt](8.2,5)(8.2,2.8)(8.6,2.8)(8.6,1.4)(8.4,1)(8.4,.6)
\psline[linecolor=pup,linewidth=1.7pt](8.6,5)(8.6,3)(9,3)(9,5)
\psline[linecolor=pup,linewidth=1.7pt](8.2,3.2)(8.6,3.2)
\psframe[linecolor=darkred](.2,.4)(9.6,.6)
\endpspicture
\label{a2abexp7}\end{eqnarray}

The following figure~\ref{a2abexp7-1} shows the nonzero admissible
stem for $\tilde D_{2,a-1}$.
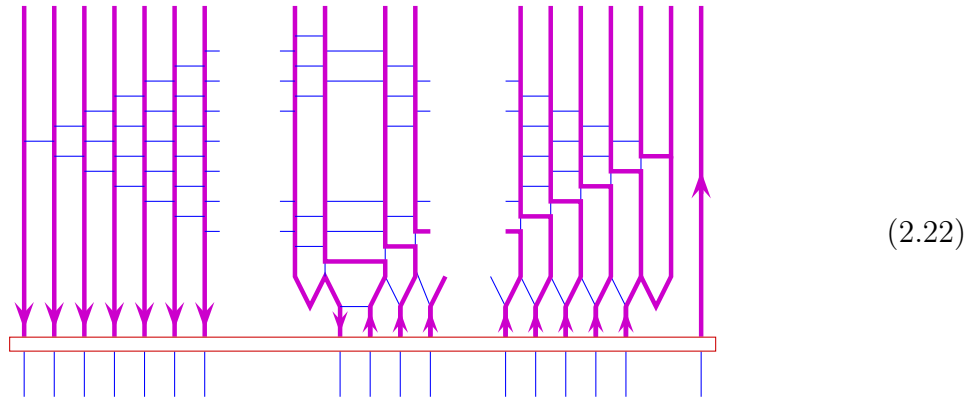
\begin{eqnarray}
\pspicture[.4](-.2,-.3)(10.1,5.1)
\psline[linecolor=pup,linewidth=1.7pt](0.4,0.6)(0.4,5.0)
\psline[arrowscale=1.5,linecolor=pup,linewidth=1.7pt]{->}(0.4,.9)(0.4,.7)
\psline[linecolor=pup,linewidth=1.7pt](0.8,0.6)(0.8,5.0)
\psline[arrowscale=1.5,linecolor=pup,linewidth=1.7pt]{->}(0.8,.9)(0.8,.7)
\psline[linecolor=pup,linewidth=1.7pt](1.2,0.6)(1.2,5.0)
\psline[arrowscale=1.5,linecolor=pup,linewidth=1.7pt]{->}(1.2,.9)(1.2,.7)
\psline[linecolor=pup,linewidth=1.7pt](1.6,0.6)(1.6,5.0)
\psline[arrowscale=1.5,linecolor=pup,linewidth=1.7pt]{->}(1.6,.9)(1.6,.7)
\psline[linecolor=pup,linewidth=1.7pt](2,0.6)(2,5.0)
\psline[arrowscale=1.5,linecolor=pup,linewidth=1.7pt]{->}(2,.9)(2,.7)
\psline[linecolor=pup,linewidth=1.7pt](2.4,0.6)(2.4,5.0)
\psline[arrowscale=1.5,linecolor=pup,linewidth=1.7pt]{->}(2.4,.9)(2.4,.7)
\psline[linecolor=pup,linewidth=1.7pt](2.8,0.6)(2.8,5.0)
\psline[arrowscale=1.5,linecolor=pup,linewidth=1.7pt]{->}(2.8,.9)(2.8,.7)
\psline[linecolor=pup,linewidth=1.7pt](4,5)(4,1.4)(4.2,1)(4.4,1.4)(4.6,1)(4.6,.6)
\psline[arrowscale=1.2,linecolor=pup,linewidth=1.7pt]{->}(4.6,.8)(4.6,.65)
\psline[linecolor=pup,linewidth=1.7pt](5.0,0.6)(5.0,1.0)
\psline[arrowscale=1.2,linecolor=pup,linewidth=1.7pt]{<-}(5.0,.94)(5.0,.8)
\psline[linecolor=pup,linewidth=1.7pt](5.4,0.6)(5.4,1.0)
\psline[arrowscale=1.2,linecolor=pup,linewidth=1.7pt]{<-}(5.4,.94)(5.4,.8)
\psline[linecolor=pup,linewidth=1.7pt](5.8,0.6)(5.8,1.0)
\psline[arrowscale=1.2,linecolor=pup,linewidth=1.7pt]{<-}(5.8,.94)(5.8,.8)
\psline[linecolor=pup,linewidth=1.7pt](6.8,0.6)(6.8,1.0)
\psline[arrowscale=1.2,linecolor=pup,linewidth=1.7pt]{<-}(6.8,.94)(6.8,.8)
\psline[linecolor=pup,linewidth=1.7pt](7.2,0.6)(7.2,1.0)
\psline[arrowscale=1.2,linecolor=pup,linewidth=1.7pt]{<-}(7.2,.94)(7.2,.8)
\psline[linecolor=pup,linewidth=1.7pt](7.6,0.6)(7.6,1.0)
\psline[arrowscale=1.2,linecolor=pup,linewidth=1.7pt]{<-}(7.6,.94)(7.6,.8)
\psline[linecolor=pup,linewidth=1.7pt](8.0,0.6)(8.0,1.0)
\psline[arrowscale=1.2,linecolor=pup,linewidth=1.7pt]{<-}(8,.94)(8,.8)
\psline[linecolor=pup,linewidth=1.7pt](8.4,0.6)(8.4,1.0)
\psline[arrowscale=1.2,linecolor=pup,linewidth=1.7pt]{<-}(8.4,.94)(8.4,.8)
\psline[linecolor=pup,linewidth=1.7pt](9.4,0.6)(9.4,5.0)
\psline[arrowscale=1.5,linecolor=pup,linewidth=1.7pt]{<-}(9.4,2.8)(9.4,2.5)
\qline(0.4,0.4)(0.4,-0.2) \qline(0.8,0.4)(0.8,-0.2)
\qline(1.2,0.4)(1.2,-0.2) \qline(1.6,0.4)(1.6,-0.2)
\qline(2,0.4)(2,-0.2) \qline(2.4,0.4)(2.4,-0.2)
\qline(2.8,0.4)(2.8,-0.2) \qline(4.6,0.4)(4.6,-0.2)
\qline(5.0,0.4)(5.0,-0.2) \qline(5.4,0.4)(5.4,-0.2)
\qline(5.8,0.4)(5.8,-0.2) \qline(6.8,0.4)(6.8,-0.2)
\qline(7.2,0.4)(7.2,-0.2) \qline(7.6,0.4)(7.6,-0.2)
\qline(8.0,0.4)(8.0,-0.2) \qline(8.4,0.4)(8.4,-0.2)
\qline(9.4,0.4)(9.4,-0.2) \qline(4.4,1.4)(4.4,1.6)
\qline(5.2,1.6)(5.2,1.8)\qline(5.6,1.8)(5.6,2.0)\qline(7.0,2)(7.0,2.2)
\qline(7.4,2.4)(7.4,2.2)\qline(7.8,2.4)(7.8,2.6)\qline(8.2,2.6)(8.2,2.8)
\qline(8.6,2.8)(8.6,3.0)
\qline(4.6,1.0)(5.0,1.0)\qline(5.2,1.4)(5.4,1.0)\qline(5.6,1.5)(5.8,1.0)
\qline(6.6,1.4)(6.8,1)\qline(7,1.4)(7.2,1)\qline(7.4,1.4)(7.6,1)\qline(7.8,1.4)(8,1)
\qline(8.2,1.4)(8.4,1)
\qline(0.4,3.2)(0.8,3.2)\qline(0.8,3.0)(1.2,3.0)\qline(0.8,3.4)(1.2,3.4)
\qline(1.2,2.8)(1.6,2.8)\qline(1.2,3.2)(1.6,3.2)\qline(1.2,3.6)(1.6,3.6)
\qline(1.6,2.6)(2.0,2.6)\qline(1.6,3.0)(2.0,3.0)\qline(1.6,3.4)(2.0,3.4)
\qline(1.6,3.8)(2.0,3.8)\qline(2.0,2.4)(2.4,2.4)\qline(2.0,2.8)(2.4,2.8)
\qline(2.0,3.2)(2.4,3.2)\qline(2.0,3.6)(2.4,3.6)\qline(2.0,4.0)(2.4,4.0)
\qline(2.4,2.2)(2.8,2.2)\qline(2.4,2.6)(2.8,2.6)\qline(2.4,3.0)(2.8,3.0)
\qline(2.4,3.4)(2.8,3.4)\qline(2.4,3.8)(2.8,3.8)\qline(2.4,4.2)(2.8,4.2)
\qline(2.8,2.0)(3.0,2.0)\qline(2.8,2.4)(3.0,2.4)\qline(2.8,2.8)(3.0,2.8)
\qline(2.8,3.2)(3.0,3.2)\qline(2.8,3.6)(3.0,3.6)\qline(2.8,4.0)(3.0,4.0)
\qline(2.8,4.4)(3.0,4.4)\qline(3.8,2.0)(4.0,2.0)\qline(3.8,2.4)(4.0,2.4)
\qline(3.8,3.6)(4.0,3.6)\qline(3.8,4.0)(4.0,4.0)\qline(3.8,4.4)(4.0,4.4)
\qline(4.0,1.8)(4.4,1.8)\qline(4.0,2.2)(4.4,2.2)\qline(4.0,3.8)(4.4,3.8)
\qline(4.0,4.2)(4.4,4.2)\qline(4.0,4.6)(4.4,4.6)\qline(5.2,2.0)(4.4,2.0)
\qline(5.2,2.4)(4.4,2.4)\qline(5.2,4.0)(4.4,4.0)\qline(5.2,4.4)(4.4,4.4)
\qline(5.2,2.2)(5.6,2.2)\qline(5.2,3.8)(5.6,3.8)\qline(5.2,4.2)(5.6,4.2)
\qline(5.2,3.4)(5.6,3.4)\qline(5.6,2.4)(5.8,2.4)\qline(5.6,4.0)(5.8,4.0)
\qline(5.6,3.6)(5.8,3.6)\qline(6.8,2.4)(7.0,2.4)\qline(7.0,4.0)(6.8,4.0)
\qline(6.8,3.6)(7.0,3.6)\qline(7.0,2.6)(7.4,2.6)\qline(7.0,3.0)(7.4,3.0)
\qline(7.0,3.4)(7.4,3.4)\qline(7.0,3.8)(7.4,3.8)\qline(7.4,2.8)(7.8,2.8)
\qline(7.4,3.2)(7.8,3.2)\qline(7.4,3.6)(7.8,3.6)\qline(7.8,3.0)(8.2,3.0)
\qline(7.8,3.4)(8.2,3.4)\qline(8.2,3.2)(8.6,3.2)
\psline[linecolor=pup,linewidth=1.7pt](4.4,5)(4.4,1.6)(5.2,1.6)(5.2,1.4)(5.0,1)(5.0,.6)
\psline[linecolor=pup,linewidth=1.7pt](5.2,5)(5.2,1.8)(5.6,1.8)(5.6,1.4)(5.4,1)(5.4,.6)
\psline[linecolor=pup,linewidth=1.7pt](5.6,5)(5.6,2.0)(5.8,2.0)
\psline[linecolor=pup,linewidth=1.7pt](6,1.4)(5.8,1)(5.8,.6)
\psline[linecolor=pup,linewidth=1.7pt](6.8,2)(7,2)(7,1.4)(6.8,1)(6.8,.6)
\psline[linecolor=pup,linewidth=1.7pt](7,5)(7,2.2)(7.4,2.2)(7.4,1.4)(7.2,1)(7.2,.6)
\psline[linecolor=pup,linewidth=1.7pt](7.4,5)(7.4,2.4)(7.8,2.4)(7.8,1.4)(7.6,1)(7.6,.6)
\psline[linecolor=pup,linewidth=1.7pt](7.8,5)(7.8,2.6)(8.2,2.6)(8.2,1.4)(8,1)(8,.6)
\psline[linecolor=pup,linewidth=1.7pt](8.2,5)(8.2,2.8)(8.6,2.8)(8.6,1.4)(8.4,1)(8.4,.6)
\psline[linecolor=pup,linewidth=1.7pt](8.6,5)(8.6,3)(9,3)(9,5)
\psline[linecolor=pup,linewidth=1.7pt](9,3)(9,1.4)(8.8,1)(8.6,1.4)
\psframe[linecolor=darkred](.2,.4)(9.6,.6)
\endpspicture
\label{a2abexp7-1}\end{eqnarray}

Note that the last figure has one loop which contributes $-[2]$.
This completes the proof of lemma.
\end{proof}

\begin{cor}
\begin{eqnarray}
\pspicture[.45](-2,-1.3)(2,1.3)
\pcline(-.5,.1)(-.5,.9)\middlearrow
\pcline(-.5,-.9)(-.5,-.1)\middlearrow \rput[b](-.5,1){$a$}
\rput[t](-.5,-1){$a$} \psarc(1,.1){.7}{0}{180}
\psarc(1,-.1){.7}{180}{360}
\pcline(1.7,.1)(1.7,-.1)\middlearrow\Aput{$1$}
\psframe[linecolor=darkred](-1,-.1)(1,.1)
\endpspicture
&=&\frac{[a+3]}{[a+1]} \pspicture[.45](-1.1,-1.3)(.2,1.3)
\pcline(-.45,.1)(-.45,.9)\middlearrow
\pcline(-.45,-.9)(-.45,-.1)\middlearrow \rput[b](-.45,1){$a$}
\rput[t](-.45,-1){$a$} \psframe[linecolor=darkred](-1,-.1)(.1,.1)
\endpspicture\nonumber\\
\pspicture[.45](-2,-1.4)(2,1.4)
\pcline(-.7,.1)(-.7,.9)\middlearrow
\pcline(-.7,-.9)(-.7,-.1)\middlearrow
\pcline(-.2,.9)(-.2,.1)\middlearrow
\pcline(-.2,-.1)(-.2,-.9)\middlearrow \rput[b](-.7,1){$a$}
\rput[b](-.2,1){$b$} \rput[t](-.7,-1){$a$} \rput[t](-.2,-1){$b$}
\psarc(1,.1){.7}{0}{180} \psarc(1,-.1){.7}{180}{360}
\pcline(1.7,-.1)(1.7,.1)\middlearrow\Aput{$1$}
\psframe[linecolor=darkred](-1,-.1)(1,.1)
\endpspicture
&=&\frac{[b+2][a+b+3]}{[b+1][a+b+2]}
\pspicture[.45](-1.1,-1.4)(.2,1.4)
\pcline(-.7,.1)(-.7,.9)\middlearrow
\pcline(-.7,-.9)(-.7,-.1)\middlearrow
\pcline(-.2,.9)(-.2,.1)\middlearrow
\pcline(-.2,-.1)(-.2,-.9)\middlearrow \rput[b](-.7,1){$a$}
\rput[b](-.2,1){$b$} \rput[t](-.7,-1){$a$} \rput[t](-.2,-1){$b$}
\psframe[linecolor=darkred](-1,-.1)(.1,.1)
\endpspicture\nonumber\\
\pspicture[.45](-1.2,-1)(2.2,1.4)
\pccurve[angleA=90,angleB=90,ncurv=1](-1,0)(1,0)\middlearrow\Aput{$a$}
\pccurve[angleA=270,angleB=270,ncurv=1](-1,0)(1,0)
\psframe[linecolor=darkred,fillstyle=solid](0,-.1)(2,.1)
\endpspicture
&=&\frac{[a+2][a+1]}{[2]} \nonumber\\
\pspicture[.45](-2,-1)(2,1) \psarc(-1,.1){.7}{0}{180}
\psarc(-1,-.1){.7}{180}{360}
\pcline(-1.7,.1)(-1.7,-.1)\middlearrow\Aput{$a$}
\psarc(1,.1){.7}{0}{180} \psarc(1,-.1){.7}{180}{360}
\pcline(1.7,-.1)(1.7,.1)\middlearrow\Aput{$b$}
\psframe[linecolor=darkred](-1,-.1)(1,.1)
\endpspicture
&=&\frac{[a+1][b+1][a+b+2]}{[2]} \label{a2triexp2}
\end{eqnarray}
\begin{proof}
After using a double clasps expansion one can get the first two
equalities with a simple calculation. The next two follow from the
previous two by induction.
\end{proof}
\label{a2trilem1}
\end{cor}

We will apply theorem~\ref{a2abexpthm1} to derive the coefficients
in equations~\ref{a2abexp1} and ~\ref{a2baexp1}. The expansion in
the proposition~\ref{a2abexp1} is known~\cite{Kuperberg:spiders},
which is only previously known expansion formula for a segregated
clasp of weight $(a,b)$ and it was used to find quantum $su(3)$
invariants in~\cite{OY:quantum}. Our proof using single clasp
expansion will be used for trihedron coefficients.

\begin{eqnarray}
\pspicture[.4](-.5,-1)(1.3,2)
\rput[t](0,-.4){$a$}\psline(0,-.3)(0,.4)\psline[arrowscale=1.5]{->}(0,0)(0,.2)
\rput[t](.7,-.4){$b$}\psline(.7,.4)(.7,-.3)\psline[arrowscale=1.5]{->}(.7,.2)(.7,0)
\psline(0,.6)(0,1.4)\psline[arrowscale=1.5]{->}(0,.9)(0,1.1)
\psline(.7,1.4)(.7,.6)\psline[arrowscale=1.5]{->}(.7,1.1)(.7,.9)
\rput[b](0,1.5){$a$} \rput[b](.7,1.5){$b$}
\psframe[linecolor=darkred](-.2,.4)(.9,.6)
\endpspicture
= \sum_{k=0}^{\mathrm{Min}(a,b)} a_k
\pspicture[.45](-2.9,-2)(2.7,2)
\psline(-1.3,-.9)(-1.3,.9)\psline[arrowscale=1.5]{->}(-1.3,-.1)(-1.3,.1)
\rput[br](-1.5,0){$a-k$}
\psline(1.3,.9)(1.3,-.9)\psline[arrowscale=1.5]{->}(1.3,.1)(1.3,-.1)
\rput[bl](1.5,0){$b-k$}
\pccurve[angleA=270,angleB=270,ncurv=1](.7,.9)(-.7,.9)\middlearrow
\pccurve[angleA=90,angleB=90,ncurv=1](-.7,-.9)(.7,-.9)\middlearrow
\rput(0,0){$k$}
\psline(-1,-1.5)(-1,-1.1)\psline[arrowscale=1.5]{->}(-1,-1.4)(-1,-1.2)
\psline(-1,1.1)(-1,1.5)\psline[arrowscale=1.5]{->}(-1,1.2)(-1,1.4)
\psline(1,-1.1)(1,-1.5)\psline[arrowscale=1.5]{->}(1,-1.2)(1,-1.4)
\psline(1,1.5)(1,1.1)\psline[arrowscale=1.5]{->}(1,1.4)(1,1.2)
\rput[b](-1,1.7){$a$}\rput[b](1,1.7){$b$}
\rput[t](-1,-1.7){$a$}\rput[t](1,-1.7){$b$}
\psframe[linecolor=darkred](-1.5,-1.1)(-.5,-.9)
\psframe[linecolor=darkred]( 1.5,-1.1)( .5,-.9)
\psframe[linecolor=darkred](-1.5, 1.1)(-.5, .9)
\psframe[linecolor=darkred]( 1.5, 1.1)( .5, .9)
\endpspicture
\label{a2abexp1}
\end{eqnarray}

\begin{prop}
The coefficients in equation~\ref{a2abexp1} is
$$a_k=(-1)^k\frac{[a]![b]![a+b-k+1]!}{[a-k]![b-k]![k]![a+b+1]!}.$$
\label{a2abexpprop1}
\end{prop}
\begin{proof}
Let me denote that a basis web in the right side of equation
~\ref{a2abexp1} by $D(k)$ which corresponding to the coefficient
$a_k$. We induct on $a+b$. It is clear for $a=0$ or $b=0$. If
$a\neq 0 \neq b$ then we use a segregated single clasp expansion
of weight $(a,b)$ in the middle. Even if we do not use entire
single clasp expansion of segregated clasp, once we attach $(a,0),
(0,b)$ clasps on the top, there is only two surviving web which
are one with one $U$ turn. One of resulting webs has some $H$'s as
in Figure~\ref{a2abexplem5} but if we push them down to $(a,b-1)$
clasp, it becomes a non-segregated clasp.

\begin{eqnarray}&
\pspicture[.47](-.5,-1)(1.3,2)
\rput[t](0,-.4){$a$}\psline(0,-.3)(0,.4)\psline[arrowscale=1.5]{->}(0,0)(0,.2)
\rput[t](.7,-.4){$b$}\psline(.7,.4)(.7,-.3)\psline[arrowscale=1.5]{->}(.7,.2)(.7,0)
\psline(0,.6)(0,1.4)\psline[arrowscale=1.5]{->}(0,.9)(0,1.1)
\psline(.7,1.4)(.7,.6)\psline[arrowscale=1.5]{->}(.7,1.1)(.7,.9)
\rput[b](0,1.5){$a$} \rput[b](.7,1.5){$b$}
\psframe[linecolor=darkred](-.2,.4)(.9,.6)
\endpspicture  =
\pspicture[.45](-1.9,-2)(2.2,2)
\psline(-1,-.3)(-1,.9)\psline[arrowscale=1.5]{->}(-1,.2)(-1,.4)
\psline(-1,-.9)(-1,-.5)\psline[arrowscale=1.5]{->}(-1,-.8)(-1,-.6)
\psline(1,.9)(1,-.3)\psline[arrowscale=1.5]{->}(1,.4)(1,.2)
\psline(1,-.5)(1,-.9)\psline[arrowscale=1.5]{->}(1,-.6)(1,-.8)
\psline(1.4,.9)(1.4,-.9)\psline[arrowscale=1.5]{->}(1.4,.1)(1.4,-.1)
\psline(-1,-1.5)(-1,-1.1)\psline[arrowscale=1.5]{->}(-1,-1.4)(-1,-1.2)
\psline(-1,1.1)(-1,1.5)\psline[arrowscale=1.5]{->}(-1,1.2)(-1,1.4)
\psline(1,-1.1)(1,-1.5)\psline[arrowscale=1.5]{->}(1,-1.2)(1,-1.4)
\psline(1,1.5)(1,1.1)\psline[arrowscale=1.5]{->}(1,1.4)(1,1.2)
\rput[b](-1,1.7){$a$}\rput[b](1,1.7){$b$}
\rput[br](-1.2,.2){$a$}\rput[br](.8,.2){$b-1$}
\rput[t](-1,-1.7){$a$}\rput[t](1,-1.7){$b$}
\psframe[linecolor=darkred](-1.5,-1.1)(-.5,-.9)
\psframe[linecolor=darkred]( 1.5,-1.1)( .5,-.9)
\psframe[linecolor=darkred](-1.5, 1.1)(-.5, .9)
\psframe[linecolor=darkred]( 1.5, 1.1)( .5, .9)
\psframe[linecolor=darkred](-1.3,-.5)(1.3, -.3)
\endpspicture
-\frac{[a]}{[a+b+1]} \pspicture[.45](-2.4,-2)(2.2,2)
\psline(-1,-.3)(-1,.9)\psline[arrowscale=1.5]{->}(-1,.2)(-1,.4)
\rput[br](-1.2,0){$a-1$}
\psline(-1,-.9)(-1,-.5)\psline[arrowscale=1.5]{->}(-1,-.8)(-1,-.6)
\psline(1,.9)(1,.2)\psline[arrowscale=1.5]{->}(1,.7)(1,.5)
\rput[bl](1.2,.4){$b-1$}
\psline(1,-.5)(1,-.9)\psline[arrowscale=1.5]{->}(1,-.6)(1,-.8)
\pccurve[angleA=90,angleB=90,ncurv=1](1.2,.2)(1.4,.2)
\pccurve[angleA=225,angleB=315,ncurv=1](.7,.9)(-.7,.9)\middlearrow
\psline(1.4,.2)(1.4,-.9)\psline[arrowscale=1.5]{->}(1.4,-.3)(1.4,-.5)
\psline(-1,-1.5)(-1,-1.1)\psline[arrowscale=1.5]{->}(-1,-1.4)(-1,-1.2)
\psline(-1,1.1)(-1,1.5)\psline[arrowscale=1.5]{->}(-1,1.2)(-1,1.4)
\psline(1,-1.1)(1,-1.5)\psline[arrowscale=1.5]{->}(1,-1.2)(1,-1.4)
\psline(1,1.5)(1,1.1)\psline[arrowscale=1.5]{->}(1,1.4)(1,1.2)
\psline(.3,-.3)(.3,0)\psline[arrowscale=1.5]{->}(.3,-.25)(.3,-.05)
\psline(.5,0)(.5,-.3)\psline[arrowscale=1.5]{->}(.5,-.05)(.5,-.25)
\psline(1.2,0)(1.2,-.3)\psline[arrowscale=1.5]{->}(1.2,-.05)(1.2,-.25)
\rput[br](1.225,-.25){$\ldots$}
\rput[b](-1,1.7){$a$}\rput[b](1,1.7){$b$}
\rput[t](-1,-1.7){$a$}\rput[t](1,-1.7){$b$}
\psframe[linecolor=darkred](-1.5,-1.1)(-.5,-.9)
\psframe[linecolor=darkred]( 1.5,-1.1)( .5,-.9)
\psframe[linecolor=darkred](-1.5, 1.1)(-.5, .9)
\psframe[linecolor=darkred]( 1.5, 1.1)( .5, .9)
\psframe[linecolor=darkred](-1.3,-.5)(1.3, -.3)
\psframe[linecolor=pup](.2,0)(1.3, .2)
\endpspicture\nonumber\\
&= \pspicture[.45](-1.9,-2)(2.2,2)
\psline(-1,-.3)(-1,.9)\psline[arrowscale=1.5]{->}(-1,.2)(-1,.4)
\psline(-1,-.9)(-1,-.5)\psline[arrowscale=1.5]{->}(-1,-.8)(-1,-.6)
\psline(1,.9)(1,-.3)\psline[arrowscale=1.5]{->}(1,.4)(1,.2)
\psline(1,-.5)(1,-.9)\psline[arrowscale=1.5]{->}(1,-.6)(1,-.8)
\psline(1.4,.9)(1.4,-.9)\psline[arrowscale=1.5]{->}(1.4,.1)(1.4,-.1)
\psline(-1,-1.5)(-1,-1.1)\psline[arrowscale=1.5]{->}(-1,-1.4)(-1,-1.2)
\psline(-1,1.1)(-1,1.5)\psline[arrowscale=1.5]{->}(-1,1.2)(-1,1.4)
\psline(1,-1.1)(1,-1.5)\psline[arrowscale=1.5]{->}(1,-1.2)(1,-1.4)
\psline(1,1.5)(1,1.1)\psline[arrowscale=1.5]{->}(1,1.4)(1,1.2)
\rput[b](-1,1.7){$a$}\rput[b](1,1.7){$b$}
\rput[br](-1.2,.2){$a$}\rput[br](.8,.2){$b-1$}
\rput[t](-1,-1.7){$a$}\rput[t](1,-1.7){$b$}
\psframe[linecolor=darkred](-1.5,-1.1)(-.5,-.9)
\psframe[linecolor=darkred]( 1.5,-1.1)( .5,-.9)
\psframe[linecolor=darkred](-1.5, 1.1)(-.5, .9)
\psframe[linecolor=darkred]( 1.5, 1.1)( .5, .9)
\psframe[linecolor=darkred](-1.3,-.5)(1.3, -.3)
\endpspicture -\frac{[a][a+1]}{[a+b+1][a+b]} \pspicture[.45](-2.4,-2)(2.2,2)
\psline(-1,-.3)(-1,.9) \psline[arrowscale=1.5]{->}(-1,.2)(-1,.4)
\rput[br](-1.2,0){$a-1$} \psline(-1,-.9)(-1,-.5)
\psline[arrowscale=1.5]{->}(-1,-.8)(-1,-.6) \psline(1,.9)(1,-.3)
\psline[arrowscale=1.5]{->}(1,.7)(1,.5) \rput[bl](1.2,0){$b-1$}
\psline(1,-.5)(1,-.9) \psline[arrowscale=1.5]{->}(1,-.6)(1,-.8)
\pccurve[angleA=30,angleB=150,ncurv=1](-.7,-.9)(.7,-.9)\middlearrow
\pccurve[angleA=210,angleB=330,ncurv=1](.7,.9)(-.7,.9)\middlearrow
\rput[b](-1,1.7){$a$} \rput[b](1,1.7){$b$} \rput[t](-1,-1.7){$a$}
\rput[t](1,-1.7){$b$} \psline(-1,-1.5)(-1,-1.1)
\psline[arrowscale=1.5]{->}(-1,-1.4)(-1,-1.2)
\psline(-1,1.1)(-1,1.5)
\psline[arrowscale=1.5]{->}(-1,1.2)(-1,1.4)
\psline(1,-1.1)(1,-1.5)
\psline[arrowscale=1.5]{->}(1,-1.2)(1,-1.4) \psline(1,1.5)(1,1.1)
\psline[arrowscale=1.5]{->}(1,1.4)(1,1.2)
\psframe[linecolor=darkred](-1.5,-1.1)(-.5,-.9)
\psframe[linecolor=darkred]( 1.5,-1.1)( .5,-.9)
\psframe[linecolor=darkred](-1.5, 1.1)(-.5, .9)
\psframe[linecolor=darkred]( 1.5, 1.1)( .5, .9)
\psframe[linecolor=darkred](-1.3,-.5)(1.3, -.3)
\endpspicture\label{a2abexplem5}
\end{eqnarray}

We can find the coefficient using the same argument using stems
and it is

$$-[2]a_{1,a-b}+a_{a,a-b+1}+\sum_{i=2}^{b} (a_{i,a-b+i-2}-
[2]a_{i,a-b+i-1}+a_{i,a-b+i})=-\frac{[a]}{[a+b+1]}$$

because $a_{i,a-b+i-2}-[2]a_{i,a-b+i-1}+a_{i,a-b+i}=0$ for all
$i=2, 3, \ldots, b$. Then we attach some $H$'s to make the middle
clasp as a non-segregated clasp of weight $(a,b-1)$. By using a
non-segregated single clasp expansion for which clasps are located
at northeast corner and by the induction hypothesis, we have

\begin{align*}
&=\sum_{k=0}^{b-1}(-1)^{k}\frac{[a]![b-1]![a+b-k]!}{[a-k]![b-1-k]![k]![a+b]!}D(k)\\
&-\frac{[a+1][a]}{[a+b+1][a+b]}\sum_{k=0}^{b-1}(-1)^k
\frac{[a-1]![b-1]![a+b-1-k]!}{[a-1-k]![b-1-k]![k]![a+b-a]!}D(k+1)\\
&=1\cdot D(0)+\sum_{k=1}^{b-1}((-1)^k\frac{[a]![b-1]![a+b-k]!}{[a-k]![b-1-k]![k]![a+b]!}\\
&+(-1)^k\frac{[a+1]![b-1]![a+b-k]!}{[a-k]![b-1-k]![k-1]![a+b]!}) D(k)\\
&-(-1)^{b-1}\frac{[a+1][a]}{[a+b+1][a+b]}\frac{[a-1]![b-1]![a]!}
{[a-b]![0]![b-1]![a+b-1]!} D(b)\\
&=D(0)+\sum_{k=1}^{b-1}(-1)^k\frac{[a]![b]![a+b+1-k]!}{[a-k]![b-k]![k]![a+b+1]!}
(\frac{[b-k][a+b+1]+[k][a+1]}{[b][a+b+1-k]})D(k)\\
&+(-1)^b\frac{[a]![b-1]![a+1]!}{[a-b]![0]![b-1]![a+b+1]!}D(b)\\
&=\sum_{k=0}^{b}(-1)^k\frac{[a]![b]![a+b+1-k]!}{[a-k]![b-k]![k]![a+b+1]!}D(k)
\end{align*}
\end{proof}

For equation~\ref{a2baexp1}, we assume $0\le a\le b$.

\begin{eqnarray}
\pspicture[.4](-.5,-1)(1.3,2) \rput[t](0,-.4){$a$}
\psline(0,-.3)(0,.4) \psline[arrowscale=1.5]{->}(0,0)(0,.2)
\rput[t](.7,-.4){$b$} \psline(.7,.4)(.7,-.3)
\psline[arrowscale=1.5]{->}(.7,.2)(.7,0) \psline(0,.6)(0,1.4)
\psline[arrowscale=1.5]{->}(0,.9)(0,1.1) \psline(.7,1.4)(.7,.6)
\psline[arrowscale=1.5]{->}(.7,1.1)(.7,.9) \rput[b](0,1.5){$a$}
\rput[b](.7,1.5){$b$} \psframe[linecolor=darkred](-.2,.4)(.9,.6)
\endpspicture
= \sum_{k=0}^{b} A(a,b,k) \pspicture[.4](-2.5,-2.2)(2.3,2.2)
\pcline(-1,.9)(-.15,.2)\middlearrow
\pcline(-1,-.9)(-.55,-.4)\middlearrow
\pcline(.15,-.2)(1,-.9)\middlearrow
\pcline(.55,.4)(1,.9)\middlearrow
\psline[linecolor=emgreen](-.7,-.2)(-.4,-.6)(.7,.2)(.4,.6)(-.7,-.2)
\pccurve[angleA=330,angleB=210,ncurv=1](-.7,.9)(.7,.9)\middlearrow
\pccurve[angleA=30,angleB=150,ncurv=1](-.7,-.9)(.7,-.9)\middlearrow
\rput(0,.95){$k$} \rput(0,-.95){$k$}
\psline(-1,-1.5)(-1,-1.1)\psline[arrowscale=1.5]{->}(-1,-1.4)(-1,-1.2)
\psline(-1,1.5)(-1,1.1)\psline[arrowscale=1.5]{->}(-1,1.4)(-1,1.2)
\psline(1,-1.1)(1,-1.5)\psline[arrowscale=1.5]{->}(1,-1.2)(1,-1.4)
\psline(1,1.1)(1,1.5)\psline[arrowscale=1.5]{->}(1,1.2)(1,1.4)
\rput[b](-1,1.7){$b$} \rput[b](1,1.7){$a$} \rput[t](-1,-1.7){$a$}
\rput[t](1,-1.7){$b$}
\psframe[linecolor=darkred](-1.5,-1.1)(-.5,-.9)
\psframe[linecolor=darkred]( 1.5,-1.1)( .5,-.9)
\psframe[linecolor=darkred](-1.5, 1.1)(-.5, .9)
\psframe[linecolor=darkred]( 1.5, 1.1)( .5, .9)
\endpspicture
\label{a2baexp1}
\end{eqnarray}

\begin{prop}
The coefficients $A(a,b,k)$ in equation~\ref{a2baexp1} satisfies
the following recurrence relation.

\begin{align*}
A(1,1,1) &=\frac{[2]}{[3]},\\
A(1,1,0) &= 1,\\
A(a,a,k) &=\frac{[2a-k+1]}{[2a+1]}(A(a-1,a,k)+A(a-1,a,k-1)),\\
A(a,a+i,k) &=\frac{[2a+1+i-k]}{[2a+1+i]}A(a,a+i-1,k).
\end{align*}
\end{prop}
\begin{proof}
Note that we assume that $A(a,b,-i)=A(a,b,a+i)=0$ for all $i>0$.
Using a non-segregated single clasp expansion at the clasp of
weight $(a,b)$, one standard and one with clasp in the northeast
corner, we have the result with two axioms of clasps. Remark that
these coefficients are not round.
\end{proof}

\section{Single Clasp Expansion for $\mathcal{U}_q({\mathfrak sp}(4))$}

It is known \cite{Kuperberg:spiders} that
$\mathcal{U}_q({\mathfrak sp}(4))$ webs are generated by a single
web

$$
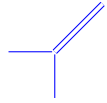
\pspicture[.4](-1,-1)(1,1) \pnode(.65;180){a1}
\rput(0,0){\rnode{a2}{$$}} \pnode(.65;270){a3} \pnode(.9;45){a4}
\ncline[nodesepA=1pt,nodesepB=1pt]{a1}{a2}
\ncline[nodesepA=1pt,nodesepB=1pt]{a2}{a3}
\ncline[doubleline=true]{a2}{a4}
\endpspicture
\label{b2gen}
$$

with the relations

\begin{eqnarray}
\pspicture[.4](-.6,-.5)(.6,.5) \pscircle(0,0){.4}
\endpspicture  &=&  - \frac{[6][2]}{[3]} \nonumber \\
\pspicture[.4](-.6,-.5)(.6,.5) \pscircle[doubleline=true](0,0){.4}
\endpspicture  &=&   \frac{[6][5]}{[3][2]} \nonumber \\
\pspicture[.4](-.6,-.5)(.6,.5) \psbezier(0,0)(.7,.7)(.7,-.7)(0,0)
\psline[doubleline=true](-.5,0)(0,0)
\endpspicture
&=&  0 \nonumber \\ \pspicture[.4](-.8,-.5)(.8,.5)
\pcarc[arcangle=45](-.3,0)(.3,0) \pcarc[arcangle=-45](-.3,0)(.3,0)
\psline[doubleline=true](-.7,0)(-.3,0)
\psline[doubleline=true](.3,0)(.7,0)
\endpspicture
&=&  -[2]^2\pspicture[.4](-.6,-.5)(.6,.5)
\psline[doubleline=true](-.5,0)(.5,0)
\endpspicture \nonumber \\
\pspicture[.4](-.9,-.9)(.9,.9) \pcarc[arcangle=-15](.4;90)(.4;210)
\pcarc[arcangle=-15](.4;210)(.4;330)
\pcarc[arcangle=-15](.4;330)(.4;90)
\psline[doubleline=true](.4;90)(.8;90)
\psline[doubleline=true](.4;210)(.8;210)
\psline[doubleline=true](.4;330)(.8;330)
\endpspicture
&= & 0 \nonumber\\
\btwohvert - \btwohhoriz  &=& \hh - \vv \label{b2rel}
\end{eqnarray}

Also it is known\cite{Kuperberg:spiders} that we can define
tetravalent vertex to achieve the same end as in
equation~\ref{b2gen2}. First we will find a single clasp expansion
of clasps of weight $(n,0)$ and $(0,n)$ and then use it to find
coefficients of double clasps expansion of clasps of weight
$(n,0)$ and $(0,n)$. Remark that the cut weight is defined
slightly different way. A cut path may cut diagonally through a
tetravalent vertex, and its weight is defined as $n\lambda_1 + (k
+ k')\lambda_2$, where $n$ is the number of type ``1'', single
strands, that it cuts, $k$ is the number of type ``2'', double
strands, that it cuts, and $k'$ is the number of tetravalent
vertices that it bisects. And there is a natural partial ordering
of the $B_2$ weight lattice given by

\begin{eqnarray}
a\lambda_1 + b\lambda_2 &\succ &(a-2)\lambda_1 + (b+1)\lambda_2 \nonumber\\
a\lambda_1 + b\lambda_2 &\succ &(a+2)\lambda_1 + (b-2)\lambda_2.
\nonumber
\end{eqnarray}

We will use the following shapes to find a single clasp expansion
because there is an ambiguity of preferred direction in the last
relation~\ref{b2rel}. We remark that the let side of the second
equality of~\ref{b2gen2} is not a crossing but a vertex where four
double edges meet.

\begin{eqnarray}
\pspicture[.4](-.8,-.8)(.8,.8) \pnode(.7;45){a1}
\pnode(.7;135){a2} \pnode(.7;225){a3} \pnode(.7;315){a4}
\ncline{a1}{a3} \ncline{a2}{a4} \endpspicture =
\pspicture[.4](-.8,-.8)(.8,.8) \pnode(.7;45){a1}
\pnode(.7;135){a2} \pnode(.7;225){a3} \pnode(.7;315){a4}
\pnode(.2;180){b1} \pnode(.2;0){b2} \ncline{a1}{b2}
\ncline{a2}{b1} \ncline{a3}{b1} \ncline{a4}{b2}
\ncline[doubleline=true]{b1}{b2}
\endpspicture
+ \pspicture[.4](-.8,-.8)(.8,.8) \pnode(.7;45){a1}
\pnode(.7;135){a2} \pnode(.7;225){a3} \pnode(.7;315){a4}
\nccurve[angleA=210,angleB=150,ncurv=1]{a1}{a4}
\nccurve[angleA=-30,angleB=30,ncurv=1]{a2}{a3}
\endpspicture \hskip .1cm  ,  \hskip .1cm
\pspicture[.4](-.8,-.8)(.8,.8) \pnode(.7;45){a1}
\pnode(.7;135){a2} \pnode(.7;225){a3} \pnode(.7;315){a4}
\ncline[doubleline=true]{a1}{a3} \ncline[doubleline=true]{a2}{a4}
\endpspicture
 =  \pspicture[.4](-.8,-.8)(.8,.8) \pnode(.7;45){a1}
\pnode(.7;135){a2} \pnode(.7;225){a3} \pnode(.7;315){a4}
\pnode(.25;45){b1} \pnode(.25;135){b2} \pnode(.25;225){b3}
\pnode(.25;315){b4} \ncline{b1}{b2} \ncline{b4}{b1}
\ncline{b3}{b2} \ncline{b4}{b3} \ncline[doubleline=true]{a1}{b1}
\ncline[doubleline=true]{a2}{b2} \ncline[doubleline=true]{a3}{b3}
\ncline[doubleline=true]{a4}{b4}
\endpspicture
\label{b2gen2}
\end{eqnarray}

By combining with the weight diagram of $V_{\lambda_1}^{\otimes
n}$ and minimal cut paths, we can find single clasp expansion of
$\mathcal{U}_q({\mathfrak sp}(4))$ of type $(n,0)$ as in
equation~\ref{b2n0exp} which has $\frac{n^2-1}{2}$ unknowns.

\begin{eqnarray}
\pspicture[.45](-.1,-.3)(2.3,2.3)
\rput[t](1,-.1){$n$}\qline(1,0)(1,.4)
\psframe[linecolor=darkred](0,.4)(2,.6)
\qline(1,.6)(1,2)\rput[b](1,2.1){$n$}
\endpspicture
= \sum_{i=0}^{n-1}\sum_{j=i+1}^{n} a_{ij}\hskip .2cm
\pspicture[.45](-.1,-.3)(3.55,2.5) \rput[t](1,-.1){$n-1$}
\qline(1,0)(1,.4) \psframe[linecolor=darkred](0,.4)(3.1,.6)
\qline(.2,.6)(.2,2) \qline(.6,.6)(.6,2) \qline(.8,.6)(.8,2)
\rput[l](1,2.3){$j$} \rput[l](2.3,2.3){$i$} \qline(3.1,2)(3.1,1.5)
\qline(3.1,1.5)(2.7,1) \qline(2.7,1)(2.7,.6)
\qline(2.5,2)(2.5,1.5) \qline(2.5,1.5)(2.1,1)
\qline(2.1,1)(2.1,.6) \qline(2.7,2)(2.7,1.5)
\qline(2.7,1.5)(2.3,1) \qline(2.3,1)(2.3,.6)
\psarc(3.1,.6){.2}{0}{180} \qline(3.3,0)(3.3,.6)
\qline(1.2,1.5)(1.2,2) \qline(1.2,1.5)(1,1)
\qline(1,1)(1,.6)\qline(1.4,1.5)(1.4,2) \qline(1.4,1.5)(1.2,1)
\qline(1.2,1)(1.2,.6) \qline(2,1.5)(2,2) \qline(2,1.5)(1.8,1)
\qline(1.8,1)(1.8,.6) \qline(1.8,1.5)(1.8,2)
\qline(1.8,1.5)(1.6,1) \qline(1.6,1)(1.6,.6) \qline(.8,1.7)(.8,2)
\rput(1,2){\rnode{a1}{$$}} \rput(2.3,2){\rnode{a6}{$$}}
\nccurve[angleA=-45,angleB=225,ncurv=1]{a1}{a6}
\endpspicture
\label{b2n0exp}
\end{eqnarray}

Since adding a $U$ turn and a $Y$ make the left side vanish, we
will have $n(n-1)$ equations.

\begin{thm}
The coefficients in Figure~\ref{b2n0exp} are
$$a_{i,j}=[2]^{i-j+1}\frac{[n+1][n-j+1][2n-2i+2]}{[n][2n+2][n-i+1]}.$$

\label{b2n0expthm}
\end{thm}

\begin{proof} To proceed the proof of the theorem, we remark that the following
useful relations in equation~\ref{b2n0help} can be easily obtained
from the relation in equation~\ref{b2rel}.

\begin{eqnarray}
\pspicture[.45](-.4,-.7)(1.2,.7) \psline(-.3,.6)(-.3,.3)
\psline(-.3,-.3)(-.3,-.6) \psarc(.3,.3){.6}{180}{270}
\psarc(.3,-.3){.6}{90}{180} \psarc(.3,0){.3}{270}{90}
\endpspicture
= \frac{[6][2]}{[3]} \hskip .3cm \pspicture[.45](-.1,-.7)(.1,.7)
\psline(-.1,.6)(-.1,-.6)
\endpspicture\hskip .5cm & , & \hskip .5cm
\pspicture[.45](-.4,-.7)(1.2,.7) \psline(-.3,.6)(-.3,.3)
\psline(-.3,-.3)(-.3,-.6) \psarc(.3,.3){.6}{180}{270}
\psarc(.3,-.3){.6}{90}{180} \psarc(.3,0){.3}{270}{90}
\psline[doubleline=true](.6,0)(1.1,0)
\endpspicture
= -[2]^2 \hskip .2cm \pspicture[.45](-.1,-.7)(1.5,.7)
\psline(-.1,.7)(.6,0) \psline(.6,0)(-.1,-.7)
\psline[doubleline=true](.6,0)(1.4,0)
\endpspicture\nonumber\\
\pspicture[.45](-1.2,-.7)(1.2,.7) \pnode(1.2;30){a1}
\pnode(1.2;150){a2} \pnode(1.2;210){a3} \pnode(1.2;330){a4}
\pnode(.6;90){b1} \pnode(.3;90){b2} \pnode(.3;270){b3}
\nccurve[angleA=315,angleB=180,ncurv=1]{a2}{b3}
\nccurve[angleA=225,angleB=0,ncurv=1]{a1}{b3}
\nccurve[angleA=135,angleB=0,ncurv=1]{a4}{b2}
\nccurve[angleA=45,angleB=180,ncurv=1]{a3}{b2}
\endpspicture
& = & -[2]^2 \pspicture[.45](-.8,-.7)(.8,.7) \pnode(.84;45){a1}
\pnode(.84;135){a2} \pnode(.84;225){a3} \pnode(.84;315){a4}
\ncline{a2}{a4} \ncline{a3}{a1} \endpspicture -[2][4]
\pspicture[.4](-.7,-.7)(.7,.7) \pnode(.84;45){a1}
\pnode(.84;135){a2} \pnode(.84;225){a3} \pnode(.84;315){a4}
\nccurve[angleA=135,angleB=225,ncurv=1]{a4}{a1}
\nccurve[angleA=45,angleB=-45,ncurv=1]{a3}{a2}
\endpspicture\nonumber\\
\pspicture[.45](-1.2,-.7)(1.2,.7) \pnode(1.2;30){a1}
\pnode(1.2;150){a2} \pnode(1.2;210){a3} \pnode(1.2;330){a4}
\pnode(.6;90){b1} \pnode(.3;90){b2} \pnode(.3;270){b3}
\nccurve[angleA=315,angleB=180,ncurv=1]{a2}{b3}
\nccurve[angleA=225,angleB=0,ncurv=1]{a1}{b3}
\nccurve[angleA=135,angleB=0,ncurv=1]{a4}{b2}
\nccurve[angleA=45,angleB=180,ncurv=1]{a3}{b2}
\ncline[doubleline=true]{b1}{b2}
\endpspicture
& = & [2]^2 \pspicture[.45](-.8,-.7)(.8,.7) \pnode(.84;45){a1}
\pnode(.84;135){a2} \pnode(.84;225){a3} \pnode(.84;315){a4}
\pnode(.6;90){b1} \pnode(.3;90){b2} \pnode(.3;270){b3}
\nccurve[angleA=-90,angleB=150,ncurv=1]{a2}{b2}
\nccurve[angleA=60,angleB=-40,ncurv=1]{a3}{b2}
\nccurve[angleA=135,angleB=225,ncurv=1]{a4}{a1}
\ncline[doubleline=true]{b1}{b2}
\endpspicture  + [2]^2 \pspicture[.45](-.8,-.7)(.8,.7) \pnode(.84;45){a1}
\pnode(.84;135){a2} \pnode(.84;225){a3} \pnode(.84;315){a4}
\pnode(.6;90){b1} \pnode(.3;90){b2} \pnode(.3;270){b3}
\nccurve[angleA=315,angleB=45,ncurv=1]{a2}{a3}
\nccurve[angleA=270,angleB=40,ncurv=1]{a1}{b2}
\nccurve[angleA=120,angleB=220,ncurv=1]{a4}{b2}
\ncline[doubleline=true]{b1}{b2}
\endpspicture\label{b2n0help}
\end{eqnarray}

Using these relations, we get the following $n-1$ equations by
adding $U$ turns from the left to right. By capping off the
generator from left to right, we have $(n-1)^2$ equations. There
are two special equations and four different shapes of equation as
follows.

$$a_{n-2,n-1}+\frac{[2][6]}{[3]} a_{n-2,n} - \frac{[2][6]}{[3]} a_{n-1,n}=0$$

$$-\frac{[2][6]}{[3]} a_{12} + \frac{[2][6]}{[3]} a_{13}
+ a_{23}+ 1+ \frac{[2][6]}{[3]} b_2 -[2][4] b_3=0$$

Type I : For $i=1, 2, \ldots, n-3$,

$$a_{i,i+1}+\frac{[2][6]}{[3]} a_{i,i+2} -[2][4]a_{i,i+3}  -
\frac{[2][6]}{[3]} a_{i+1,i+2} +\frac{[2][6]}{[3]}
a_{i+1,i+3}+a_{i+2,i+3}=0.$$

Type II : For $i=0, 1, \ldots, n-2$,

$$a_{i,n-1}-[2]^2a_{i,n}=0.$$

Type III : For $i= 0, 1, 2, \ldots, n-3$, $k=2, 3, \ldots, n-i-1$,

$$a_{i,n-k}-[2]^2a_{i,n-k+1}+[2]^2a_{i,n-k+2}=0.$$

Type IV : For $i=3, 4, \ldots, n$, $k=n-i+3, n-i+4, \ldots, n$,

$$[2]^2 a_{n-k,i}-[2]^2a_{n-k+1,i}+a_{n-k+2,i}=0.$$

Then we check the answer in the proposition satisfies the
equations. Since these webs in equation~\ref{b2n0exp} form a
basis, the coefficients are unique. Therefore, it completes the
proof.
\end{proof}

\begin{eqnarray}
\pspicture[.45](-.2,-.3)(1.2,1.3) \rput[t](.5,-.1){$n$}
\qline(.5,0)(.5,.4) \qline(.5,.6)(.5,1) \rput[b](.5,1.1){$n$}
\psframe[linecolor=darkred](0,.4)(1,.6)
\endpspicture
= \pspicture[.45](-.2,-.3)(1.4,1.3) \rput[t](.5,-0.1){$n-1$}
\qline(.5,0)(.5,.4) \qline(.5,.6)(.5,1) \rput[b](.5,1.1){$n-1$}
\qline(1.3,0)(1.3,1) \psframe[linecolor=darkred](0,.4)(1,.6)
\endpspicture
+ a_{12} \pspicture[.45](-.2,-.3)(1.45,2.3)
\rput[t](.5,-.1){$n-1$} \qline(.5,0)(.5,.4)
\psframe[linecolor=darkred](0,.4)(1,.6) \qline(.25,.6)(.25,1.4)
\rput[l](.35,1){$n-2$} \qline(.5,1.6)(.5,2)
\rput[b](.5,2.1){$n-1$} \psarc(1,.6){.2}{0}{180}
\qline(1.2,0)(1.2,.6) \psarc(1,1.4){.2}{180}{0}
\qline(1.2,1.4)(1.2,2) \psframe[linecolor=darkred](0,1.4)(1,1.6)
\endpspicture
+a_{02} \pspicture[.45](-.2,-.3)(1.45,2.3) \rput[t](.5,-.1){$n-1$}
\qline(.5,0)(.5,.4) \qline(.25,.6)(.25,1.4) \qline(.5,1.6)(.5,2)
\qline(1.2,1.4)(1.2,2) \rput[b](.5,2.1){$n-1$}
\qline(.8,.6)(1.2,1.4) \qline(1.2,0)(1.2,.6)
\qline(.8,1.4)(1.2,.6) \psframe[linecolor=darkred](0,.4)(1,.6)
\psframe[linecolor=darkred](0,1.4)(1,1.6)
\endpspicture
\label{b2n02exp}
\end{eqnarray}

\begin{cor}
The double clasp expansion of $B_2$ of type $(n,0)$ can be
obtained as in the equation~\ref{b2n02exp} where $a_{12}, a_{02}$
are from the Theorem~\ref{b2n0expthm}. \label{b2n0expcor}
\end{cor}

\begin{eqnarray}
\pspicture[.45](-.1,-.3)(2.3,2.3) \rput[t](1,-.1){$n$}
\psline[doubleline=true](1,0)(1,.4)
\psline[doubleline=true](1,.6)(1,2) \rput[b](1,2.1){$n$}
\psframe[linecolor=darkred](0,.4)(2,.6)
\endpspicture
= \sum_{i=0}^{n-1}\sum_{j=i+1}^{n} a_{ij}\hskip .2cm
\pspicture[.45](-.3,-.3)(3.55,2.5) \rput[t](1,-.1){$n-1$}
\psline[doubleline=true](1,0)(1,.4)
\psline[doubleline=true](.2,.6)(.2,2)
\psline[doubleline=true](.6,.6)(.6,2)
\psline[doubleline=true](.8,.6)(.8,2) \rput[l](1,2.3){$j$}
\rput[l](2.3,2.3){$i$} \psline[doubleline=true](3.1,2)(3.1,1.5)
\psline[doubleline=true](3.1,1.5)(2.7,1)
\psline[doubleline=true](2.7,1)(2.7,.6)
\psline[doubleline=true](2.5,2)(2.5,1.5)
\psline[doubleline=true](2.5,1.5)(2.1,1)
\psline[doubleline=true](2.1,1)(2.1,.6)
\psline[doubleline=true](2.7,2)(2.7,1.5)
\psline[doubleline=true](2.7,1.5)(2.3,1)
\psline[doubleline=true](2.3,1)(2.3,.6)
\psline[doubleline=true](3.3,0)(3.3,.6)
\psline[doubleline=true](1.2,1.5)(1.2,2)
\psline[doubleline=true](1.2,1.5)(1,1)
\psline[doubleline=true](1,1)(1,.6)
\psline[doubleline=true](1.4,1.5)(1.4,2)
\psline[doubleline=true](1.4,1.5)(1.2,1)
\psline[doubleline=true](1.2,1)(1.2,.6)
\psline[doubleline=true](2,1.5)(2,2)
\psline[doubleline=true](2,1.5)(1.8,1)
\psline[doubleline=true](1.8,1)(1.8,.6)
\psline[doubleline=true](1.8,1.5)(1.8,2)
\psline[doubleline=true](1.8,1.5)(1.6,1)
\psline[doubleline=true](1.6,1)(1.6,.6)
\psline[doubleline=true](.8,1.7)(.8,2) \rput(1,2){\rnode{a1}{$$}}
\rput(2.3,2){\rnode{a2}{$$}}
\nccurve[doubleline=true,angleA=-60,angleB=240,ncurv=1]{a1}{a2}
\psarc[doubleline=true](3.1,.6){.2}{0}{180}
\psframe[linecolor=darkred](0,.4)(3.1,.6)
\endpspicture
\label{b20nexp}
\end{eqnarray}

Then we look for $(0,n)$ case. The main idea for $(n,0)$ works
exactly same except we replace the base as in the
equation~\ref{b20nexp}. By capping off $U$ turns and a lower
weight cap, we get the following coefficients and we can solve
them successively as in Theorem~\ref{b20nexpthm}. Also the
equation~\ref{b20nhelp} is useful to find the following equations.

\begin{eqnarray}
\pspicture[.45](-.4,-.7)(.4,.7) \psline(-.3,.6)(0,.3)
\psline(0,-.3)(-.3,-.6)\psline(0,.3)(0,-.3)
\psarc[doubleline=true](0,0){.3}{270}{90}
\endpspicture
 = [5] \hskip .2cm \pspicture[.45](-.1,-.7)(.1,.7)
\psline(0,.6)(0,-.6)
\endpspicture &,& \hskip .2cm \pspicture[.45](-.4,-.7)(1.2,.7)
\psline[doubleline=true](-.3,.6)(-.3,.3)
\psline[doubleline=true](-.3,-.3)(-.3,-.6)
\psarc[doubleline=true](.3,.3){.6}{180}{270}
\psarc[doubleline=true](.3,-.3){.6}{90}{180}
\psarc[doubleline=true](.3,0){.3}{270}{90}
\endpspicture
= -[2]^2[5] \hskip .6cm \pspicture[.45](-.1,-.7)(.1,.7)
\psline[doubleline=true](-.1,.6)(-.1,-.6)
\endpspicture \nonumber\\
\pspicture[.45](-.4,-.7)(1,.7) \psline(-.3,.6)(0,.3)
\psline(0,-.3)(-.3,-.6) \psline(0,.3)(0,-.3)
\rput(0,.3){\rnode{a1}{$$}} \rput(.9,-.6){\rnode{a2}{$$}}
\rput(0,-.3){\rnode{a3}{$$}} \rput(.9,.6){\rnode{a4}{$$}}
\nccurve[doubleline=true,angleA=30,angleB=135,ncurv=1]{a1}{a2}
\nccurve[doubleline=true,angleA=-30,angleB=225,ncurv=1]{a3}{a4}
\endpspicture
& = &-[2][4] \pspicture[.45](-.4,-.7)(.4,.7) \psline(-.3,.6)(0,.3)
\psline(0,-.3)(-.3,-.6) \psline(0,.3)(0,-.3)
\psline[doubleline=true](0,-.3)(.3,-.6)
\psline[doubleline=true](0,.3)(.3,.6)
\endpspicture
-[2]^2[3] \pspicture[.45](-1,-.7)(1,.7)
\rput(.7,-.6){\rnode{a1}{$$}} \rput(.7,.6){\rnode{a2}{$$}}
\rput(-.7,-.6){\rnode{a3}{$$}} \rput(-.7,.6){\rnode{a4}{$$}}
\nccurve[doubleline=true,angleA=135,angleB=225,ncurv=1]{a1}{a2}
\nccurve[angleA=45,angleB=-45,ncurv=1]{a3}{a4}
\endpspicture\nonumber\\
\pspicture[.45](-1.2,-.7)(1.2,.7) \pnode(1.2;30){a1}
\pnode(1.2;150){a2} \pnode(1.2;210){a3} \pnode(1.2;330){a4}
\pnode(.6;90){b1} \pnode(.3;90){b2} \pnode(.3;270){b3}
\nccurve[doubleline=true,angleA=315,angleB=180,ncurv=1]{a2}{b3}
\nccurve[doubleline=true,angleA=225,angleB=0,ncurv=1]{a1}{b3}
\nccurve[doubleline=true,angleA=135,angleB=0,ncurv=1]{a4}{b2}
\nccurve[doubleline=true,angleA=45,angleB=180,ncurv=1]{a3}{b2}
\endpspicture
& = & -[2][4] \pspicture[.45](-.8,-.7)(.8,.7) \pnode(.84;45){a1}
\pnode(.84;135){a2} \pnode(.84;225){a3} \pnode(.84;315){a4}
\ncline[doubleline=true]{a2}{a4} \ncline[doubleline=true]{a3}{a1}
\endpspicture +[2]^4[3] \pspicture[.4](-.7,-.7)(.7,.7)
\pnode(.84;45){a1} \pnode(.84;135){a2} \pnode(.84;225){a3}
\pnode(.84;315){a4}
\nccurve[doubleline=true,angleA=135,angleB=225,ncurv=1]{a4}{a1}
\nccurve[doubleline=true,angleA=45,angleB=-45,ncurv=1]{a3}{a2}
\endpspicture
\label{b20nhelp}
\end{eqnarray}

$$a_{n-2,n-1}-[5][2]^2a_{n-2,n}+\frac{[6][5]}{[3][2]}a_{n-1,n}=0$$
$$-[3][2]^2a_{n-2,n}+[5]a_{n-1,n}=0$$

Type I : For $i=0, 1, \ldots, n-3$,
$$a_{i,i+1}-[5][2]^2a_{i,i+2}+[3][2]^4a_{i,i+3}
+\frac{[6][5]}{[3][2]}a_{i+1,i+2}-[5][2]^2a_{i+1,i+3}+a_{i+2,i+3}=0$$

Type II : For $i=0, 1, \ldots, n-2$,
$$a_{i,n-1}-[4][2]a_{i,n}=0$$

Type III : For $i=0,1,\ldots ,n-3$ and $j=i+1,i+2,\ldots, n-2$,
$$a_{i,j}-[4][2]a_{i,j+1}+[2]^4a_{i,j+2}=0$$

Type IV : For $i=0,1,\ldots, n-3$ and $j=i+3,i+4,\ldots, n$,
$$[2]^4a_{i,j}-[4][2]a_{i+1,j}+a_{i+2,j}=0$$

Type V : For $i=1,2,\ldots, n-2$
$$-[3][2]^2a_{i-1,i+1}+[2]^4a_{i-1,i+2}+[5]a_{i,i+1}-[3][2]^2a_{i,i+2}=0.$$

\begin{thm}
For  $n\ge 2$,
$$a_{i,j}=[2]^{2(1+i-j)}\frac{[2n+1-2i][2n-2j+2]}{[2n][2n+1]}.$$
\label{b20nexpthm}
\end{thm}

\begin{eqnarray}
\pspicture[.45](-.2,-.3)(1.2,1.3)
\rput[t](.5,-.1){$n$}\psline[doubleline=true](.5,0)(.5,.4)
\psline[doubleline=true](.5,.6)(.5,1)
\psframe[linecolor=darkred](0,.4)(1,.6) \rput[b](.5,1.1){$n$}
\endpspicture
= \pspicture[.45](-.2,-.3)(1.4,1.3)
\rput[t](.5,-0.1){$n-1$}\psline[doubleline=true](.5,0)(.5,.4)
\psline[doubleline=true](.5,.6)(.5,1)\rput[b](.5,1.1){$n-1$}
\psline[doubleline=true](1.3,0)(1.3,1)
\psframe[linecolor=darkred](0,.4)(1,.6)
\endpspicture
+ a_{12} \pspicture[.45](-.2,-.3)(1.45,2.3)
\rput[t](.5,-.1){$n-1$} \psline[doubleline=true](.5,0)(.5,.4)
\psline[doubleline=true](.25,.6)(.25,1.4)\rput[l](.35,1){$n-2$}
\psline[doubleline=true](.5,1.6)(.5,2) \rput[b](.5,2.1){$n-1$}
\psarc[doubleline=true](1,.6){.2}{0}{180}
\psline[doubleline=true](1.2,0)(1.2,.6)
\psarc[doubleline=true](1,1.4){.2}{180}{0}
\psline[doubleline=true](1.2,1.4)(1.2,2)
\psframe[linecolor=darkred](0,.4)(1,.6)
\psframe[linecolor=darkred](0,1.4)(1,1.6)
\endpspicture
+a_{02} \pspicture[.45](-.2,-.3)(1.45,2.3) \rput[t](.5,-.1){$n-1$}
\psline[doubleline=true](.5,0)(.5,.4)
\psline[doubleline=true](.25,.6)(.25,1.4)
\psline[doubleline=true](.5,1.6)(.5,2) \rput[b](.5,2.1){$n-1$}
\psline[doubleline=true](.8,.6)(1.2,1.4)
\psline[doubleline=true](1.2,0)(1.2,.6)
\psline[doubleline=true](.8,1.4)(1.2,.6)
\psline[doubleline=true](1.2,1.4)(1.2,2)
\psframe[linecolor=darkred](0,.4)(1,.6)
\psframe[linecolor=darkred](0,1.4)(1,1.6)
\endpspicture
\label{b20n2exp}
\end{eqnarray}

\begin{cor}
The double clasp expansion of $B_2$ of type $(0,n)$ can be
obtained as in the equation~\ref{b20n2exp} where $a_{12}, a_{02}$
are from the Theorem~\ref{b20nexpthm}. \label{b20nexpcor}
\end{cor}

\newpage
\pagestyle{myheadings} \markright{  \rm \normalsize CHAPTER 3.
\hspace{0.5cm} Trihedron coefficients for
$\mathcal{U}_q(\mathfrak{sl}(3,\mathbb{C}))$}
\chapter{Trihedron Coefficients for $\mathcal{U}_q(\mathfrak{sl}(3,\mathbb{C}))$}
\thispagestyle{myheadings}

By~\cite{Lickorish:gtm}~\cite{MV:3valent}~\cite{Turaev:quantum},
we define a trivalent vertex as follow. A triple integers
$(a,b,c)$ is admissible if $a+b+c$ is even and $|a-b|\le c\le
a+b$. This is equivalent to the following. For
$\mathfrak{sl}(2,\mathbb{C})$, $\dim(\Inv(V_{a}\otimes
V_{b}\otimes V_{c}))$ is $1$ if $(a, b, c)$ is an admissible
triple or $0$ otherwise, where $V_a$ is an irreducible
representation of highest weight $a$. Given an admissible triple,
we define a trivalent vertex

$$
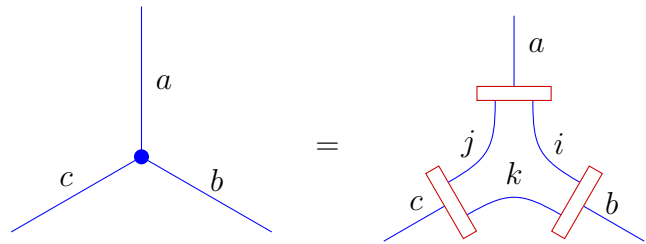
\pspicture[.4](-2.2,-1.4)(2.2,2.2) \qline(0,2)(0,0)
\qline(0,0)(1.732,-1) \qline(0,0)(-1.732,-1) \pscircle*(0,0){.1}
\rput(.3,1){$a$} \rput(1,-.3){$b$} \rput(-1,-.3){$c$}
\endpspicture
=\pspicture[.4](-2.2,-1.2)(2.2,2.2) \qline(0,2)(0,1.066)
\rput(.6,.3){$i$}\rput(-.6,.3){$j$}\rput(0,-.1){$k$}
\rput(.3,1.6){$a$} \rput(1.3,-.5){$b$} \rput(-1.3,-.5){$c$}
\pccurve[angleA=-90,angleB=30,ncurv=1](-.25,.866)(-.875,-.2165)
\pccurve[angleA=-90,angleB=150,ncurv=1](.25,.866)(.875,-.2165)
\pccurve[angleA=30,angleB=150,ncurv=1](-.625,-.6495)(.625,-.6495)
\qline(-.9232,-.533)(-1.732,-1) \qline(.9232,-.533)(1.732,-1)
\psframe[fillstyle=solid,fillcolor=white,linecolor=darkred](-.5,.866)(.5,1.066)
\pspolygon[fillstyle=solid,fillcolor=white,linecolor=darkred]
(.5,-.866)(.6732,-.966)(1.1732,-.1)(1,0)
\pspolygon[fillstyle=solid,fillcolor=white,linecolor=darkred]
(-.5,-.866)(-.6732,-.966)(-1.1732,-.1)(-1,0)
\endpspicture
$$

where $i=(a+b-c)/2$, $j=(a+c-b)/2$ and $(b+c-a)/2$. Then the
trihedron coefficient for $\mathcal{U}_q(\mathfrak{sl}(2,\mathbb{C}))$ is
known~\cite{Lickorish:gtm}~\cite{MV:3valent}~\cite{Turaev:quantum}
as

$$
\pspicture[.48](-1.2,-2.2)(1.2,2.2) \pscircle*(-1,0){.1}
\pscircle*(1,0){.1}
\pccurve[angleA=60,angleB=120,ncurv=1](-1,0)(1,0)
\pccurve[angleA=0,angleB=180,ncurv=1](-1,0)(1,0)
\pccurve[angleA=-60,angleB=-120,ncurv=1](-1,0)(1,0)
\rput(0,1){$a$}\rput(0,.25){$b$}\rput(0,-.89){$c$}
\endpspicture=
\pspicture[.48](-2.2,-2.2)(2.2,2.2) \pcline(-1.3,-.9)(-.3,-.1)
\pcline(1.3,.9)(.3,.1) \pcline(1.3,-.9)(.3,-.1)
\pcline(-1.3,.9)(-.3,.1)
\pccurve[angleA=225,angleB=315,ncurv=1](.7,.9)(-.7,.9)
\pccurve[angleA=45,angleB=135,ncurv=1](-.7,-.9)(.7,-.9)\rput(0,.8){$j$}
\rput(0,-.85){$j$} \qline(-1,-1.5)(-1,-1.1)\qline(-1,1.1)(-1,1.5)
\qline(1,-1.5)(1,-1.1) \qline(1,1.1)(1,1.5) \qline(1,-1.5)(2,-1.5)
\qline(2,-1.5)(2,1.5) \qline(2,1.5)(1,1.5) \qline(-1,1.5)(-2,1.5)
\qline(-2,1.5)(-2,-1.5) \qline(-2,-1.5)(-1,-1.5)
\rput[b](-1,.2){$i$}\rput[b](1,.2){$k$}
\psframe[linecolor=darkred](-1.5,-1.1)(-.5,-.9)
\psframe[linecolor=darkred]( 1.5,-1.1)( .5,-.9)
\psframe[linecolor=darkred](-1.5, 1.1)(-.5, .9)
\psframe[linecolor=darkred]( 1.5, 1.1)( .5, .9)
\psframe[linecolor=darkred](-.5,.1)(.5,-.1)
\endpspicture
=(-1)^{i+j+k}\frac{[i+j+k+1]![i]![j]![k]!}{[i+j]![j+k]![i+k]!}.
\label{a1triexp2}
$$

In previous chapter we found a recursive formula for generalized
Jones-Wenzl projectors. So we study how we generalize trihedron
coefficients to $\mathcal{U}_q(\mathfrak{sl}(3,\mathbb{C}))$. But,
the definition of the trivalent vertex is a little subtle. We will
prove the following statement in lemma~\ref{a2trilem11}. Let
$\lambda_1, \lambda_2$ be the fundamental dominant weights of
$\mathfrak{sl}(3,\mathbb{C})$ (mainly we will use
$\mathfrak{sl}(3,\mathbb{C})$-modules because it is known that
representation theories of $\mathfrak{sl}(3,\mathbb{C})$ and
$\mathcal{U}_q(\mathfrak{sl}(3,\mathbb{C}))$ are parallel: see
theorem~\ref{quantumtensor} and the representation theory of
$\mathfrak{sl}(3,\mathbb{C})$ is well known
in~\cite{Humphreys:gtm}). Let $V_{a\lambda_1 +b\lambda_2}$ be an
irreducible representation of $\mathfrak{sl}(3,\mathbb{C})$ of
highest weight $a\lambda_1 +b\lambda_2$. Each edge is decorated by
an irreducible representation of $\mathfrak{sl}(3,\mathbb{C})$,
denoted by $V_{a_1\lambda_1+b_1\lambda_2}$,
$V_{a_2\lambda_1+b_2\lambda_2}$ and
$V_{a_3\lambda_1+b_3\lambda_2}$ where $a_i, b_j$ are nonnegative
integers. Let $d=$ Min $\{a_1,a_2,a_3,b_1,b_2,b_3\}$. Then
$\dim(\Inv(V_{a_1\lambda_1+b_1\lambda_2}\otimes
V_{a_2\lambda_1+b_2\lambda_2}\otimes
V_{a_3\lambda_1+b_3\lambda_2}))$ is $d+1$ if there exist non
negative integers $i, j, k, l, m, n, o, p, q$ such that $a_1=i+p$,
$a_2=j+n$, $a_3=k+l$, $b_1=j+o$, $b_2=k+m$, $b_3=i+q$ and
$0=le^{\frac{\pi}{3}i}$ $+me^{\frac{\pi}{6}i}$ $+n+$
$oe^{-\frac{\pi}{6}i}$ $+pe^{-\frac{\pi}{3}i}$ $-q$. Otherwise, it
is zero.

If $\dim(\Inv(V_{a_1\lambda_1+b_1\lambda_2}\otimes
V_{a_2\lambda_1+b_2\lambda_2}\otimes
V_{a_3\lambda_1+b_3\lambda_2}))$ is nonzero, we say a triple of
ordered pairs $((a_1,b_1),(a_2,b_2),(a_3,b_3))$ is {\it
admissible}. It has shown that for a fixed boundary, there are
fillings which are cut outs from the the hexagonal tiling of the
plane~\cite{Kuperberg:spiders}. A general trihedron shape is given
in the following figure where $\alpha+\beta=a_1$, $i+j+k=b_1$,
$k+l=a_2$, $\beta+j+m=b_2$, $i+m=a_3$ and $\alpha+j+l=b_3$. The
top and the bottom part are actually the same after some
modifications which we will discuss later.

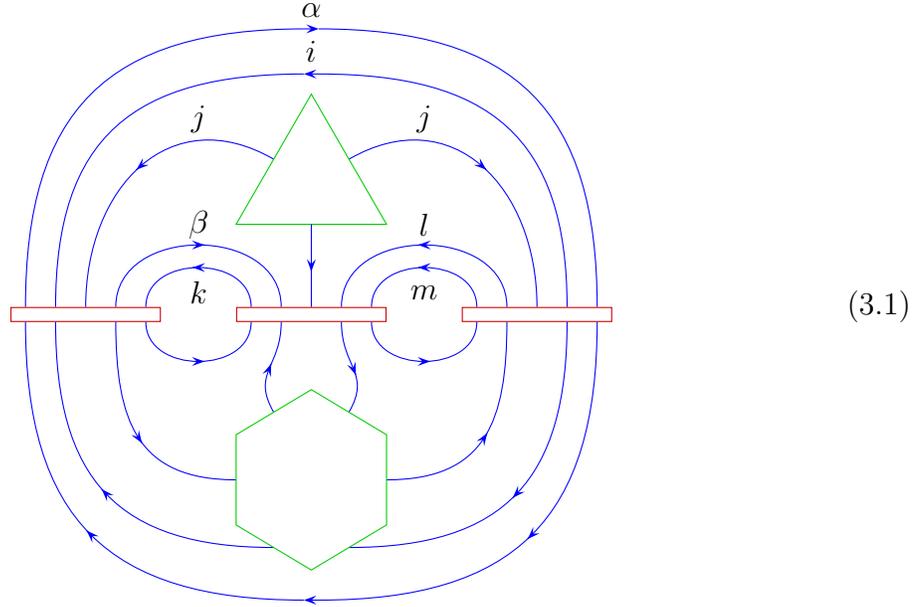
\begin{eqnarray}
\pspicture[.5](-4.1,-4.2)(4.1,4.2) \rput(0,4.05){$\alpha$}
\rput(-1.5,1.2){$\beta$} \rput(0,3.5){$i$} \rput(-1.5,2.6){$j$}
\rput(1.5,2.6){$j$} \rput(-1.5,.3){$k$} \rput(1.5,1.2){$l$}
\rput(1.5,.3){$m$} \psline[arrowscale=1.5]{->}(-.1,3.8)(.1,3.8)
\psline[arrowscale=1.5]{->}(.1,3.2)(-.1,3.2)
\pccurve[angleA=90,angleB=180,ncurv=1](-3.8,.1)(-.1,3.8)
\pccurve[angleA=90,angleB=180,ncurv=1](-3.4,.1)(-.1,3.2)
\pccurve[angleA=150,angleB=90,ncurv=1](-.5,2.066)(-3,.1)\middlearrow
\pccurve[angleA=90,angleB=90,ncurv=1](-.8,.1)(-2.2,.1)\middlearrow
\pccurve[angleA=90,angleB=90,ncurv=1](-2.6,.1)(-.4,.1)\middlearrow
\pccurve[angleA=90,angleB=90,ncurv=1](2.6,.1)(.4,.1)\middlearrow
\pccurve[angleA=90,angleB=90,ncurv=1](2.2,.1)(.8,.1)\middlearrow
\pccurve[angleA=30,angleB=90,ncurv=1](.5,2.066)(3,.1)\middlearrow
\pccurve[angleA=0,angleB=90,ncurv=1](.1,3.2)(3.4,.1)
\pccurve[angleA=0,angleB=90,ncurv=1](.1,3.8)(3.8,.1)
\pspolygon[linecolor=emgreen](1,1.2)(-1,1.2)(0,2.932)
\psline[arrowscale=1.5]{->}(.1,-3.8)(-.1,-3.8)
\pcline(0,1.2)(0,.1)\middlearrow
\pccurve[angleA=-180,angleB=-90,ncurv=1](-.1,-3.8)(-3.8,-.1)\middlearrow
\pccurve[angleA=-180,angleB=-90,ncurv=1](-.5,-3.1)(-3.4,-.1)\middlearrow
\pccurve[angleA=-90,angleB=180,ncurv=1](-2.6,-.1)(-1,-2.2)\middlearrow
\pccurve[angleA=-90,angleB=-90,ncurv=1](-2.2,-.1)(-.8,-.1)\middlearrow
\pccurve[angleA=120,angleB=-90,ncurv=1](-.5,-1.3)(-.4,-.1)\middlearrow
\pccurve[angleA=-90,angleB=60,ncurv=1](.4,-.1)(.5,-1.3)\middlearrow
\pccurve[angleA=-90,angleB=-90,ncurv=1](.8,-.1)(2.2,-.1)\middlearrow
\pccurve[angleA=0,angleB=-90,ncurv=1](1,-2.2)(2.6,-.1)\middlearrow
\pccurve[angleA=-90,angleB=0,ncurv=1](3.4,-.1)(.5,-3.1)\middlearrow
\pccurve[angleA=-90,angleB=0,ncurv=1](3.8,-.1)(.1,-3.8)\middlearrow
\pspolygon[linecolor=emgreen](0,-3.4)(1,-2.8)(1,-1.6)(0,-1)(-1,-1.6)(-1,-2.8)
\psframe[linecolor=darkred](-1,-.1)(1,.1)
\psframe[linecolor=darkred](-4,-.1)(-2,.1)
\psframe[linecolor=darkred](2,-.1)(4,.1)
\endpspicture
\label{thetadef}
\end{eqnarray}

So we can write trihedron coefficients as a $(d+1)\times (d+1)$
matrix. Let us denote it by $M_{\Theta}(a_1,b_1,a_2,b_2,a_3,b_3)$
or $M_{\Theta}(\lambda)$ where $\lambda
=(a_1,b_1,a_2,b_2,a_3,b_3)$. Also we denotes its $(i,j)$ entry by
$\Theta_{i,j}(a_1,b_1,a_2,b_2,a_3,b_3)$ or
$\Theta_{i,j}(\lambda)$. It is obvious that $M_{\Theta}(\lambda)$
is symmetric. Unfortunately the trihedron coefficient of this
shape is no longer rational expression composed of monomials
of quantum integers (if so, we say it to be {\it round}) in a
simple case $((1,1),(1,1),(1,1))$.
So we start to look the case $a_1=0$. Then we have found
the trihedron coefficients for the case $\alpha=\beta=0$ and
either $k=0$ or $j=0$ from the general shape.

\begin{thm}
$M_{\Theta}(0,i+j,l,j+m,i+m,j+l)$ is
$$(-1)^{j}[j+l+1][i+j+l+m+2]\frac{[i+j+m+1]![i]![j]![m]!}{[i+j]![j+m]![i+m]![2]}.$$
\label{onezerothm1}
\end{thm}

\begin{thm}
$M_{\Theta}(0,i+k,k+l,m,i+m,l)$ is
$$\sum_{n=0}^{\mathrm{min}\{ l,k,m\}}
a_n \frac{[i+l+1][i+l+m+2]}{[i+n+1][i+n+m+2]}
\frac{[i+k+m+2]![k-n+1]![i+m]![m+1]!}{[i+k]![i+m+n]![k+m-n+1]![2]}.$$
where $$a_n=
(-1)^n[k]^{2n}\frac{[k+l-n]![m]![k+m+l-n+1]!}{[k+l]![m-n]![k+m+l+1]!}.$$
\label{onezerothm2}
\end{thm}

In section $2$, we show all possible trihedron shapes and some
properties of trihedron shapes. In section $3$, we prove the main
theorems.

\section{General Shapes}

Now we will look at the trihedron coefficients.
The general shape is given in figure~(\ref{thetadef})
where the weight of clasps are $a_1\lambda_1+b_1\lambda_2$,
$a_2\lambda_1+b_2\lambda_2$ and $a_3\lambda_1+b_3\lambda_2$ and
$a_i, b_j$ are nonnegative integers.

For $\mathfrak{sl}(2,\mathbb{C})$, $\dim(\Inv(V_{i}\otimes
V_{j}\otimes V_{k}))$ is $1$ if $(i, j, k)$ is an admissible
triple or $0$ otherwise, where $V_i$ be an irreducible
representation of highest weight $i$. So there is a unique way to
fill in the triangle. But for $\mathfrak{sl}(3,\mathbb{C})$, we
could have more than one ways. Thus the shape of the polygon that
we are filling in might vary depending on the weight of clasps.
First we will discuss how we find a general shape in
figure~(\ref{thetadef}). Instead of finding a cut out from the
hexagonal tiling, we find a way to put three clasps into the
hexagonal tiling. Since it bounds a polygon and we have a set of
restrictions how clasps can be bent, we can change this problem to
an elementary geometry problem. Since all clasps are segregated,
until we reverse the direction of arrows, all possible interior
angles are $60, 180$ or $300$. Since an $Y$ makes the web zero, we
can exclude $60$. We would not count $180$ because it can be seen
as a subdivision. Thus it has to be $300$ if we do not change the
direction of arrows. When we change the direction of arrows, there
are also three possible interior angle either $0, 120$, $240$ or
$360$, let us denote this angle by $\alpha_i$. Since $U$ turn
makes the web zero, we can also exclude $0$. Let us denote the
angle between clasps by $\beta_i$. If there is a cut out bounded
by three given clasps of which all $a_i, b_j$ are nonzero, we can
see the following equality from the sum of interior angles.
$$180(n-2)= 300(n-6) +\sum_i^3 \alpha_i + \beta_i.$$
And we can simplify it to have $1440-120n=\sum \alpha_i+\beta_i$.
Since $\beta_i$ is either $120$ or $240$, $n$ can be either $6, 7$
or $8$. For each $n$, we look at the all possible combinations of
$\alpha_i, \beta_j$ up to symmetries. Then we check whether it
actually bounds a polygon. For example, $n=6$ there are $6$
possible combinations of angles but one does not bound a polygon,
$(\alpha_i)= (0,120,240), (\beta_i)=(120,120,120)$. For $n>6$, one
has to use $n-6$ times of $300$ angles. If some of $a_i, b_j$ are
zero then we can play the same game to find all shapes, we might
have $60$ for some $\beta_j$. The polygon in the middle might be a
triangle, a rectangle or a pentagon but we will consider them as a
special case of a hexagon. But we can obtain all these possible
shapes from the general shape by substituting some zeros. Then we
prove the following lemma. Let $d=$ Min
$\{a_1,a_2,a_3,b_1,b_2,b_3\}$.

\begin{lem}
$\dim(\Inv(V_{a_1\lambda_1+b_1\lambda_2}\otimes
V_{a_2\lambda_1+b_2\lambda_2}\otimes
V_{a_3\lambda_1+b_3\lambda_2}))$ is $d+1$ if there exist non
negative integers $i, j, k, l, m, n, o, p, q$ such that $a_1=i+p$,
$a_2=j+n$, $a_3=k+l$, $b_1=j+o$, $b_2=k+m$, $b_3=i+q$ and
$0=le^{\frac{\pi}{3}i}$ $+me^{\frac{\pi}{6}i}$ $+n+$
$oe^{-\frac{\pi}{6}i}$ $+pe^{-\frac{\pi}{3}i}$ $-q$. Otherwise, it
is zero.
\label{a2trilem11}
\end{lem}
\begin{proof}
One might directly find the answer by using the decomposition of
the tensor product $V_{a_1\lambda_1+b_1\lambda_2}\otimes
V_{a_2\lambda_1+b_2\lambda_2}\otimes
V_{a_3\lambda_1+b_3\lambda_2}$ into irreducible representations.
But one can see that these are conditions we can easily obtain
from the figure~(\ref{thetadef}) and the last one is one condition
that the hexagon in the middle does exist.
\end{proof}

Since the parallelogram in figure~(\ref{changeofdir}) changes the
directions of edges (the interior is filled by the unique maximal
hexagonal cut out), we can push the hexagon to an equilateral
triangle (possibly empty). The size of this equilateral triangle
is the minimum of differences of lengths of three pairs of
parallel edges of the hexagon.

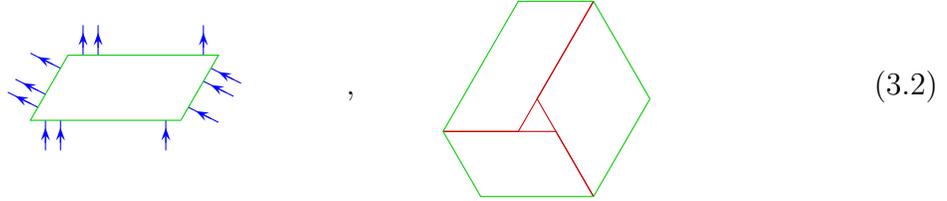
\begin{eqnarray}
\pspicture[.5](-.5,-.6)(3.2,1.3)
\psline(.2,-.4)(.2,0)\psline[arrowscale=1.5]{->}(.2,-.3)(.2,-.1)
\psline(.4,-.4)(.4,0)\psline[arrowscale=1.5]{->}(.4,-.3)(.4,-.1)
\psline(1.8,-.4)(1.8,0)\psline[arrowscale=1.5]{->}(1.8,-.3)(1.8,-.1)
\psline(.7,.866)(.7,1.266)\psline[arrowscale=1.5]{->}(.7,.966)(.7,1.166)
\psline(.9,.866)(.9,1.266)\psline[arrowscale=1.5]{->}(.9,.966)(.9,1.166)
\psline(2.3,.866)(2.3,1.266)\psline[arrowscale=1.5]{->}(2.3,.966)(2.3,1.166)
\psline(.1,.1732)(-.3,.3732)\psline[arrowscale=1.5]{->}(0,.2232)(-.2,.3232)
\psline(.2,.3463)(-.2,.5463)\psline[arrowscale=1.5]{->}(.1,.3963)(-.1,.4963)
\psline(.4,.6928)(0,.8928)\psline[arrowscale=1.5]{->}(.3,.7428)(.1,.8428)
\psline(2.1,.1732)(2.5,-.0268)\psline[arrowscale=1.5]{<-}(2.2,.1232)(2.4,.0232)
\psline(2.3,.5196)(2.7,.3196)\psline[arrowscale=1.5]{<-}(2.4,.4696)(2.6,.3696)
\psline(2.4,.6928)(2.8,.4928)\psline[arrowscale=1.5]{<-}(2.5,.6428)(2.7,.5428)
\pspolygon[linecolor=emgreen](0,0)(2,0)(2.5,.866)(.5,.866)
\endpspicture
\hskip 1cm , \hskip 1cm \pspicture[.5](-.1,-1)(3.2,2)
\psline[linecolor=darkred](0,0)(1.5,0)(2,-.866)
\psline[linecolor=darkred](2,-.866)(1.25,.433)(2,1.732)
\psline[linecolor=darkred](0,0)(1,0)(2,1.732)
\pspolygon[linecolor=emgreen](0,0)(.5,-.866)(2,-.866)(2.75,.433)(2,1.732)(1,1.732)
\endpspicture
\label{changeofdir}
\end{eqnarray}

But the resulting shape has mixed aspects. First it might contain
some non-segregated clasps which did not exist before.
When we apply single or double clasp
expansions to these non segregated clasps, it is possible to get
one by transforming the basis by adding $H$'s but usually it
becomes very difficult to deal with. If we keep the hexagonal
shape, it usually produces multiple non-vanishing terms in a
single clasp expansion. So we will use both shapes depend on the
feasibility.

\section{Proofs of Theorems}

Since $M_{\Theta}(\lambda)$ is an $1\times 1$ matrix, we will
write it as a scalar through this section. For nontrivial planar
$\Theta$ shapes with $d=0$, we could have at most three zeros for
$a_i, b_j$. Before we prove the main theorems we start with cases
of two zeros.

\subsection{Two Zeros}

If we have two zeros, up to symmetries, it is either one of these
subcases: 1) $a_1=b_3=0$, 2) $a_1=a_3=0$.

If $a_1=b_3=0$, there exists $j$ such that $b_1=a_2+j$,
$a_3=b_2+j$ and its shape is

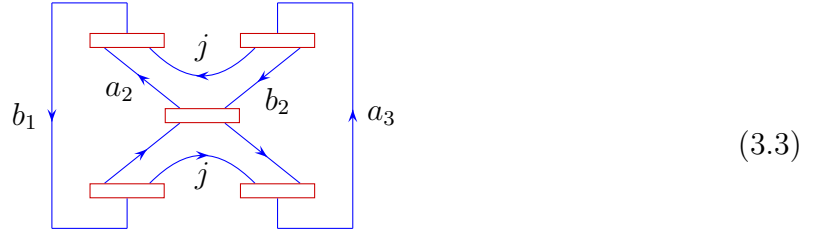
\begin{eqnarray}
\pspicture[.4](-2.9,-2.5)(2.9,2.5)
\pcline(-1.3,-.9)(-.3,-.1)\middlearrow
\pcline(1.3,.9)(.3,.1)\middlearrow
\pcline(.3,-.1)(1.3,-.9)\middlearrow
\pcline(-.3,.1)(-1.3,.9)\middlearrow
\pccurve[angleA=225,angleB=315,ncurv=1](.7,.9)(-.7,.9)\middlearrow
\pccurve[angleA=45,angleB=135,ncurv=1](-.7,-.9)(.7,-.9)\middlearrow
\psline(-1,-1.5)(-1,-1.1)\psline(-1,1.1)(-1,1.5)
\psline(1,-1.5)(1,-1.1) \psline(1,1.1)(1,1.5)
\psline(1,-1.5)(2,-1.5)
\pcline(-2,1.5)(-2,-1.5)\middlearrow\Bput{$b_1$}
\pcline(2,-1.5)(2,1.5)\middlearrow\Bput{$a_3$}
\psline(2,1.5)(1,1.5) \psline(-1,1.5)(-2,1.5)
\psline(-2,-1.5)(-1,-1.5) \rput(0,.9){$j$}
\rput[b](-1.1,.2){$a_2$} \rput(0,-.8){$j$} \rput[c](1,.2){$b_2$}
\psframe[linecolor=darkred](-1.5,-1.1)(-.5,-.9)
\psframe[linecolor=darkred]( 1.5,-1.1)( .5,-.9)
\psframe[linecolor=darkred](-1.5, 1.1)(-.5, .9)
\psframe[linecolor=darkred]( 1.5, 1.1)( .5, .9)
\psframe[linecolor=darkred](-.5,.1)(.5,-.1)
\endpspicture
\label{twozeroshape1}
\end{eqnarray}

By the same idea of the proof of proposition~\ref{a2abexpprop1},
we can easily get the following proposition~\ref{twozerothm1}. Let
$i=a_2$ and $k=b_2$.

\begin{prop}
$M_{\Theta}(0,i+j,i,k,j+k,0)$ is
$$\frac{[i+j+k+2]![i+1]![j]![k+1]!}
{[i+j]![j+k]![i+k+1]![2]}.$$
\label{twozerothm1}
\end{prop}

\begin{proof}
Let us use the notation $[i,j,k]$ for the trihedron coefficient of
this $\Theta$ shape. We induct on $k$. If $k=0$,
$$[i,j,0]=\frac{[i+j+2][i+j+1]}{[2]}.$$
For $k\neq 0$, by the idea of the proof of
proposition~\ref{a2abexpprop1} and the induction hypothesis, we
have
\begin{align*}
[i,j,k]&=\frac{[j+k+2]}{[j+k]}[i.j,k-1]-\frac{[i+1][i]}{[i+k+1][i+k]}[i-1,j+1,k-1]\\
&=\frac{[i+j+k+2]![i+1]![j]![k+1]!}{[i+j]![i+k+1]![j+k]!}(\frac{[j+k+2][i+k+1]-[j+1][i]}{[i+j+k+2][k+1]})\\
&=\frac{[i+j+k+2]![i+1]![j]![k+1]!}{[i+j]![i+k+1]![j+k]!}
\end{align*}
because $[j+1+k+1][i+k+1]=[j+1][i]+[i+j+k+2][k+1]$.
\end{proof}

If $a_1=a_3=0$, there exists $k$ such that $b_1=b_2+k$,
$a_2=b_3-b_2+k$ and its shape is

\begin{eqnarray}
\pspicture[.5](-2.1,-2.2)(2.1,2.2)
\pccurve[angleA=180,angleB=90,ncurv=1](-.5,1.1)(-1.3,.1)\middlearrow
\pccurve[angleA=90,angleB=90,ncurv=1](-.4,.1)(-1.1,.1)\middlearrow
\pccurve[angleA=90,angleB=-120,ncurv=1](-.2,.1)(-.25,.65)\middlearrow
\pccurve[angleA=-60,angleB=90,ncurv=1](.25,.65)(.2,.1)\middlearrow
\pccurve[angleA=30,angleB=90,ncurv=1](0,1.1)(1.7,.1)\middlearrow
\pspolygon[linecolor=emgreen](.5,.8)(0,.5)(-.5,.8)(-.5,1.4)
\pccurve[angleA=-90,angleB=180,ncurv=1](-1.3,-.1)(-.5,-1.1)\middlearrow
\pccurve[angleA=-90,angleB=-90,ncurv=1](-1.1,-.1)(-.4,-.1)\middlearrow
\pccurve[angleA=120,angleB=-90,ncurv=1](-.25,-.65)(-.2,-.1)\middlearrow
\pccurve[angleA=-90,angleB=60,ncurv=1](.2,-.1)(.25,-.65)\middlearrow
\pccurve[angleA=-90,angleB=-30,ncurv=1](1.7,-.1)(0,-1.1)\middlearrow
\pspolygon[linecolor=emgreen](.5,-.8)(0,-.5)(-.5,-.8)(-.5,-1.4)
\psframe[linecolor=darkred](-.5,-.1)(.5,.1)
\psframe[linecolor=darkred](-2,-.1)(-1,.1)
\psframe[linecolor=darkred](1,-.1)(2,.1)
\endpspicture
= \pspicture[.5](-2.1,-2)(2.1,2)
\psline(0,.1)(0,.5)\psline[arrowscale=1.5]{<-}(0,.2)(0,.4)
\psline(0,-.5)(0,-.1)\psline[arrowscale=1.5]{<-}(0,-.4)(0,-.2)
\pccurve[angleA=150,angleB=90,ncurv=1](-.375,1.15)(-1.8,.1)\middlearrow
\pccurve[angleA=-90,angleB=-150,ncurv=1](-1.8,-.1)(-.375,-1.15)\middlearrow
\pccurve[angleA=30,angleB=90,ncurv=1](.375,1.15)(1.8,.1)\middlearrow
\pccurve[angleA=-90,angleB=-30,ncurv=1](1.8,-.1)(.375,-1.15)\middlearrow
\pccurve[angleA=150,angleB=30,ncurv=1](-.2,.1)(-1.2,.1)\middlearrow
\pccurve[angleA=-30,angleB=-150,ncurv=1](-1.2,-.1)(-.2,-.1)\middlearrow
\pccurve[angleA=30,angleB=150,ncurv=1](.2,.1)(1.2,.1)\middlearrow
\pccurve[angleA=-150,angleB=-30,ncurv=1](1.2,-.1)(.2,-.1)\middlearrow
\pspolygon[linecolor=emgreen](0,1.8)(.75,.5)(-.75,.5)
\pspolygon[linecolor=emgreen](0,-1.8)(.75,-.5)(-.75,-.5)
\psframe[linecolor=darkred](-.5,-.1)(.5,.1)
\psframe[linecolor=darkred](-2,-.1)(-1,.1)
\psframe[linecolor=darkred](1,-.1)(2,.1)
\endpspicture
\label{twozeroshape2}
\end{eqnarray}

We need to prove a sequence of lemmas. First, we prove that the
clasp in the middle is not essential.

\begin{lem} Let $n\ge 1$, then
\begin{eqnarray}
\pspicture[.5](-2.1,-2)(2.1,2)
\rput(-1.6,1){$n$}\rput(1.6,1){$n$}\rput(.3,.3){$n$}
\psline(0,.1)(0,.5)\psline[arrowscale=1.5]{<-}(0,.2)(0,.4)
\psline(0,-.5)(0,-.1)\psline[arrowscale=1.5]{<-}(0,-.4)(0,-.2)
\pccurve[angleA=150,angleB=90,ncurv=1](-.375,1.15)(-1.5,.1)\middlearrow
\pccurve[angleA=-90,angleB=-150,ncurv=1](-1.5,-.1)(-.375,-1.15)\middlearrow
\pccurve[angleA=30,angleB=90,ncurv=1](.375,1.15)(1.5,.1)\middlearrow
\pccurve[angleA=-90,angleB=-30,ncurv=1](1.5,-.1)(.375,-1.15)\middlearrow
\pspolygon[linecolor=emgreen](0,1.8)(.75,.5)(-.75,.5)
\pspolygon[linecolor=emgreen](0,-1.8)(.75,-.5)(-.75,-.5)
\psframe[linecolor=darkred](-.5,-.1)(.5,.1)
\psframe[linecolor=darkred](-2,-.1)(-1,.1)
\psframe[linecolor=darkred](1,-.1)(2,.1)
\endpspicture
= \pspicture[.5](-2.1,-2)(2.1,2)
\rput(-1.6,1){$n$}\rput(1.6,1){$n$}\rput(.3,0){$n$}
\psline(0,-.5)(0,.5)\psline[arrowscale=1.5]{<-}(0,-.1)(0,.1)
\pccurve[angleA=150,angleB=90,ncurv=1](-.375,1.15)(-1.5,.1)\middlearrow
\pccurve[angleA=-90,angleB=-150,ncurv=1](-1.5,-.1)(-.375,-1.15)\middlearrow
\pccurve[angleA=30,angleB=90,ncurv=1](.375,1.15)(1.5,.1)\middlearrow
\pccurve[angleA=-90,angleB=-30,ncurv=1](1.5,-.1)(.375,-1.15)\middlearrow
\pspolygon[linecolor=emgreen](0,1.8)(.75,.5)(-.75,.5)
\pspolygon[linecolor=emgreen](0,-1.8)(.75,-.5)(-.75,-.5)
\psframe[linecolor=darkred](-2,-.1)(-1,.1)
\psframe[linecolor=darkred](1,-.1)(2,.1)
\endpspicture
= \pspicture[.5](-2.1,-2)(2.8,2)
\rput(-1.6,1){$n$}\rput(2,0){$n$}\rput(.3,0){$n$}
\psline(0,-.5)(0,.5)\psline[arrowscale=1.5]{<-}(0,-.1)(0,.1)
\pccurve[angleA=150,angleB=90,ncurv=1](-.375,1.15)(-1.5,.1)\middlearrow
\pccurve[angleA=-90,angleB=-150,ncurv=1](-1.5,-.1)(-.375,-1.15)\middlearrow
\psline(1,1.15)(.375,1.15)\psline[arrowscale=1.5]{<-}(.8,1.15)(.6,1.15)
\psline(1,-1.15)(.375,-1.15)\psline[arrowscale=1.5]{->}(.8,-1.15)(.6,-1.15)
\psarc(1.2,0.1){1.05}{0}{90} \psarc(1.2,-.1){1.05}{270}{0}
\psline[arrowscale=1.5]{<-}(2.25,-.1)(2.25,.1)
\pspolygon[linecolor=emgreen](0,1.8)(.75,.5)(-.75,.5)
\pspolygon[linecolor=emgreen](0,-1.8)(.75,-.5)(-.75,-.5)
\psframe[linecolor=darkred](1,.65)(1.2,1.65)
\psframe[linecolor=darkred](1,-1.65)(1.2,-.65)
\psframe[linecolor=darkred](-2,-.1)(-1,.1)
\endpspicture
\label{twozerolem1-1}
\end{eqnarray}
\label{twozerolem1}
\end{lem}

\begin{proof}
The idea of the proof is that if we have any $Y$'s in the single
clasp expansion at the middle clasp, it becomes zero. The
argument, we used to find the general shape, leads us that there
does not exist a filling with boundary $(0,n,0,n,1,n-2)$. Thus, it
has to vanish once we have any $Y$'s.
\end{proof}

Form the third figure in equation~(\ref{twozerolem1-1}), we apply
a single clasp expansion to the clasp in the left. For the
following equation, it should be clear without the direction of
edges.

\begin{lem} Let $n\ge 1$, then
\begin{eqnarray}
\pspicture[.5](-2.1,-2)(1.3,2) \rput(-1.6,1){$n$}\rput(.3,0){$n$}
\psline(0,-.5)(0,.5)\psline[arrowscale=1.5]{<-}(0,-.1)(0,.1)
\pccurve[angleA=150,angleB=90,ncurv=1](-.375,1.15)(-1.5,.1)\middlearrow
\pccurve[angleA=-90,angleB=-150,ncurv=1](-1.5,-.1)(-.375,-1.15)\middlearrow
\psline(1,1.15)(.375,1.15)\psline[arrowscale=1.5]{<-}(.8,1.15)(.6,1.15)
\psline(1,-1.15)(.375,-1.15)\psline[arrowscale=1.5]{->}(.8,-1.15)(.6,-1.15)
\pspolygon[linecolor=emgreen](0,1.8)(.75,.5)(-.75,.5)
\pspolygon[linecolor=emgreen](0,-1.8)(.75,-.5)(-.75,-.5)
\psframe[linecolor=darkred](1,.65)(1.2,1.65)
\psframe[linecolor=darkred](1,-1.65)(1.2,-.65)
\psframe[linecolor=darkred](-2,-.1)(-1,.1)
\endpspicture
= (-1)^{n}[n+1] \pspicture[.5](-.5,-2)(1.3,2) \rput(.4,0){$n$}
\pccurve[angleA=180,angleB=180,ncurv=1](1,-1.15)(1,1.15)\middlearrow
\psframe[linecolor=darkred](1,.65)(1.2,1.65)
\psframe[linecolor=darkred](1,-1.65)(1.2,-.65)
\endpspicture
\label{twozerolem3-2}
\end{eqnarray}
\label{twozerolem3}
\end{lem}

\begin{proof}
We induct on $n$. If $n=1$, the coefficient is $-[2]=(-1)^1[2]$.
We apply a single clasp expansion at the left clasp which gives us
the first equality in the following
equation~(\ref{twozerolem3-1}).

\begin{eqnarray}
&\pspicture[.5](-2.1,-2)(1.3,2) \rput(-1.7,1){$n$}
\psline(0,-.5)(0,.5)\psline[arrowscale=1.5]{<-}(0,-.1)(0,.1)
\pccurve[angleA=150,angleB=90,ncurv=1](-.375,1.15)(-1.5,.1)\middlearrow
\pccurve[angleA=-90,angleB=-150,ncurv=1](-1.5,-.1)(-.375,-1.15)\middlearrow
\psline(1,1.15)(.375,1.15)\psline[arrowscale=1.5]{<-}(.8,1.15)(.6,1.15)
\psline(1,-1.15)(.375,-1.15)\psline[arrowscale=1.5]{->}(.8,-1.15)(.6,-1.15)
\pspolygon[linecolor=emgreen](0,1.799)(.75,.5)(-.75,.5)
\pspolygon[linecolor=emgreen](0,-1.799)(.75,-.5)(-.75,-.5)
\psframe[linecolor=darkred](1,.65)(1.2,1.65)
\psframe[linecolor=darkred](1,-1.65)(1.2,-.65)
\psframe[linecolor=darkred](-2,-.1)(-1,.1)
\endpspicture
=\pspicture[.5](-2.5,-2)(1.3,2)\rput(-.9,.2){$n-1$}
\rput(-2,1.4){$1$}
\psline(0,-.5)(0,.5)\psline[arrowscale=1.5]{<-}(0,-.1)(0,.1)
\pccurve[angleA=150,angleB=90,ncurv=1](-.375,1.15)(-1.5,-.1)\middlearrow
\pccurve[angleA=-90,angleB=-150,ncurv=1](-1.5,-.3)(-.375,-1.15)\middlearrow
\pccurve[angleA=150,angleB=90,ncurv=1](-.125,1.5833)(-2.25,-.1)\middlearrow
\pccurve[angleA=-90,angleB=-150,ncurv=1](-2.25,-.3)(-.125,-1.5833)\middlearrow
\psline(-2.25,-.1)(-2.25,-.3)
\psline(1,1.15)(.375,1.15)\psline[arrowscale=1.5]{<-}(.8,1.15)(.6,1.15)
\psline(1,-1.15)(.375,-1.15)\psline[arrowscale=1.5]{->}(.8,-1.15)(.6,-1.15)
\pspolygon[linecolor=emgreen](0,1.8)(.75,.5)(-.75,.5)
\pspolygon[linecolor=emgreen](0,-1.8)(.75,-.5)(-.75,-.5)
\psframe[linecolor=darkred](1,.65)(1.2,1.65)
\psframe[linecolor=darkred](1,-1.65)(1.2,-.65)
\psframe[linecolor=darkred](-2,-.1)(-1,-.3)
\endpspicture +\frac{[n-1]}{[n]}
\pspicture[.5](-2.5,-2)(1.3,2)\rput(-.9,.2){$n-2$}
\psline(0,-.5)(0,.5)\psline[arrowscale=1.5]{<-}(0,-.1)(0,.1)
\pccurve[angleA=150,angleB=90,ncurv=1](-.375,1.15)(-1.5,-.1)\middlearrow
\pccurve[angleA=-90,angleB=-150,ncurv=1](-1.5,-.3)(-.375,-1.15)\middlearrow
\psline(1,1.15)(.375,1.15)\psline[arrowscale=1.5]{<-}(.8,1.15)(.6,1.15)
\psline(1,-1.15)(.375,-1.15)\psline[arrowscale=1.5]{->}(.8,-1.15)(.6,-1.15)
\psarc(-2,.1){.25}{0}{180} \psarc(-2,1.1){.25}{180}{0}
\psline(-2.25,-.1)(-2.25,.1) \psline(-2,.35)(-2,.85)
\qline(-1.75,-.1)(-1.75,.1)
\pccurve[angleA=-90,angleB=-150,ncurv=1](-2.25,-.1)(-.25,-1.3667)\middlearrow
\pccurve[angleA=150,angleB=90,ncurv=1](-.125,1.5833)(-2.25,1.1)\middlearrow
\pccurve[angleA=150,angleB=90,ncurv=1](-.25,1.3667)(-1.75,1.1)\middlearrow
\pspolygon[linecolor=emgreen](0,1.8)(.75,.5)(-.75,.5)
\pspolygon[linecolor=emgreen](0,-1.8)(.75,-.5)(-.75,-.5)
\psframe[linecolor=darkred](1,.65)(1.2,1.65)
\psframe[linecolor=darkred](1,-1.65)(1.2,-.65)
\psframe[linecolor=darkred](-2,-.1)(-1,-.3)
\endpspicture \\ \nonumber&= \pspicture[.5](-2.4,-2)(1.2,2)\rput(-.9,.2){$n-1$}
\rput(.8,0){$1$}
\psline(0.6,-.5)(0.6,.5)\psline[arrowscale=1.5]{<-}(0.6,-.1)(0.6,.1)
\psline(-.25,-.5)(-.25,.5)\psline[arrowscale=1.5]{<-}(-.25,-.1)(-.25,.1)
\pccurve[angleA=150,angleB=90,ncurv=1](-.375,1.15)(-1.5,-.1)\middlearrow
\pccurve[angleA=-90,angleB=-150,ncurv=1](-1.5,-.3)(-.375,-1.15)\middlearrow
\pccurve[angleA=150,angleB=90,ncurv=1](-.125,1.65)(-2.25,-.1)\middlearrow
\pccurve[angleA=-90,angleB=-150,ncurv=1](-2.25,-.3)(-.125,-1.65)\middlearrow
\psline(-2.25,-.3)(-2.25,-.1)
\qline(-.625,.55)(.24,1.15)\qline(-.625,-.55)(.24,-1.15)
\psline(1,1.15)(.375,1.15)\psline[arrowscale=1.5]{<-}(.8,1.15)(.6,1.15)
\psline(1,-1.15)(.375,-1.15)\psline[arrowscale=1.5]{->}(.8,-1.15)(.6,-1.15)
\pspolygon[linecolor=emgreen,fillcolor=white,fillstyle=solid](0,1.8)(-.1875,1.475)(0.375,.5)(.75,.5)
\pspolygon[linecolor=emgreen,fillcolor=white,fillstyle=solid](0.25,0.5)(-.25,1.3667)(-.75,.5)
\pspolygon[linecolor=emgreen,fillcolor=white,fillstyle=solid](0,-1.8)(-.1875,-1.5125)(0.4,-.5)(.75,-.5)
\pspolygon[linecolor=emgreen,fillcolor=white,fillstyle=solid](0.25,-0.5)(-.25,-1.3667)(-.75,-.5)
\psframe[linecolor=darkred](1,.65)(1.2,1.65)
\psframe[linecolor=darkred](1,-1.65)(1.2,-.65)
\psframe[linecolor=darkred](-2,-.1)(-1,-.3)
\endpspicture +\frac{[n-1]}{[n]}
\pspicture[.5](-2.4,-2)(1.2,2)\rput(-.05,0){$n-2$}\rput(.9,0){$2$}
\pccurve[angleA=-70,angleB=70,ncurv=1](0.5,.5)(0.5,-.5)\middlearrow
\pccurve[angleA=-105,angleB=105,ncurv=1](-.5625,.5)(-.5625,-.5)\middlearrow
\pccurve[angleA=150,angleB=90,ncurv=1](-.625,.7167)(-1,.1)\middlearrow
\pccurve[angleA=-90,angleB=-150,ncurv=1](-1,-.1)(-.625,-.7167)\middlearrow
\psline(1,1.15)(.375,1.15)\psline[arrowscale=1.5]{<-}(.8,1.15)(.6,1.15)
\psline(1,-1.15)(.375,-1.15)\psline[arrowscale=1.5]{->}(.8,-1.15)(.6,-1.15)
\psarc(-2,-.7){.25}{0}{180} \psarc(-2,.7){.25}{180}{0}
\psline(-2,-.45)(-2,.45) \qline(-.625,.55)(.24,1.15)
\qline(-.625,-.55)(.24,-1.15)
\pccurve[angleA=-90,angleB=-150,ncurv=1](-1.75,-.7)(-.25,-1.3667)\middlearrow
\pccurve[angleA=-90,angleB=-150,ncurv=1](-2.25,-.7)(-.125,-1.5833)\middlearrow
\pccurve[angleA=150,angleB=90,ncurv=1](-.125,1.5833)(-2.25,.7)\middlearrow
\pccurve[angleA=150,angleB=90,ncurv=1](-.25,1.3667)(-1.75,.7)\middlearrow
\pspolygon[linecolor=emgreen,fillcolor=white,fillstyle=solid](0,1.8)(-.5,.9333)(-.25,.5)(.75,.5)
\pspolygon[linecolor=emgreen,fillcolor=white,fillstyle=solid](-.375,.5)(-.5625,.825)(-.75,.5)
\pspolygon[linecolor=emgreen,fillcolor=white,fillstyle=solid](0,-1.8)(-.5,-.9333)(-.25,-.5)(.75,-.5)
\pspolygon[linecolor=emgreen,fillcolor=white,fillstyle=solid](-.375,-0.5)(-.5625,-.825)(-.75,-.5)
\psframe[linecolor=darkred](1,.65)(1.2,1.65)
\psframe[linecolor=darkred](1,-1.65)(1.2,-.65)
\psframe[linecolor=darkred](-1.4,-.1)(-.8,.1)
\endpspicture
+\frac{[n-2]}{[n]}
\pspicture[.5](-2.9,-2)(1.2,2)\rput(-.05,0){$n-2$}\rput(.9,0){$2$}
\pccurve[angleA=-70,angleB=70,ncurv=1](0.5,.5)(0.5,-.5)\middlearrow
\pccurve[angleA=-105,angleB=105,ncurv=1](-.5625,.5)(-.5625,-.5)\middlearrow
\pccurve[angleA=150,angleB=90,ncurv=1](-.625,.7167)(-1.3,.1)\middlearrow
\pccurve[angleA=-90,angleB=-150,ncurv=1](-1.3,-.1)(-.6875,-.6083)\middlearrow
\psline(1,1.15)(.375,1.15)\psline[arrowscale=1.5]{<-}(.8,1.15)(.6,1.15)
\psline(1,-1.15)(.375,-1.15)\psline[arrowscale=1.5]{->}(.8,-1.15)(.6,-1.15)
\psarc(-2.5,.1){.25}{0}{180} \psarc(-2.5,1.1){.25}{180}{0}
\psline(-2.75,-.1)(-2.75,.1) \psline(-2.75,-.1)(-2.75,.1)
\psline(-2.5,.35)(-2.5,.85) \psarc(-2,-.1){.25}{180}{0}
\psarc(-2,-1.1){.25}{0}{180} \psline(-2.25,-.1)(-2.25,.1)
\psline(-2,-.35)(-2,-.85) \qline(-.625,.55)(.24,1.15)
\qline(-.625,-.55)(.24,-1.15)
\pccurve[angleA=-90,angleB=-150,ncurv=1](-1.75,-1.1)(-.625,-.7167)\middlearrow
\pccurve[angleA=-90,angleB=-150,ncurv=1](-2.25,-1.1)(-.25,-1.4)\middlearrow
\qline(-2.75,-.1)(-2.75,-1.1)
\pccurve[angleA=-90,angleB=-150,ncurv=1](-2.75,-1.1)(-.125,-1.5833)\middlearrow
\pccurve[angleA=150,angleB=90,ncurv=1](-.125,1.5833)(-2.75,1.1)\middlearrow
\pccurve[angleA=150,angleB=90,ncurv=1](-.25,1.3667)(-2.25,1.1)\middlearrow
\pspolygon[linecolor=emgreen,fillcolor=white,fillstyle=solid](0,1.8)(-.5,.9333)(-.25,.5)(.75,.5)
\pspolygon[linecolor=emgreen,fillcolor=white,fillstyle=solid](-.375,.5)(-.5625,.825)(-.75,.5)
\pspolygon[linecolor=emgreen,fillcolor=white,fillstyle=solid](0,-1.8)(-.5,-.9333)(-.25,-.5)(.75,-.5)
\pspolygon[linecolor=emgreen,fillcolor=white,fillstyle=solid](-.375,-0.5)(-.5625,-.825)(-.75,-.5)
\psframe[linecolor=darkred](1,.65)(1.2,1.65)
\psframe[linecolor=darkred](1,-1.65)(1.2,-.65)
\psframe[linecolor=darkred](-2,-.1)(-1,.1)
\endpspicture \\\nonumber
&=
\alpha\pspicture[.5](-1.5,-2)(1.3,2)\rput(-.55,0){$n-1$}\rput(.8,0){$1$}
\psline(0.6,-.5)(0.6,.5)\psline[arrowscale=1.5]{<-}(0.6,-.1)(0.6,.1)
\pccurve[angleA=150,angleB=90,ncurv=1](-.125,1.5833)(-1.25,0)\middlearrow
\pccurve[angleA=-90,angleB=-150,ncurv=1](-1.25,0)(-.125,-1.5833)\middlearrow
\psline(1,1.15)(.375,1.15)\psline[arrowscale=1.5]{<-}(.8,1.15)(.6,1.15)
\psline(1,-1.15)(.375,-1.15)\psline[arrowscale=1.5]{->}(.8,-1.15)(.6,-1.15)
\pccurve[angleA=-100,angleB=100,ncurv=1](.1,1)(.1,-1)\middlearrow
\pspolygon[linecolor=emgreen,fillcolor=white,fillstyle=solid](0,1.8)(-.1875,1.475)(0.375,.5)(.75,.5)
\pspolygon[linecolor=emgreen,fillcolor=white,fillstyle=solid](0,-1.8)(-.1875,-1.5125)(0.4,-.5)(.75,-.5)
\psframe[linecolor=darkred](1,.65)(1.2,1.65)
\psframe[linecolor=darkred](1,-1.65)(1.2,-.65)
\endpspicture +\beta\frac{[n-1]}{[n]}
\pspicture[.5](-2.5,-2)(1.3,2)\rput(-1.1,0){$n-2$}\rput(.8,0){$2$}
\psline(0.5,-.5)(0.5,.5)\psline[arrowscale=1.5]{<-}(0.5,-.1)(0.5,.1)
\psline(1,1.15)(.375,1.15)\psline[arrowscale=1.5]{<-}(.8,1.15)(.6,1.15)
\psline(1,-1.15)(.375,-1.15)\psline[arrowscale=1.5]{->}(.8,-1.15)(.6,-1.15)
\psarc(-2,-.7){.25}{0}{180} \psarc(-2,.7){.25}{180}{0}
\psline(-2,-.45)(-2,.45)
\pccurve[angleA=-90,angleB=-150,ncurv=1](-1.75,-.7)(-.25,-1.3667)\middlearrow
\pccurve[angleA=-90,angleB=-150,ncurv=1](-2.25,-.7)(-.125,-1.5833)\middlearrow
\pccurve[angleA=150,angleB=90,ncurv=1](-.125,1.5833)(-2.25,.7)\middlearrow
\pccurve[angleA=150,angleB=90,ncurv=1](-.25,1.3667)(-1.75,.7)\middlearrow
\pccurve[angleA=-100,angleB=100,ncurv=1](-.375,.85)(-.375,-.85)\middlearrow
\pspolygon[linecolor=emgreen,fillcolor=white,fillstyle=solid](0,1.8)(-.5,.9333)(-.25,.5)(.75,.5)
\pspolygon[linecolor=emgreen,fillcolor=white,fillstyle=solid](0,-1.8)(-.5,-.9333)(-.25,-.5)(.75,-.5)
\psframe[linecolor=darkred](1,.65)(1.2,1.65)
\psframe[linecolor=darkred](1,-1.65)(1.2,-.65)
\endpspicture + \frac{[n-2]}{[n]}
\pspicture[.5](-2.5,-2)(1.3,2)\rput(-.9,1.2){$n-2$}
\psline(0,-.2)(0,.2)\psline[arrowscale=1.5]{<-}(0,-.1)(0,.1)
\pccurve[angleA=150,angleB=90,ncurv=1](-.375,.85)(-1.5,.1)\middlearrow
\pccurve[angleA=-90,angleB=-150,ncurv=1](-1.5,-.1)(-.375,-.85)\middlearrow
\pccurve[angleA=90,angleB=180,ncurv=1](-2.25,.1)(0,1.9)\middlearrow
\pccurve[angleA=180,angleB=-90,ncurv=1](0,-1.9)(-2.25,-.1)\middlearrow
\pccurve[angleA=0,angleB=180,ncurv=1](0,1.9)(1.2,1.5)
\pccurve[angleA=180,angleB=0,ncurv=1](1.2,-1.5)(0,-1.9)
\psline(-2.25,.1)(-2.25,-.1) \psline(1.2,.8)(1,.8)(1,-.8)(1.2,-.8)
\psline[arrowscale=1.5]{->}(1,-.1)(1,.1)
\psline[angleA=30,angleB=180,ncurv=1](.375,.85)(1.2,1.175)\middlearrow
\psline[angleA=180,angleB=-30,ncurv=1](1.2,-1.175)(.375,-.85)\middlearrow
\pspolygon[linecolor=emgreen](0,1.5)(.75,.2)(-.75,.2)
\pspolygon[linecolor=emgreen](0,-1.5)(.75,-.2)(-.75,-.2)
\psframe[linecolor=darkred](1.2,.65)(1.4,1.65)
\psframe[linecolor=darkred](1.2,-1.65)(1.4,-.65)
\psframe[linecolor=darkred](-2,-.1)(-1,.1)
\endpspicture\nonumber
\label{twozerolem3-1}
\end{eqnarray}
where $\alpha=(-1)^{n-1}[n]$ and $\beta=(-1)^{n-2}[n-1]$.

For the second equality, we apply a single clasp expansion for the
resulting clasp of weight $(n-1)\lambda_2$ for which clasps of
weight $(n-2)\lambda_2$ are located at the northeast corner. It is not
difficult to see the following equality. The number of strings
coming from the trapezoid can be any integer between $0$ and $n$
where $n$ is the weight of the clasp given the left side of
equality.

$$
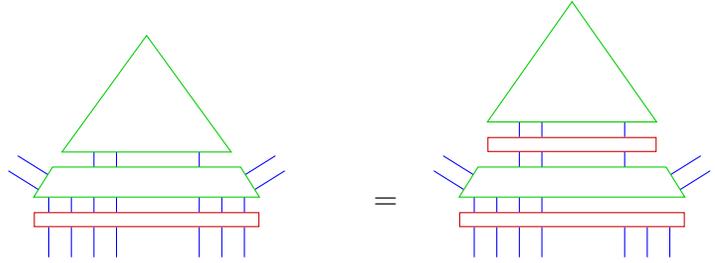
\pspicture[.2](-1.9,-.5)(2.9,3.2) \psline(-1.3,0)(-1.3,-.4)
\psline(-1,0)(-1,-.4) \psline(-.7,0)(-.7,-.4)
\psline(-.4,0)(-.4,-.4) \psline(.7,0)(.7,-.4) \psline(1,0)(1,-.4)
\psline(1.3,0)(1.3,-.4) \psline(-1.3,0.2)(-1.3,.4)
\psline(-1,0.2)(-1,.4) \psline(-.7,0.2)(-.7,.4)
\psline(-.4,0.2)(-.4,.4) \psline(.7,0.2)(.7,.4)
\psline(1,0.2)(1,.4) \psline(1.3,0.2)(1.3,.4)
\psline(-.7,1)(-.7,.8) \psline(-.4,1)(-.4,.8) \psline(.7,1)(.7,.8)
\psline(-1.4375,.5)(-1.8375,.75) \psline(-1.3125,.7)(-1.7125,.95)
\psline(1.4375,.5)(1.8375,.75) \psline(1.3125,.7)(1.7125,.95)
\pspolygon[linecolor=emgreen](-1.5,.4)(1.5,.4)(1.25,.8)(-1.25,.8)
\pspolygon[linecolor=emgreen](-1.125,1)(1.125,1)(0,2.55)
\psframe[linecolor=darkred](-1.5,0)(1.5,.2)
\endpspicture
= \pspicture[.2](-2.2,-.5)(1.9,3.2) \psline(-1.3,0)(-1.3,-.4)
\psline(-1,0)(-1,-.4) \psline(-.7,0)(-.7,-.4)
\psline(-.4,0)(-.4,-.4) \psline(.7,0)(.7,-.4) \psline(1,0)(1,-.4)
\psline(1.3,0)(1.3,-.4) \psline(-1.3,0.2)(-1.3,.4)
\psline(-1,0.2)(-1,.4) \psline(-.7,0.2)(-.7,.4)
\psline(-.4,0.2)(-.4,.4) \psline(-.7,1)(-.7,.8)
\psline(-.4,1)(-.4,.8) \psline(.7,1)(.7,.8)
\psline(-.7,1.2)(-.7,1.4) \psline(-.4,1.2)(-.4,1.4)
\psline(.7,1.2)(.7,1.4) \psline(-1.4375,.5)(-1.8375,.75)
\psline(-1.3125,.7)(-1.7125,.95) \psline(1.4375,.5)(1.8375,.75)
\psline(1.3125,.7)(1.7125,.95)
\pspolygon[linecolor=emgreen](-1.5,.4)(1.5,.4)(1.25,.8)(-1.25,.8)
\pspolygon[linecolor=emgreen](-1.125,1.4)(1.125,1.4)(0,3)
\psframe[linecolor=darkred](-1.5,0)(1.5,.2)
\psframe[linecolor=darkred](-1.125,1)(1.125,1.2)
\endpspicture.
$$

Thus, we can put a clasp of weight $(n-1)\lambda_2((n-2)\lambda_2, (n-2)\lambda_2)$ at the
gap between an equilateral triangle and a trapezoid at the first
(second and third, respectively) figure in the second line of
equation~(\ref{twozerolem3-1}). Therefore, we can use induction
to get the third equality in
equation~(\ref{twozerolem3-1}). Note that the size of the third equilateral
triangle in the third line is $n-2$.

Last step is to count how many $-[2]$'s will be produced when we
change it to multiple of the web in the right hand side of
equation~(\ref{twozerolem3-2}). But it is fairly easy to see that
each of them has just one factor of $-[2]$ in the first two in the
third line. If we add up the all coefficients, we have

\begin{eqnarray}
&-[2](-1)^{n-1}[n]-[2](-1)^{n-2}\frac{[n-1][n-1]}{[n]}+
(-1)^{n-2}\frac{[n-1][n-2]}{[n]}\nonumber\\
&=(-1)^{n}\frac{1}{[n]}([2][n]^2-[2][n-1]^2+[n-2][n-1])\nonumber\\
&=(-1)^{n}\frac{1}{[n]}[n+1][n]=(-1)^{n}[n+1]\nonumber
\end{eqnarray}
\end{proof}

The lemma~\ref{twozerolem3} can be generalized to the following
lemma.

\begin{lem}
\begin{eqnarray}
\pspicture[.5](-3.1,-2.5)(1.8,2.5)
\rput(.8,2.3){$i$}\rput(.8,1.5){$j$}
\psline(-.5,-.5)(-.5,-.1)\psline[arrowscale=1.5]{<-}(-.5,-.4)(-.5,-.2)
\psline(-.5,.5)(-.5,.1)\psline[arrowscale=1.5]{->}(-.5,.4)(-.5,.2)
\pccurve[angleA=150,angleB=90,ncurv=1](-.875,1.1495)(-2.2,.1)\middlearrow
\pccurve[angleA=-90,angleB=-150,ncurv=1](-2.2,-.1)(-.875,-1.1495)\middlearrow
\pccurve[angleA=180,angleB=90,ncurv=1](-.5,2)(-2.8,.1)\middlearrow
\pccurve[angleA=-90,angleB=-180,ncurv=1](-2.8,-.1)(-.5,-2)\middlearrow
\psline(-.5,2)(1.5,2) \psline(-.5,-2)(1.5,-2)
\psline(1.5,1.1495)(-.125,1.1495)\psline[arrowscale=1.5]{<-}(.8,1.1495)(.6,1.1495)
\psline(1.5,-1.1495)(-.125,-1.1495)\psline[arrowscale=1.5]{->}(.8,-1.1495)(.6,-1.1495)
\pspolygon[linecolor=emgreen](-.5,1.799)(.25,.5)(-1.25,.5)
\pspolygon[linecolor=emgreen](-.5,-1.799)(.25,-.5)(-1.25,-.5)
\psframe[linecolor=darkred](1.5,.8)(1.7,2.35)
\psframe[linecolor=darkred](1.5,-2.35)(1.7,-.8)
\psframe[linecolor=darkred](-3,-.1)(-2,.1)
\psframe[linecolor=darkred](-1.2,-.1)(.2,.1)
\endpspicture
=(-1)^{j}\frac{[i+j+1]}{[i+1]} \pspicture[.5](-.5,-2.5)(1.8,2.5)
\pccurve[angleA=180,angleB=180,ncurv=1](1.5,-1.35)(1.5,1.35)\middlearrow
\pccurve[angleA=180,angleB=180,ncurv=1](1.5,2)(1.5,-2)\middlearrow
\psframe[linecolor=darkred](1.5,1)(1.7,2.35)
\psframe[linecolor=darkred](1.5,-2.35)(1.7,-1)
\endpspicture
\label{twozerolem4-1}
\end{eqnarray}
\label{twozerolem4}
\end{lem}

\begin{proof}
We induct on $i+j$. If $i=1,j=0$, the coefficient is
$1=(-1)^1\frac{[2]}{[2]}$. If $i=0$, it follows from the previous
lemma. The first equality can be proven with a single clasp
expansion at the left middle clasp in left figure. The second
equality can be proven the same argument we use in lemma. The rest
of proof follows by induction.
\end{proof}

\begin{lem}
\begin{eqnarray}
\pspicture[.5](-3.1,-2.5)(1.8,2.5)
\rput(.5,2.3){$i$}\rput(.5,1.5){$j$}\rput(1.1,.5){$k$}
\psline(-.5,-.5)(-.5,-.1)\psline[arrowscale=1.5]{<-}(-.5,-.4)(-.5,-.2)
\psline(-.5,.5)(-.5,.1)\psline[arrowscale=1.5]{->}(-.5,.4)(-.5,.2)
\pccurve[angleA=150,angleB=90,ncurv=1](-.875,1.1495)(-2.2,.1)\middlearrow
\pccurve[angleA=-90,angleB=-150,ncurv=1](-2.2,-.1)(-.875,-1.1495)\middlearrow
\pccurve[angleA=180,angleB=90,ncurv=1](1.5,.7)(.5,.1)\middlearrow
\pccurve[angleA=-90,angleB=-180,ncurv=1](.5,-.1)(1.5,-.7)\middlearrow
\pccurve[angleA=180,angleB=90,ncurv=1](-.5,2)(-2.8,.1)\middlearrow
\pccurve[angleA=-90,angleB=-180,ncurv=1](-2.8,-.1)(-.5,-2)\middlearrow
\psline(-.5,2)(1.5,2) \psline(-.5,-2)(1.5,-2)
\psline(1.5,1.1495)(-.125,1.1495)\psline[arrowscale=1.5]{<-}(.8,1.1495)(.6,1.1495)
\psline(1.5,-1.1495)(-.125,-1.1495)\psline[arrowscale=1.5]{->}(.8,-1.1495)(.6,-1.1495)
\pspolygon[linecolor=emgreen](-.5,1.799)(.25,.5)(-1.25,.5)
\pspolygon[linecolor=emgreen](-.5,-1.799)(.25,-.5)(-1.25,-.5)
\psframe[linecolor=darkred](1.5,.35)(1.7,2.35)
\psframe[linecolor=darkred](1.5,-2.35)(1.7,-.35)
\psframe[linecolor=darkred](-3,-.1)(-2,.1)
\psframe[linecolor=darkred](-.7,-.1)(.7,.1)
\endpspicture
=(-1)^{j+1}\frac{[i+j+k+1]![i]![j]![k]!}{[i+j]![j+k]![i+k+1]!}
\pspicture[.5](-.5,-2.5)(1.8,2.5)
\pccurve[angleA=180,angleB=180,ncurv=1](1.5,.7)(1.5,-.7)\middlearrow
\pccurve[angleA=180,angleB=180,ncurv=1](1.5,-1.35)(1.5,1.35)\middlearrow
\pccurve[angleA=180,angleB=180,ncurv=1](1.5,2)(1.5,-2)\middlearrow
\psframe[linecolor=darkred](1.5,.35)(1.7,2.35)
\psframe[linecolor=darkred](1.5,-2.35)(1.7,-.35)
\endpspicture
\label{twozerolem5-1}
\end{eqnarray}
\label{twozerolem5}
\end{lem}
\begin{proof}
We induct on $k$. If $k=0$, it follows from the previous lemma. If
$k\neq 0$, we use a single clasp expansion at the left middle
clasp. Then we get the following equality.

By induction, the coefficient is equal to
\begin{align*}
&(-1)^{j}\frac{[i+j+k]![i]![j]![k-1]!}{[i+j]![i+k]![j+k-1]!}\\
&+\frac{[j][j]}{[j+k][j+k-1]}(-1)^{j-1}\frac{[i+j+k]![i+1]![j-1]![k-1]!}
{[i+j]![i+j+1]![j+k-2]!}\nonumber\\
&=(-1)^j\frac{[i+j+k+1]![i]![j]![k]!}{[i+j]![j+k]![i+k+1]!}
(\frac{[i+k+1][j+k]-[i+1][j]}{[i+j+k+1][k]})\\
&=(-1)^j\frac{[i+j+k+1]![i]![j]![k]!}{[i+j]![j+k]![i+k+1]!}
\end{align*}
because $[i+1+k][j+k]=[i+1][j]+[i+j+k+1][k]$.
\end{proof}

\begin{prop}
$M_{\Theta}(0,b_1,a_2,b_2,0,b_3)$ is
$$(-1)^{b_2+1}\frac{[b_1-b_2+b_3+2]![b_1-b_2]![b_3-b_2]![b_2+1]!}
{[b_1]![b_1-2b_2+b_3]![b_3]![2]}.$$ or
$M_{\Theta}(0,i+j,i+k,j,0,j+k)$ is
$$(-1)^{j}\frac{[i+j+k+2]![i]![j+1]![k]!}
{[i+j]![j+k]![i+k]![2]}.$$ \label{twozeroprop2}
\end{prop}

\begin{cor}
$M_{\Theta}(0,n,0,n,0,n)$ is
$$(-1)^{n}\frac{[n+1]^2[n+2]}
{[2]}.$$ \label{twozerocor1}
\end{cor}

\subsection{Proof of Theorem~\ref{onezerothm1}}

We start to prove the following lemma to find the trihedron
coefficient.

\begin{lem}
\begin{eqnarray}
\pspicture[.5](-3.1,-2.5)(2,2.5)
\rput(.5,2.3){$i$}\rput(.5,1.5){$j$}\rput(1.2,.3){$l$}
\rput(1.2,.8){$m$}
\psline(-.5,-.5)(-.5,-.1)\psline[arrowscale=1.5]{<-}(-.5,-.4)(-.5,-.2)
\psline(-.5,.5)(-.5,.1)\psline[arrowscale=1.5]{->}(-.5,.4)(-.5,.2)
\pccurve[angleA=150,angleB=90,ncurv=1](-.875,1.1495)(-2.2,.1)\middlearrow
\pccurve[angleA=-90,angleB=-150,ncurv=1](-2.2,-.1)(-.875,-1.1495)\middlearrow
\pccurve[angleA=180,angleB=90,ncurv=1](1.5,1)(.4,.1)\middlearrow
\pccurve[angleA=-90,angleB=180,ncurv=1](.4,-.1)(1.5,-1)\middlearrow
\pccurve[angleA=90,angleB=180,ncurv=1](.8,.1)(1.5,.5)\middlearrow
\pccurve[angleA=180,angleB=-90,ncurv=1](1.5,-.5)(.8,-.1)\middlearrow
\pccurve[angleA=180,angleB=90,ncurv=1](-.5,2)(-2.8,.1)\middlearrow
\pccurve[angleA=-90,angleB=-180,ncurv=1](-2.8,-.1)(-.5,-2)\middlearrow
\psline(-.5,2)(1.5,2) \psline(-.5,-2)(1.5,-2)
\pcline(-.125,1.1495)(1.5,1.5)\middlearrow
\pcline(1.5,-1.5)(-.125,-1.1495)\middlearrow
\pspolygon[linecolor=emgreen](-.5,1.799)(.25,.5)(-1.25,.5)
\pspolygon[linecolor=emgreen](-.5,-1.799)(.25,-.5)(-1.25,-.5)
\psframe[linecolor=darkred](1.5,.35)(1.7,2.35)
\psframe[linecolor=darkred](1.5,-2.35)(1.7,-.35)
\psframe[linecolor=darkred](-3,-.1)(-1.8,.1)
\psframe[linecolor=darkred](-1.4,-.1)(1,.1)
\endpspicture
= \pspicture[.5](-3.1,-2.5)(2,2.5)
\rput(.5,2.3){$i$}\rput(.5,1.5){$j$}\rput(1.2,.3){$l$}
\rput(1.2,.8){$m$}
\psline(-.5,-.5)(-.5,-.1)\psline[arrowscale=1.5]{<-}(-.5,-.4)(-.5,-.2)
\psline(-.5,.5)(-.5,.1)\psline[arrowscale=1.5]{->}(-.5,.4)(-.5,.2)
\pccurve[angleA=150,angleB=90,ncurv=1](-.875,1.1495)(-2.2,.1)\middlearrow
\pccurve[angleA=-90,angleB=-150,ncurv=1](-2.2,-.1)(-.875,-1.1495)\middlearrow
\pccurve[angleA=180,angleB=90,ncurv=1](1.5,1)(.4,.1)\middlearrow
\pccurve[angleA=-90,angleB=180,ncurv=1](.4,-.1)(1.5,-1)\middlearrow
\pccurve[angleA=90,angleB=180,ncurv=1](.8,.1)(1.5,.5)\middlearrow
\pccurve[angleA=180,angleB=-90,ncurv=1](1.5,-.5)(.8,-.1)\middlearrow
\pccurve[angleA=180,angleB=90,ncurv=1](-.5,2)(-2.8,.1)\middlearrow
\pccurve[angleA=-90,angleB=-180,ncurv=1](-2.8,-.1)(-.5,-2)\middlearrow
\psline(-.5,2)(1.5,2) \psline(-.5,-2)(1.5,-2)
\psline(.8,-.1)(.8,.1)  \pcline(-.125,1.1495)(1.5,1.5)\middlearrow
\pcline(1.5,-1.5)(-.125,-1.1495)\middlearrow
\pspolygon[linecolor=emgreen](-.5,1.799)(.25,.5)(-1.25,.5)
\pspolygon[linecolor=emgreen](-.5,-1.799)(.25,-.5)(-1.25,-.5)
\psframe[linecolor=darkred](1.5,.35)(1.7,2.35)
\psframe[linecolor=darkred](1.5,-2.35)(1.7,-.35)
\psframe[linecolor=darkred](-3,-.1)(-1.8,.1)
\psframe[linecolor=darkred](-1.4,-.1)(.6,.1)
\endpspicture
\label{onezerolem1-1}
\end{eqnarray}
\label{onezerolem1}
\end{lem}
\begin{proof}
We use an expansion in equation~(\ref{a2abexp1}) at the right
middle clasp. Then there are min$\{l, m+j\}$ terms in the
expansion. But once we have a $U$ turn, we will show that it
becomes zero. If there is a $U$ turn we use a single clasp
expansion at the top-left clasp of weight $m+j$. Then all terms
vanish except one term which has $Y$ joining the top right clasp
of weight $(i+l,j+m)$ and the triangle in the top center. Then,
there is a sequence of $H$'s we can push that move the entire
shape by one string. Eventually $Y$ has to join two strings from
the triangle but we knew from lemma~\ref{onezerolem2} that it
becomes zero. Therefore, we can free $l$ strings from the right
middle clasp of weight $(l,j+m)$.
\end{proof}

Unfortunately this lemma~\ref{onezerolem1} is not true if $k\neq
0$. Actually only two terms survive but there ia a layer of $H$'s
which makes the problem difficult in this approach.
Continuing the prrof the theorem, lemma~\ref{onezerolem1} implies that
$M_{\Theta}(0,i+j,l,j+m,i+m,j+l)$ is
$$\frac{[j+l+1][i+j+l+m+2]}{[j+1][i+j+m+2]}
M_{\Theta}(0,i+j,0,j+m,i+m,j).$$ Then the result follows by
proposition~\ref{twozeroprop2}.

\subsection{Proof of Theorem ~\ref{onezerothm2}}

We use the same idea of lemma~\ref{twozerolem1} but for $j=0$, all
terms in this expansion do not vanish. For next step we need to
show the following lemma.

\begin{lem}
Let $\alpha= 0$ if $n> \mathrm{min}\{ l, k\}$,
$(\frac{[k]^n[l+k-n]!}{[l+k]!})^2$ if $n\le \mathrm{min}\{ l,
k\}$. Then
\begin{eqnarray}
\pspicture[.48](-3.1,-3.4)(3.9,3.3)
\rput(3,3){$i$}\rput(-1.5,1){$k$}
\rput(3,1.9){$m$}\rput(3,2.4){$l$}\rput(.5,.9){$n$}\rput(2.5,0){$m-n$}
\pccurve[angleA=180,angleB=90,ncurv=1](3.5,1.7)(1.5,1.2)\middlearrow
\pccurve[angleA=-90,angleB=180,ncurv=1](1.5,-1.2)(3.5,-1.7)\middlearrow
\pccurve[angleA=90,angleB=90,ncurv=1](-.7,1.2)(-2.3,.1)\middlearrow
\pccurve[angleA=-90,angleB=-90,ncurv=1](-2.3,-.1)(-.7,-1.2)\middlearrow
\pccurve[angleA=180,angleB=90,ncurv=1](-.5,2.7)(-2.7,.1)\middlearrow
\pccurve[angleA=-90,angleB=-180,ncurv=1](-2.7,-.1)(-.5,-2.7)\middlearrow
\pccurve[angleA=90,angleB=180,ncurv=1](-.3,1.2)(.5,2.2)\middlearrow
\pccurve[angleA=-180,angleB=-90,ncurv=1](.5,-2.2)(-.3,-1.2)\middlearrow
\pccurve[angleA=-135,angleB=-45,ncurv=1](1.3,1)(-.3,1)\middlearrow
\pccurve[angleA=45,angleB=135,ncurv=1](-.3,-1)(1.3,-1)\middlearrow
\psline(.5,2.2)(3.5,2.2) \psline(.5,-2.2)(3.5,-2.2)
\psline(-.5,2.7)(3.5,2.7) \psline(-.5,-2.7)(3.5,-2.7)
\psline(-.7,-1)(-.7,1)\psline[arrowscale=1.5]{->}(-.7,-.1)(-.7,.1)
\psline(1.7,-1)(1.7,1)\psline[arrowscale=1.5]{->}(1.7,.1)(1.7,-.1)
\psframe[linecolor=darkred](3.5,1.4)(3.7,2.9)
\psframe[linecolor=darkred](3.5,-2.9)(3.7,-1.4)
\psframe[linecolor=darkred](-3,-.1)(-2,.1)
\psframe[linecolor=darkred](-1,1)(0,1.2)
\psframe[linecolor=darkred](-1,-1)(0,-1.2)
\psframe[linecolor=darkred](1,1)(2,1.2)
\psframe[linecolor=darkred](1,-1)(2,-1.2)
\endpspicture
= \alpha\pspicture[.48](-3.5,-2.5)(2,2.5)
\rput(.5,2.3){$i+n$}\rput(-1.5,.3){$k-n$} \rput(1,.7){$m$}
\pccurve[angleA=180,angleB=90,ncurv=1](1.5,1)(-.3,.1)\middlearrow
\pccurve[angleA=-90,angleB=180,ncurv=1](-.3,-.1)(1.5,-1)\middlearrow
\pccurve[angleA=90,angleB=90,ncurv=1](-.7,.1)(-2.3,.1)\middlearrow
\pccurve[angleA=-90,angleB=-90,ncurv=1](-2.3,-.1)(-.7,-.1)\middlearrow
\pccurve[angleA=180,angleB=90,ncurv=1](-.5,2)(-2.7,.1)\middlearrow
\pccurve[angleA=-90,angleB=-180,ncurv=1](-2.7,-.1)(-.5,-2)\middlearrow
\psline(-.5,2)(1.5,2) \psline(-.5,-2)(1.5,-2)
\psframe[linecolor=darkred](1.5,.7)(1.7,2.3)
\psframe[linecolor=darkred](1.5,-2.3)(1.7,-.7)
\psframe[linecolor=darkred](-3,-.1)(-2,.1)
\psframe[linecolor=darkred](-1,-.1)(0,.1)
\endpspicture
\label{onezerolem2-1}
\end{eqnarray}
\label{onezerolem2}
\end{lem}

\begin{proof}
The clasp of weight $m\lambda_2$ can be pushed into the clasp of
weight $l\lambda_1 +$ $(i+m)\lambda_2$. For the clasp of weight
$(k+n)\lambda_1$, we use a single clasp expansion.
\end{proof}

First we use the equation~(\ref{a2abexp1}) at the middle clasp of
weight $(k+l)\lambda_1+m\lambda_2$. By lemma~\ref{onezerolem2} we
can transform each web to the web in the righthand side of
lemma~\ref{onezerolem2}.

\begin{eqnarray}
\pspicture[.5](-3.3,-2.5)(1.5,2.5)
\rput(.9,2.3){$i$}\rput(-1.5,.3){$k$} \rput(.9,1.3){$m$}
\rput(.9,1.8){$l$}
\pccurve[angleA=180,angleB=90,ncurv=1](1.3,1)(-.3,.1)\middlearrow
\pccurve[angleA=-90,angleB=180,ncurv=1](-.3,-.1)(1.3,-1)\middlearrow
\pccurve[angleA=90,angleB=180,ncurv=1](-.5,.1)(1.3,1.5)\middlearrow
\pccurve[angleA=180,angleB=-90,ncurv=1](1.3,-1.5)(-.5,-.1)\middlearrow
\pccurve[angleA=90,angleB=90,ncurv=1](-.7,.1)(-2.3,.1)\middlearrow
\pccurve[angleA=-90,angleB=-90,ncurv=1](-2.3,-.1)(-.7,-.1)\middlearrow
\pccurve[angleA=180,angleB=90,ncurv=1](-.5,2)(-2.7,.1)\middlearrow
\pccurve[angleA=-90,angleB=-180,ncurv=1](-2.7,-.1)(-.5,-2)\middlearrow
\psline(-.5,2)(1.3,2) \psline(-.5,-2)(1.3,-2)
\psframe[linecolor=darkred](1.3,.7)(1.5,2.3)
\psframe[linecolor=darkred](1.3,-2.3)(1.5,-.7)
\psframe[linecolor=darkred](-3,-.1)(-2,.1)
\psframe[linecolor=darkred](-1,-.1)(0,.1)
\endpspicture
=\sum_{n=0}^{\mathrm{min}\{ l+k,m\}} a_n
\pspicture[.5](-3.2,-3.4)(3.7,3.3)
\rput(3,3){$i$}\rput(-1.5,1){$k$}
\rput(3,1.9){$m$}\rput(3,2.4){$l$}\rput(.5,.9){$n$}\rput(2.5,0){$m-n$}
\pccurve[angleA=180,angleB=90,ncurv=1](3.5,1.7)(1.5,1.2)\middlearrow
\pccurve[angleA=-90,angleB=180,ncurv=1](1.5,-1.2)(3.5,-1.7)\middlearrow
\pccurve[angleA=90,angleB=90,ncurv=1](-.7,1.2)(-2.3,.1)\middlearrow
\pccurve[angleA=-90,angleB=-90,ncurv=1](-2.3,-.1)(-.7,-1.2)\middlearrow
\pccurve[angleA=180,angleB=90,ncurv=1](-.5,2.7)(-2.7,.1)\middlearrow
\pccurve[angleA=-90,angleB=-180,ncurv=1](-2.7,-.1)(-.5,-2.7)\middlearrow
\pccurve[angleA=90,angleB=180,ncurv=1](-.3,1.2)(.5,2.2)\middlearrow
\pccurve[angleA=-180,angleB=-90,ncurv=1](.5,-2.2)(-.3,-1.2)\middlearrow
\pccurve[angleA=-135,angleB=-45,ncurv=1](1.3,1)(-.3,1)\middlearrow
\pccurve[angleA=45,angleB=135,ncurv=1](-.3,-1)(1.3,-1)\middlearrow
\psline(.5,2.2)(3.5,2.2) \psline(.5,-2.2)(3.5,-2.2)
\psline(-.5,2.7)(3.5,2.7) \psline(-.5,-2.7)(3.5,-2.7)
\psline(-.7,-1)(-.7,1)\psline[arrowscale=1.5]{->}(-.7,-.1)(-.7,.1)
\psline(1.7,-1)(1.7,1)\psline[arrowscale=1.5]{->}(1.7,.1)(1.7,-.1)
\psframe[linecolor=darkred](3.5,1.4)(3.7,2.9)
\psframe[linecolor=darkred](3.5,-2.9)(3.7,-1.4)
\psframe[linecolor=darkred](-3,-.1)(-2,.1)
\psframe[linecolor=darkred](-1,1)(0,1.2)
\psframe[linecolor=darkred](-1,-1)(0,-1.2)
\psframe[linecolor=darkred](1,1)(2,1.2)
\psframe[linecolor=darkred](1,-1)(2,-1.2)
\endpspicture
\label{onezerothm2-1}
\end{eqnarray}

Since its shape can be written as $[k-n,i+n,m]$, by
proposition~\ref{twozerothm1} it has value

$$
\alpha \frac{[i+k+m+2]![k-n+1]![i+m]![m+1]!}{[i+k]![i+m+n]![k+m-n+1]![2]}.$$
Therefore, it completes the proof.

\newpage
\pagestyle{myheadings}
\markright{  \rm \normalsize CHAPTER 4. \hspace{0.5cm}
A Set of Complete Relations of $\mathcal{U}_q(\mathfrak{sl}(4,\mathbb{C}))$}
\chapter{A Complete Set of Relations of $\mathcal{U}_q(\mathfrak{sl}(4,\mathbb{C}))$}
\thispagestyle{myheadings}

Our webs are generated by the two shapes of trivalent vertices.

$$
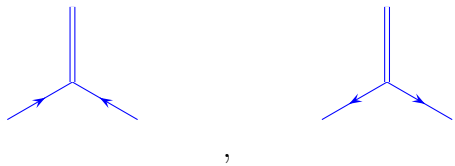
\pspicture(-1,-1)(1,1) \pnode(1;90){a1} \pnode(1;210){a2}
\pnode(1;330){a3} \pnode(0;0){b1} \ncline{a2}{b1}\middlearrow
\ncline{a3}{b1}\middlearrow \ncline[doubleline=true]{a1}{b1}
\endpspicture
\hskip 1cm , \hskip 1cm \pspicture(-1,-1)(1,1) \pnode(1;90){a1}
\pnode(1;210){a2} \pnode(1;330){a3} \pnode(0;0){b1}
\ncline{b1}{a2}\middlearrow \ncline{b1}{a3}\middlearrow
\ncline[doubleline=true]{a1}{b1}
\endpspicture
$$

And the following is our conjecture for a complete set of relations for
$\mathcal{U}_q(\mathfrak{sl}(4,\mathbb{C}))$.

\begin{eqnarray}
\pspicture[.4](-.6,-.5)(.6,.5)
\pscircle(0,0){.4}\psline[arrows=->,arrowscale=1.5](.1,.4)(.11,.4)
\endpspicture
&=& \left[ \begin{matrix} 4\\
1 \end{matrix} \right] \label{4single}\\ \doubleloop
&=& \left[ \begin{matrix} 4\\
2 \end{matrix} \right] \label{4double}\\
\pspicture[.4](-1.5,-.8)(1.5,.8) \pnode(1.2;180){a1}
\pnode(.4;180){a2} \pnode(.4;0){a3} \pnode(1.2;0){a4}
\nccurve[angleA=120,angleB=60,nodesep=1pt]{a2}{a3}\middlearrow
\nccurve[angleA=-120,angleB=-60,nodesep=1pt]{a2}{a3}\middlearrow
\ncline[doubleline=true]{a3}{a4} \ncline[doubleline=true]{a1}{a2}
\endpspicture
& =&  [2] \pspicture[.4](-.8,-.8)(.8,.8) \pnode(.6;180){a1}
\pnode(.6;0){a2}  \ncline[doubleline=true]{a1}{a2}
\endpspicture \label{4bigon1}\\
\pspicture[.4](-1.5,-.8)(1.5,.8) \pnode(1.2;180){a1}
\pnode(.4;180){a2} \pnode(.4;0){a3} \pnode(1.2;0){a4}
\ncline{a1}{a2}\middlearrow
\nccurve[angleA=-120,angleB=-60,nodesep=1pt]{a3}{a2}\middlearrow
\ncline{a3}{a4}\middlearrow
\nccurve[doubleline=true,angleA=120,angleB=60]{a2}{a3}
\endpspicture
& = &  [3]\pspicture[.4](-.8,-.8)(.8,.8) \pnode(.6;180){a1}
\pnode(.6;0){a2} \ncline{a1}{a2}\middlearrow
\endpspicture \label{4bigon2}\\
\pspicture[.4](-.8,-1.2)(1.4,1.2) \rput(-.6,1.0){\rnode{c1}{$$}}
\rput(-.6,-1.0){\rnode{c2}{$$}}\rput(.6,-1.0){\rnode{c3}{$$}}
\rput(.6,1.0){\rnode{c4}{$$}}\pnode(.4;90){b1} \pnode(.4;270){b2}
\ncline{c1}{b1}\middlearrow \ncline{c2}{b2}\middlearrow
\ncline{c3}{b2}\middlearrow \ncline{c4}{b1}\middlearrow
\ncline[doubleline=true]{b2}{b1}
\endpspicture
& = & \pspicture[.4](-1.2,-.8)(1.2,.8)
\rput(-1.0,.7){\rnode{c1}{$$}}
\rput(-1.0,-.7){\rnode{c2}{$$}}\rput(1.0,-.7){\rnode{c3}{$$}}
\rput(1.0,.7){\rnode{c4}{$$}}\pnode(.4;180){b1} \pnode(.4;0){b2}
\ncline{c1}{b1}\middlearrow \ncline{c2}{b1}\middlearrow
\ncline{c3}{b2}\middlearrow \ncline{c4}{b2}\middlearrow
\ncline[doubleline=true]{b1}{b2}
\endpspicture \label{4fourid1}\\
\pspicture[.4](-.8,-1.2)(1.4,1.2) \rput(-.6,1.0){\rnode{c1}{$$}}
\rput(-.6,-1.0){\rnode{c2}{$$}}\rput(.6,-1.0){\rnode{c3}{$$}}
\rput(.6,1.0){\rnode{c4}{$$}}\pnode(.4;90){b1} \pnode(.4;270){b2}
\ncline{b1}{c1}\middlearrow \ncline{b2}{c2}\middlearrow
\ncline{b2}{c3}\middlearrow \ncline{b1}{c4}\middlearrow
\ncline[doubleline=true]{b2}{b1}
\endpspicture
&= &\pspicture[.4](-1.2,-.8)(1.2,.8)
\rput(-1.0,.7){\rnode{c1}{$$}}
\rput(-1.0,-.7){\rnode{c2}{$$}}\rput(1.0,-.7){\rnode{c3}{$$}}
\rput(1.0,.7){\rnode{c4}{$$}}\pnode(.4;180){b1} \pnode(.4;0){b2}
\ncline{b1}{c1}\middlearrow \ncline{b1}{c2}\middlearrow
\ncline{b2}{c3}\middlearrow \ncline{b2}{c4}\middlearrow
\ncline[doubleline=true]{b1}{b2}
\endpspicture \label{4fourid2}\\
\pspicture[.4](-1.1,-1.1)(1.1,1.1)
\pnode(1;45){a1}\pnode(1;135){a2}\pnode(1;225){a3}\pnode(1;315){a4}
\pnode(.5;45){b1}\pnode(.5;135){b2}\pnode(.5;225){b3}\pnode(.5;315){b4}
\ncline{a1}{b1}\middlearrow \ncline{b2}{a2}\middlearrow
\ncline{a3}{b3}\middlearrow \ncline{b4}{a4}\middlearrow
\ncline[doubleline=true]{b1}{b2} \ncline{b2}{b3}\middlearrow
\ncline[doubleline=true]{b3}{b4} \ncline{b4}{b1}\middlearrow
\endpspicture
&= & [2] \pspicture[.4](-1.1,-1.1)(1.1,1.1)
\pnode(1;45){a1}\pnode(1;135){a2}\pnode(1;225){a3}\pnode(1;315){a4}
\nccurve[angleA=225,angleB=315]{a1}{a2}\middlearrow
\nccurve[angleA=45,angleB=135]{a3}{a4}\middlearrow
\endpspicture + \pspicture[.45](-1.1,-1.1)(1.1,1.1)
\pnode(1;45){a1}\pnode(1;135){a2}\pnode(1;225){a3}\pnode(1;315){a4}
\nccurve[angleA=225,angleB=135]{a1}{a4}\middlearrow
\nccurve[angleA=45,angleB=315]{a3}{a2}\middlearrow
\endpspicture
\label{4fourexp1}\\
\pspicture[.4](-1.1,-1.1)(1.1,1.1)
\pnode(1;45){a1}\pnode(1;135){a2}\pnode(1;225){a3}\pnode(1;315){a4}
\pnode(.5;45){b1}\pnode(.5;135){b2}\pnode(.5;225){b3}\pnode(.5;315){b4}
\ncline{b2}{a2}\middlearrow \ncline{a3}{b3}\middlearrow
\ncline{b2}{b1}\middlearrow \ncline{b4}{b3}\middlearrow
\ncline{b4}{b1}\middlearrow \ncline[doubleline=true]{a1}{b1}
\ncline[doubleline=true]{b4}{a4} \ncline[doubleline=true]{b3}{b2}
\endpspicture
&= &\pspicture[.4](-1.1,-1.1)(1.1,1.1)
\pnode(1;45){a1}\pnode(1;135){a2}\pnode(1;225){a3}\pnode(1;315){a4}
\pnode(.3;90){b1}\pnode(.3;270){b2} \ncline{b1}{a2}\middlearrow
\ncline{b1}{b2}\middlearrow \ncline{a3}{b2}\middlearrow
\ncline[doubleline=true]{a1}{b1} \ncline[doubleline=true]{a4}{b2}
\endpspicture + \pspicture[.45](-1.1,-1.1)(1.1,1.1)
\pnode(1;45){a1}\pnode(1;135){a2}\pnode(1;225){a3}\pnode(1;315){a4}
\pnode(.3;90){b1}\pnode(.3;270){b2}
\nccurve[angleA=45,angleB=315]{a3}{a2}\middlearrow
\nccurve[doubleline=true,angleA=225,angleB=135]{a1}{a4}
\endpspicture
\label{4fourexp2}\\
\pspicture[.4](-1.1,-1.1)(1.1,1.1)
\pnode(1;45){a1}\pnode(1;135){a2}\pnode(1;225){a3}\pnode(1;315){a4}
\pnode(.5;45){b1}\pnode(.5;135){b2}\pnode(.5;225){b3}\pnode(.5;315){b4}
\ncline{b2}{b1}\middlearrow \ncline{b4}{b1}\middlearrow
\ncline{b2}{b3}\middlearrow \ncline{b4}{b3}\middlearrow
\ncline[doubleline=true]{a1}{b1} \ncline[doubleline=true]{a2}{b2}
\ncline[doubleline=true]{a3}{b3} \ncline[doubleline=true]{a4}{b4}
\endpspicture
& = &\pspicture[.4](-1.1,-1.1)(1.1,1.1)
\pnode(1;45){a1}\pnode(1;135){a2}\pnode(1;225){a3}\pnode(1;315){a4}
\pnode(.5;45){b1}\pnode(.5;135){b2}\pnode(.5;225){b3}\pnode(.5;315){b4}
\ncline{b1}{b2}\middlearrow \ncline{b1}{b4}\middlearrow
\ncline{b3}{b2}\middlearrow \ncline{b3}{b4}\middlearrow
\ncline[doubleline=true]{a1}{b1} \ncline[doubleline=true]{a2}{b2}
\ncline[doubleline=true]{a3}{b3} \ncline[doubleline=true]{a4}{b4}
\endpspicture  \label{4fourexp3}\\
\pspicture[.4](-1,-1)(1,1) \pnode(.5;30){a1} \pnode(.5; 90){a2}
\pnode(.5;150){a3}\pnode(.5;210){a4}
\pnode(.5;270){a5}\pnode(.5;330){a6} \rput(1; 30){\rnode{b1}{$$}}
\rput(1; 90){\rnode{b2}{$$}} \rput(1;150){\rnode{b3}{$$}}
\rput(1;210){\rnode{b4}{$$}} \rput(1;270){\rnode{b5}{$$}}
\rput(1;330){\rnode{b6}{$$}} \ncline{a2}{a1}\middlearrow
\ncline{a4}{a3}\middlearrow \ncline{a6}{a5}\middlearrow
\ncline{b1}{a1}\middlearrow \ncline{a2}{b2}\middlearrow
\ncline{b3}{a3}\middlearrow \ncline{a4}{b4}\middlearrow
\ncline{b5}{a5}\middlearrow \ncline{a6}{b6}\middlearrow
\ncline[doubleline=true]{a4}{a5} \ncline[doubleline=true]{a2}{a3}
\ncline[doubleline=true]{a6}{a1}
\endpspicture
& = &\pspicture[.4](-1,-1.1)(1,1.1) \pnode(.5; 30){a1}\pnode(.5;
90){a2} \pnode(.5;150){a3}\pnode(.5;210){a4}
\pnode(.5;270){a5}\pnode(.5;330){a6} \rput(1; 30){\rnode{b1}{$$}}
\rput(1; 90){\rnode{b2}{$$}} \rput(1;150){\rnode{b3}{$$}}
\rput(1;210){\rnode{b4}{$$}} \rput(1;270){\rnode{b5}{$$}}
\rput(1;330){\rnode{b6}{$$}} \ncline{a2}{a3}\middlearrow
\ncline{a4}{a5}\middlearrow \ncline{a6}{a1}\middlearrow
\ncline{b1}{a1}\middlearrow \ncline{a2}{b2}\middlearrow
\ncline{b3}{a3}\middlearrow \ncline{a4}{b4}\middlearrow
\ncline{b5}{a5}\middlearrow \ncline{a6}{b6}\middlearrow
\ncline[doubleline=true]{a1}{a2} \ncline[doubleline=true]{a3}{a4}
\ncline[doubleline=true]{a5}{a6}
\endpspicture
-  \pspicture[.4](-1.3,-1.1)(1.3,1.1)  \rput(1;
30){\rnode{b1}{$$}} \rput(1; 90){\rnode{b2}{$$}}
\rput(1;150){\rnode{b3}{$$}} \rput(1;210){\rnode{b4}{$$}}
\rput(1;270){\rnode{b5}{$$}} \rput(1;330){\rnode{b6}{$$}}
\nccurve[angleA=210,angleB=270]{b1}{b2}\middlearrow
\nccurve[angleA=330,angleB=30]{b3}{b4}\middlearrow
\nccurve[angleA=90,angleB=150]{b5}{b6}\middlearrow
\endpspicture
+ \pspicture[.4](-1.3,-1.1)(1.3,1.1) \rput(1;
30){\rnode{b1}{$$}}\rput(1; 90){\rnode{b2}{$$}}
\rput(1;150){\rnode{b3}{$$}}\rput(1;210){\rnode{b4}{$$}}
\rput(1;270){\rnode{b5}{$$}}\rput(1;330){\rnode{b6}{$$}}
\nccurve[angleA=330,nodesepA=3pt,angleB=270]{b3}{b2}\middlearrow
\nccurve[angleA=90,angleB=30,nodesepB=3pt]{b5}{b4}\middlearrow
\nccurve[angleA=210,angleB=150]{b1}{b6}\middlearrow
\endpspicture\label{4sixexp}
\end{eqnarray}

 First of all, we compute the dimension of
invariant space of all tensor products of $4, 6$ fundamental
representation of $\mathcal{U}_q(\mathfrak{sl}(4,\mathbb{C}))$.
There is a general way to find a basis webs with a fixed boundary,
all of them are fundamental representation, for
$\mathcal{U}_q(\mathfrak{sl}(2,\mathbb{C}))$ and
$\mathcal{U}_q(\mathfrak{sl}(3,\mathbb{C}))$~\cite{KK:notdual}. It
is still unknown how we can actually find all basis webs of all
possible boundaries of fundamental representations. But for our
case, there are only few boundaries and we can find the dimension
and even find a basis webs without difficulty, most of basis webs do
not have any faces.

Since all webs are generated by two trivalent vertices, by
multiplying a complex number, we can have a different set of
generators. Therefore, we have two choices of freedom to set any
two independent coefficients. Let $a, b, c, d, e$ and $f$ be
unknowns in equation~\ref{4bigon2}, \ref{4fourid2},
\ref{4fourexp1} and \ref{4fourexp2}. By the quantum Weyl formula,
we do know the value of the first two equations. We use the first
choice of freedom to have the equation~\ref{4bigon1}. The
following equality implies $a=[3]$.

\begin{align*}
\pspicture[.4](-1.3,-1.1)(1.3,1.1)
\rput(1;0){\rnode{b1}{$$}}\rput(1; 180){\rnode{b2}{$$}}
\nccurve[doubleline=true]{b1}{b2}
\nccurve[angleA=110,angleB=70]{b1}{b2}\middlearrow
\nccurve[angleA=250,angleB=300]{b1}{b2}\middlearrow
\endpspicture &=a
\pspicture[.4](-.6,-.5)(.6,.5)
\pscircle(0,0){.4}\psline[arrows=->,arrowscale=1.5](.1,.4)(-.1,.4)
\endpspicture =[4]a\\
&=[2]\doubleloop=[2]\frac{[4][3]}{[2]}=[4][3]
\end{align*}

We use the last choice of freedom to have equation~\ref{4fourid1}.
To get the equation~\ref{4fourid2} (which is actually the dual of
the equation~\ref{4fourid1}), we start from the following
equations.

$$
\pspicture[.4](-1.1,-1.1)(1.1,1.1)
\pnode(1;45){a1}\pnode(1;135){a2}\pnode(1;225){a3}\pnode(1;315){a4}
\pnode(.5;45){b1}\pnode(.5;135){b2}\pnode(.5;225){b3}\pnode(.5;315){b4}
\ncline{b1}{a1}\middlearrow \ncline{a2}{b2}\middlearrow
\ncline{a3}{b3}\middlearrow \ncline{b4}{a4}\middlearrow
\ncline[doubleline=true]{b1}{b4} \ncline{b1}{b2}\middlearrow
\ncline[doubleline=true]{b3}{b2} \ncline{b4}{b3}\middlearrow
\endpspicture
= \alpha \pspicture[.4](-1.1,-1.1)(1.1,1.1)
\pnode(1;45){a1}\pnode(1;135){a2}\pnode(1;225){a3}\pnode(1;315){a4}
\nccurve[angleA=-45,angleB=-135]{a2}{a1}\middlearrow
\nccurve[angleA=45,angleB=135]{a3}{a4}\middlearrow
\endpspicture + \beta \pspicture[.45](-1.1,-1.1)(1.1,1.1)
\pnode(1;45){a1}\pnode(1;135){a2}\pnode(1;225){a3}\pnode(1;315){a4}
\pnode(.4;0){b1}\pnode(.4;180){b2} \ncline{b1}{a1}\middlearrow
\ncline{a2}{b2}\middlearrow \ncline{b1}{a4}\middlearrow
\ncline{a3}{b2}\middlearrow \ncline[doubleline=true]{b1}{b2}
\endpspicture
\label{4fourexp5}
$$

Using~\ref{4fourid1} and~\ref{4fourid2} (with an unknown variable
$a$), we found $\alpha =0$ and $\beta = [2]b$. By attaching a
$U$ turn on the top of each webs in the equation, we get
$[2][3]=[4]\alpha +[3]\beta$. Thus $b=1$. For
equation~\ref{4fourexp1}, we attach $U$ turns on the top and right
side of each webs in the equation. Then the resulting web can be
expanded as a linear combination of a basis webs of different
boundary. By comparing the coefficients, we get $[3][3]=[4]c+d,
[2][3]=c+[4]d$. It is easy to find that $c=[2], d=1$. For
equation~\ref{4fourexp2} we attach

\begin{eqnarray}
\pspicture[.4](-1.1,-1.1)(1.1,1.1)
\pnode(1;45){a1}\pnode(1;135){a2}\pnode(1;225){a3}\pnode(1;315){a4}
\pnode(.3;90){b1}\pnode(.3;270){b2} \ncline{b1}{a1}\middlearrow
\ncline{b1}{b2}\middlearrow \ncline{a4}{b2}\middlearrow
\ncline[doubleline=true]{a2}{b1} \ncline[doubleline=true]{a3}{b2}
\endpspicture
\label{4extra1}
\end{eqnarray}

to right side of each webs to get $e=f=1$.

For last two equations~\ref{4fourexp3} and~\ref{4sixexp}, we need
to start from the following equations.

\begin{eqnarray}
\pspicture[.4](-1.1,-1.1)(1.1,1.1)
\pnode(1;45){a1}\pnode(1;135){a2}\pnode(1;225){a3}\pnode(1;315){a4}
\pnode(.5;45){b1}\pnode(.5;135){b2}\pnode(.5;225){b3}\pnode(.5;315){b4}
\ncline{b2}{b1}\middlearrow \ncline{b4}{b1}\middlearrow
\ncline{b2}{b3}\middlearrow \ncline{b4}{b3}\middlearrow
\ncline[doubleline=true]{a1}{b1} \ncline[doubleline=true]{a2}{b2}
\ncline[doubleline=true]{a3}{b3} \ncline[doubleline=true]{a4}{b4}
\endpspicture
& = & g\pspicture[.4](-1.1,-1.1)(1.1,1.1)
\pnode(1;45){a1}\pnode(1;135){a2}\pnode(1;225){a3}\pnode(1;315){a4}
\pnode(.5;45){b1}\pnode(.5;135){b2}\pnode(.5;225){b3}\pnode(.5;315){b4}
\ncline{b1}{b2}\middlearrow \ncline{b1}{b4}\middlearrow
\ncline{b3}{b2}\middlearrow \ncline{b3}{b4}\middlearrow
\ncline[doubleline=true]{a1}{b1} \ncline[doubleline=true]{a2}{b2}
\ncline[doubleline=true]{a3}{b3} \ncline[doubleline=true]{a4}{b4}
\endpspicture  + h \pspicture[.4](-1.1,-1.1)(1.1,1.1)
\pnode(1;45){a1}\pnode(1;135){a2}\pnode(1;225){a3}\pnode(1;315){a4}
\pnode(.5;45){b1}\pnode(.5;135){b2}\pnode(.5;225){b3}\pnode(.5;315){b4}
\nccurve[doubleline=true,angleA=-135,angleB=135]{a1}{a4}
\nccurve[doubleline=true,angleA=-45,angleB=45]{a2}{a3}
\endpspicture + i \pspicture[.4](-1.1,-1.1)(1.1,1.1)
\pnode(1;45){a1}\pnode(1;135){a2}\pnode(1;225){a3}\pnode(1;315){a4}
\pnode(.5;45){b1}\pnode(.5;135){b2}\pnode(.5;225){b3}\pnode(.5;315){b4}
\nccurve[doubleline=true,angleA=-135,angleB=-45]{a1}{a2}
\nccurve[doubleline=true,angleA=45,angleB=135]{a3}{a4}
\endpspicture
\label{4fourexp32}\\ \nonumber \pspicture[.4](-1,-1)(1,1)
\pnode(.5;30){a1} \pnode(.5; 90){a2}
\pnode(.5;150){a3}\pnode(.5;210){a4}
\pnode(.5;270){a5}\pnode(.5;330){a6} \rput(1; 30){\rnode{b1}{$$}}
\rput(1; 90){\rnode{b2}{$$}} \rput(1;150){\rnode{b3}{$$}}
\rput(1;210){\rnode{b4}{$$}} \rput(1;270){\rnode{b5}{$$}}
\rput(1;330){\rnode{b6}{$$}} \ncline{a2}{a1}\middlearrow
\ncline{a4}{a3}\middlearrow \ncline{a6}{a5}\middlearrow
\ncline{b1}{a1}\middlearrow \ncline{a2}{b2}\middlearrow
\ncline{b3}{a3}\middlearrow \ncline{a4}{b4}\middlearrow
\ncline{b5}{a5}\middlearrow \ncline{a6}{b6}\middlearrow
\ncline[doubleline=true]{a4}{a5} \ncline[doubleline=true]{a2}{a3}
\ncline[doubleline=true]{a6}{a1}
\endpspicture
& = & j \pspicture[.4](-1,-1.1)(1,1.1) \pnode(.5;
30){a1}\pnode(.5; 90){a2} \pnode(.5;150){a3}\pnode(.5;210){a4}
\pnode(.5;270){a5}\pnode(.5;330){a6} \rput(1; 30){\rnode{b1}{$$}}
\rput(1; 90){\rnode{b2}{$$}} \rput(1;150){\rnode{b3}{$$}}
\rput(1;210){\rnode{b4}{$$}} \rput(1;270){\rnode{b5}{$$}}
\rput(1;330){\rnode{b6}{$$}} \ncline{a2}{a3}\middlearrow
\ncline{a4}{a5}\middlearrow \ncline{a6}{a1}\middlearrow
\ncline{b1}{a1}\middlearrow \ncline{a2}{b2}\middlearrow
\ncline{b3}{a3}\middlearrow \ncline{a4}{b4}\middlearrow
\ncline{b5}{a5}\middlearrow \ncline{a6}{b6}\middlearrow
\ncline[doubleline=true]{a1}{a2} \ncline[doubleline=true]{a3}{a4}
\ncline[doubleline=true]{a5}{a6}
\endpspicture
+ k  \pspicture[.4](-1.3,-1.1)(1.3,1.1)  \rput(1;
30){\rnode{b1}{$$}} \rput(1; 90){\rnode{b2}{$$}}
\rput(1;150){\rnode{b3}{$$}} \rput(1;210){\rnode{b4}{$$}}
\rput(1;270){\rnode{b5}{$$}} \rput(1;330){\rnode{b6}{$$}}
\nccurve[angleA=210,angleB=270]{b1}{b2}\middlearrow
\nccurve[angleA=330,angleB=30]{b3}{b4}\middlearrow
\nccurve[angleA=90,angleB=150]{b5}{b6}\middlearrow
\endpspicture
+l \pspicture[.4](-1.3,-1.1)(1.3,1.1) \rput(1;
30){\rnode{b1}{$$}}\rput(1; 90){\rnode{b2}{$$}}
\rput(1;150){\rnode{b3}{$$}}\rput(1;210){\rnode{b4}{$$}}
\rput(1;270){\rnode{b5}{$$}}\rput(1;330){\rnode{b6}{$$}}
\nccurve[angleA=330,nodesepA=3pt,angleB=270]{b3}{b2}\middlearrow
\nccurve[angleA=90,angleB=30,nodesepB=3pt]{b5}{b4}\middlearrow
\nccurve[angleA=210,angleB=150]{b1}{b6}\middlearrow
\endpspicture \\
&+& m \pspicture[.4](-1.3,-1.1)(1.3,1.1) \rput(1;
30){\rnode{b1}{$$}}\rput(1; 90){\rnode{b2}{$$}}
\rput(1;150){\rnode{b3}{$$}}\rput(1;210){\rnode{b4}{$$}}
\rput(1;270){\rnode{b5}{$$}}\rput(1;330){\rnode{b6}{$$}}
\ncline{b1}{b4}\middlearrow
\nccurve[angleA=-30,angleB=-90,nodesepB=3pt]{b3}{b2}\middlearrow
\nccurve[angleA=90,angleB=150]{b5}{b6}\middlearrow
\endpspicture
 +n \pspicture[.4](-1.3,-1.1)(1.3,1.1)
\rput(1; 30){\rnode{b1}{$$}}\rput(1; 90){\rnode{b2}{$$}}
\rput(1;150){\rnode{b3}{$$}}\rput(1;210){\rnode{b4}{$$}}
\rput(1;270){\rnode{b5}{$$}}\rput(1;330){\rnode{b6}{$$}}
\ncline{b5}{b2}\middlearrow
\nccurve[angleA=-30,angleB=30,nodesepB=3pt]{b3}{b4}\middlearrow
\nccurve[angleA=210,angleB=150]{b1}{b6}\middlearrow
\endpspicture
+p \pspicture[.4](-1.3,-1.1)(1.3,1.1) \rput(1;
30){\rnode{b1}{$$}}\rput(1; 90){\rnode{b2}{$$}}
\rput(1;150){\rnode{b3}{$$}}\rput(1;210){\rnode{b4}{$$}}
\rput(1;270){\rnode{b5}{$$}}\rput(1;330){\rnode{b6}{$$}}
\ncline{b3}{b6}\middlearrow
\nccurve[angleA=90,angleB=30,nodesepB=3pt]{b5}{b4}\middlearrow
\nccurve[angleA=210,angleB=-90]{b1}{b2}\middlearrow
\endpspicture
\label{4sixexp2}
\end{eqnarray}

By attaching $U$ turns and $H$ (as in~\ref{4extra1}) for
equation~\ref{4fourexp3}, we get

\begin{align*}
[2] [3] &=[2] [3]g + h +
\frac{[4][3]}{[2]}i\\
 [2] [3] &=[2] [3]g +\frac{[4][3]}{[2]} h+
i\\ [2] &=[2]g+[3]h\\ [2] &=[2]g+i
\end{align*}

One can solve them to have
$g=1$ and $h=i=0$. For the equation~\ref{4sixexp}, we just need to
attach $H$ (as in~\ref{4extra1}) to right top side of each basis
webs in the equation~\ref{4sixexp2}. Then we follow the same
procedure to get the following six equations: $j=1$, $[2]^2j$
$+[3]k$ $=1$, $[2]j$ $+[3]p$ $=[2]$, $l=1$, $m=0$ and $n=0$.

\bibliographystyle{hamsalpha}
\bibliography{dongseok}

\end{document}